%% file: DGart.tex
\documentclass{amsproc}

\usepackage{asa}
\usepackage{amscd}
\usepackage{amssymb}
\usepackage{array}
\usepackage{epsfig}

\input{DGdefs}

\begin{document}

\title{Basics of the b--Calculus}
 
\author{Daniel Grieser}
\address{Institut f\"ur Mathematik, Humboldt-Universit\"at zu Berlin, Unter den
Linden 6, 10099 Berlin, Germany}
\email{grieser@mathematik.hu-berlin.de}
\thanks{2000 {\em Mathematics Subject Classification. } Primary
58-01,    % Global Analysis, expository
58G15;    % PsiDO
Secondary
58G18,    % perturbations, asympt.
35B25,    % singular perturb.
35B40,    % asympt. beh. of solutions
35C20,    % asympt. exp.
35S05}    % PsiDOs

\thanks{The author was supported by the Deutsche Forschungsgemeinschaft.}      

\keywords{pseudodifferential analysis, manifolds with corners, 
blow-up, asymptotic analysis}
%\subjclass[2000]{Primary

%\date{June 21, 2000}

\begin{abstract}
R.\ B.\ Melrose's \Mb-calculus provides a framework for dealing with
problems of partial differential equations that arise in singular
or degenerate geometric situations. This article is a somewhat informal
short course introducing many of the basic ideas of this world,
assuming little more than a basic analysis and manifold background.
As examples, classical pseudodifferential operators on manifolds and
\Mb-pseudodifferential (also known as totally characteristic)
operators on manifolds with boundary are discussed.
\end{abstract}

\maketitle
%
%%%%%%%%%%%%%%%%%%%%%%%%%%%%%%%%%%%%%%%%%%%%%%%%%%%%%%%%%%%%
%%  Insert here the body of your article
%%%%%%%%%%%%%%%%%%%%%%%%%%%%%%%%%%%%%%%%%%%%%%%%%%%%%%%%%%%%
%
\setcounter{tocdepth}{3}  % put subsubsections in table of contents

\tableofcontents
\listoffigures

\input{DGintro}

\input{DGgeometry}

\input{DGanalysis}

\input{DGpde}

\input{DGappendix}

\bibliographystyle{abbrv}
\bibliography{lib}
%%%%%%%%%%%%%%%%%%%%%%%%%%%%%%%%%%%%%%%%%%%%%%%%%%%%%%%%%%%%

\end{document}

%% file: DGdefs.tex
\newenvironment{boldlist}
{\begin{list}{}
  {
  }
}
{\end{list}
}

\newenvironment{principle}{%\stepcounter{princctr}
                           \medskip\par
                           {\bf Principle% \arabic{princctr}
                           }\begin{quote} }
                          {\end{quote}  \par}

\newcommand{\Mb}{$b$}  % Melrose b 
\newcommand{\inthisbook} {in this book}    %% to be replaced by empty string
\newcommand{\ff}{{\rm ff}}
\newcommand{\rb}{{\rm rb}}
\newcommand{\lb}{{\rm lb}}
\newcommand{\fff}{{\rm fff}}
\newcommand{\bface} {{\rm b}}
\newcommand{\Rplus}{\R_+}
\newcommand{\utilde}{\tilde{u}}
\newcommand{\xtilde}{\tilde{x}}
\newcommand{\ytilde}{\tilde{y}}

\newcommand{\gtilde}{\tilde{g}}
\newcommand{\fbar}{\overline{f}}

\newcommand{\pitilde}{\tilde{\pi}}
\newcommand{\pitildebar}{\overline{\pitilde}}

\newcommand{\Dcal}{{\mathcal D}}
\newcommand{\Ecal}{{\mathcal E}}
\newcommand{\Fcal}{{\mathcal F}}
\newcommand{\extunion}{\overline{\cup}}
\newcommand{\Extunion}{\overline{\bigcup}}
\newcommand{\codim}{\mathop{\rm codim}}
\newcommand{\interior}[1]{{#1}^{\rm o}}
\newcommand{\Wint}{\interior{W}}
\newcommand{\Zint}{\interior{Z}}
\newcommand{\singsupp}{{\rm sing supp}\,}
\newcommand{\PDO}{$\Psi$DO}
\newcommand{\Diff}{{\rm Diff}}
\newcommand{\dxx}{\frac{dx}x}
\newcommand{\dss}{\frac{ds}s}
\newcommand{\dxshalf}{\left|\dxx\dss\right|^{1/2}}
\newcommand{\dxxprime}{\frac{dx'}{x'}}
\newcommand{\Spec}{{\rm Spec}}

\newcommand{\Mbook}{\cite{DGMel:APSIT}}
\newcommand{\Mnotes}{\cite{DGMel:DAMWC}}
\def\qed{\vbox{\hrule
  \hbox{\vrule\hbox to 5pt{\vbox to 8pt{\vfil}\hfil}\vrule}\hrule}}

\newcommand {\eps}{\varepsilon}
\newcommand\dist{{\rm dist\,}}

\newcommand\Id{{\rm Id}}
\newcommand\id{{\rm id}}

\renewcommand\Re{{\rm Re\,}}

\newcommand\supp{{\rm supp\,}}

\newcommand {\R}{\mathbb{R}}
\newcommand {\Z}{\mathbb{Z}}
\newcommand {\C}{\mathbb{C}}
\newcommand {\N}{\mathbb{N}}

%% file: DGintro.tex
\section{Introduction} \label{secintro}
This article gives a leisurely introduction to the \Mb-calculus
of R.\ B.\ Melrose. Here, we use the
term '\Mb-calculus'  in a broad
sense: A geometrically inspired
 way of viewing and solving problems about smooth functions and
distributions (especially their asymptotic and singular behavior)
and differential equations (especially as they arise in singular geometric
situations); a set of concepts introduced to realize
this view mathematically; and a set of basic and general theorems
about these concepts. The \Mb-calculus in the narrower, technical sense (as a set of operators)
will also be discussed.

The style of this article is rather informal. 
We emphasize examples, motivations and intuition and often refer to the
literature for full proofs and the most general definitions. While the
ultimate goal is to extend the classical pseudodifferential operator (\PDO)
calculus, large parts (Sections \ref{DGsecgeom} and \ref{DGsecana})
are interesting in other contexts as well.
Knowledge of the classical \PDO\ calculus is not a prerequisite. Rather,
it will be introduced, if sketchily, as  the simplest instance of the more
general theory to be developed.

We begin with some general considerations
on  solving  linear partial differential equations
(PDE), to show in which direction we aim. 
Since it is usually impossible to get an explicit solution, one wants
to study existence and uniqueness, and qualitative
properties of solutions. The PDE may contain parameters, then one wants to study
how these things depend on the parameters. (For example, spectral problems
are of this kind.)

To fix ideas, let us look at the case of an {\em elliptic} partial
differential operator $P$ with smooth coefficients on some manifold $X_0$
(for example, the Laplacian on a Riemannian manifold)
and at the equation 
\begin{equation}\label{DGeqPDE}
Pu=f.
\end{equation}
Solving for $u$ in terms of $f$ means finding an inverse $Q$ of $P$
(which we assume to exist for the moment, between suitably chosen
function spaces). 'Knowing' $Q$ would mean knowing its Schwartz
 kernel, that is the distribution, also denoted $Q$, on $X_0\times X_0$
satisfying $(Qf)(x)=\int_{X_0} Q(x,x')\,f(x')\,dx'$
 (also known
as Green's function).
Many important properties of equation \eqref{DGeqPDE} may be read off from
certain partial information on $Q$:

\begin{enumerate}
\item[A.] The location and nature of the singularities of $Q$ (i.e. places
where the distribution $Q$ is not a $C^\infty$ function).
\item[B.] The asymptotic 
behavior of $Q$ when approaching the 'boundary' of $X_0\times X_0$ (i.e.\ when
leaving any compact subset).
\end{enumerate}
We will refer to this information as the {\em singularity structure} of $Q$.

Then we reformulate problem \eqref{DGeqPDE}  as:
\begin{quote}
{\bf Main Problem: }{\em
                Given the singularity structure of $P$,
  determine the singularity structure of $Q=P^{-1}$.}
\end{quote}
(Here $P$ is also identified with its Schwartz kernel.)
If $P$ depends smoothly on parameters in a space $T$, then one wants to
find the singularity structure of $Q$ on the space $X_0\times X_0\times T$.%
\footnote{
In applications, if we start with an operator 
that has 'singular' coefficients (at some place or parameter value)   then
we take $X_0$ and $T$ to be the set  where the coefficients are smooth. 
Similarly, if we start with an operator on a 'singular space' (e.g.\ a 
manifold with boundary) then $X_0$ is the smooth part of that space
(the interior of the manifold with boundary).
}
\footnote{\label{DGfnspec}%
Note that many problems of linear analysis are specializations of
the Main Problem, for example: asymptotics of eigenvalues and eigenfunctions
under singular perturbations,  mapping and Fredholm
properties of elliptic operators, heat kernel asymptotics (the latter
in the analogous parabolic setup).
}
Since the coefficients of $P$ are smooth functions on $X_0\times T$, 
the singularity structure of $P$ only depends
on their asymptotic behavior 'near the boundary'  of $X_0\times T$
(i.e.\ when leaving any compact subset); for
example, they may blow up or degenerate ('non-uniform ellipticity')
in various ways.
If $P$ is not invertible then one asks the same questions for approximate
inverses (parametrices) of $P$.

For example, if $X_0$ is compact (and $P$ elliptic as before) and there are no
parameters  then the singularity structure
is given by  point A above only,  and
the classical pseudodifferential  calculus tells us that $Q$ is
smooth outside the diagonal and has 'conormal'  singularities
on the diagonal, and gives a recipe for calculating 
these modulo smooth functions (i.e.\ the complete symbol of $Q$).
As an example with parameter, consider $P_z=-\Delta+z$ on a compact Riemannian
manifold $X_0$, for $z\in[1,\infty)$. Then the singularity structure of
$P_z^{-1}$ describes not only the conormal singularity at the diagonal, but also
the asymptotic behavior of the resolvent kernel as $z\to\infty$.

The goal of the \Mb-calculus is to solve the Main Problem  for a fairly
broad class of singularity structures of $P$, the so-called boundary
fibration structures.
This general goal still seems out of reach, but a growing list of instances
shows the versatility of the \Mb-calculus in treating problems arising in
geometric analysis (see the references given below).

Figure \ref{DGfigbcalc} shows a rough outline of the \Mb-calculus approach to
the Main Problem. 
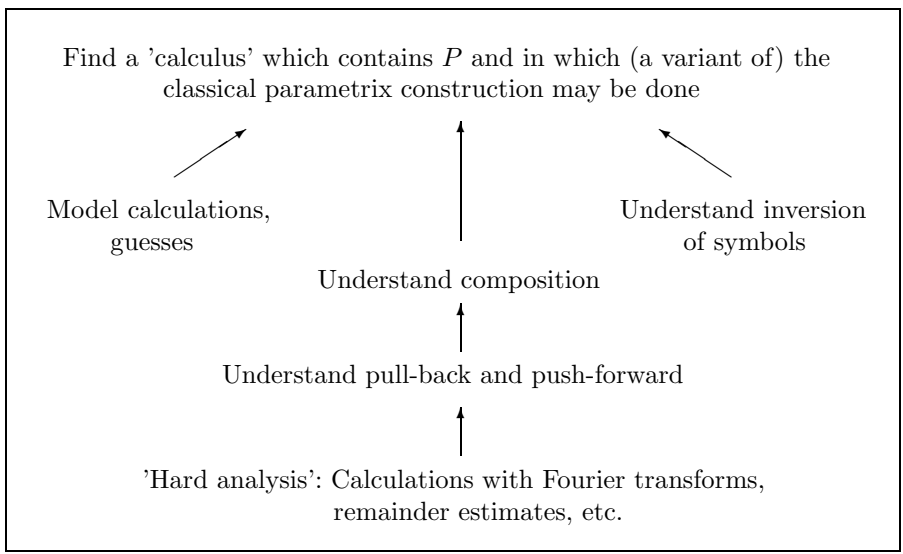
\begin{figure}[htbp]
\framebox{
\setlength{\unitlength}{3945sp}%
\begin{picture}(5400,3300)(2550,-4300)
\thinlines
\put(7000,-2000){\vector(-3, 2){450}}
\put(3500,-2000){\vector( 3, 2){450}}
\put(5300,-2400){\vector( 0, 1){750}}
\put(5300,-3100){\vector( 0, 1){300}}
\put(5300,-3750){\vector( 0, 1){300}}
%%%
\put(2800,-1300){\makebox(0,0)[lb]{\smash{Find
a 'calculus' which contains $P$ and in which (a variant of) the}}}
\put(3430,-1500){\makebox(0,0)[lb]{\smash{classical
parametrix construction may be done}}}
\put(6300,-2260){\makebox(0,0)[lb]{\smash{Understand inversion}}}
\put(6700,-2460){\makebox(0,0)[lb]{\smash{of symbols}}}
\put(2700,-2260){\makebox(0,0)[lb]{\smash{Model calculations,}}}
\put(3100,-2460){\makebox(0,0)[lb]{\smash{guesses}}}
\put(4400,-2700){\makebox(0,0)[lb]{\smash{Understand composition}}}
\put(3800,-3300){\makebox(0,0)[lb]{\smash{Understand pull-back and push-forward}}}    
\put(3300,-3960){\makebox(0,0)[lb]{\smash{'Hard analysis': Calculations with Fourier transforms,}}} 
\put(4500,-4160){\makebox(0,0)[lb]{\smash{remainder estimates, etc.}}}
\end{picture}
}
%%%
\caption{The \Mb-calculus hierarchy for inverting $P$}\label{DGfigbcalc}
\end{figure}
An arrow means 'is used for'.
A {\em calculus} is a set of operators with a fixed singularity structure, which
is closed (at least conditionally) under composition, together with a collection
of {\em symbols}, i.e.\ rules that assign to each operator certain 'simplified' operators,
usually by some sort of (partial) freezing of coefficients.
'Understanding' composition etc.\ means determining the singularity structure of the
composition from the singularity structures of the factors.
In particular, understanding inversion of symbols is another instance of the
Main Problem itself, but for a simpler class of operators, and this shows the
iterative nature of the problem. The lowest level of the iteration is inversion of
constant coefficient operators, which may be done directly using the
Fourier transform.

Therefore, in the construction of a calculus one expects the Fourier transform
to play a central role, and this is reflected in the definition and occurrence of conormal
distributions. However, calculations involving the Fourier transform tend to
be messy and to obscure essential structures; this is why it is banished to
the bottom level in
Figure \ref{DGfigbcalc}: The Fourier transform  is
 only used explicitly in analyzing pull-back and push-forward
of conormal distributions; composition is then  reduced to 
a combination of pull-back and push-forward operations.%
\footnote{
This should be taken with a grain of salt, but gives a general
guideline.  In the 'full \Mb-calculus' in Section \ref{DGsecpde} the
Fourier transform, in the guise of the Mellin transform,
is also used for inverting the 'indicial operator'; this belongs to
the 'Understand inversion of symbols' part in Figure \ref{DGfigbcalc}.
}
%\footnote{
%More precisely, it is used in the inversion of symbols, and then the way in
%which the result is 
%Actually not quite: Analysis is also needed in the inversion of symbols.
%}
The 'essential structures' are added on the higher levels and encoded
geometrically. 
This parallels the \Mb-calculus way of describing singularity structure (see below)
and is one of its fundamental characteristics:
\medskip

\begin{center}{\bf Fundamental principles of the \Mb-calculus}
\end{center}
\begin{enumerate}
\item 
%\begin{principle}
 \label{DGprincgeom}
Many complications may be understood geometrically,  the analysis
may be reduced to  a few fundamentals.
%\end{principle}
%\begin{principle}
\item \label{DGprincinv}
All concepts which are introduced should be defined in a coordinate-invariant
way. If they depend on choices, the exact freedom in these choices should
be determined. This helps in understanding the concepts themselves.
%\end{principle}
%\begin{principle}
\item \label{DGprincschwartz}
Operators are always described by their Schwartz kernels.
%\end{principle}
%\begin{principle}
\item \label{DGprincb}
All differential objects (e.g.\ densities, differential operators) should
be written as \Mb-objects (i.e.\ using $dx/x, x\partial/\partial x$ etc.\
instead of $dx, \partial/\partial x$ near a boundary $\{x=0\}$).
(This is specific to the \Mb-calculus in the narrow sense.)
%\end{principle}
\end{enumerate}
\medskip

The aim of this article is to explain Figure \ref{DGfigbcalc} and to illustrate
the use and power of these principles.
Our first task is to elucidate what we mean by 'singularity structure'.
The main point will be that complicated behavior of a function 
(or distribution) may often
be described economically by 'blowing up' the underlying space and
then looking at a rather 'simple' function on the new space. The resulting
spaces are {\em manifolds with corners}, and this is the reason for the
central role they play in the \Mb-calculus: They are simultaneously simple
enough to allow for simple analysis, and general enough to describe many
phenomena.

In Section \ref{DGsecgeom} we introduce manifolds with corners and
discuss the singularity structure of smooth functions.

In Section \ref{DGsecana} we discuss the lower
three lines of Figure \ref{DGfigbcalc} and conormal distributions. 
We spend some time to explain the central role played by the Push-Forward
Theorem. As an illustration of the second arrow from below in Figure
\ref{DGfigbcalc}, we define classical pseudodifferential operators (\PDO s) 
and show that they are closed under composition.

Finally, the top levels of  Figure \ref{DGfigbcalc} are addressed in Section
\ref{DGsecpde}.
We first recall the essential ingredients of the classical \PDO\
calculus that permit the construction  of a parametrix for elliptic operators.
We then show how this may be generalized to operators $P$
with  the simplest non-trivial singularity structure, the so-called
\Mb-differential operators on (the interior of) a manifold with boundary. 
In applications, these occur in the context of manifolds
with infinite cylindrical ends or with conical singularities.
Starting from a model calculation we construct the small and full \Mb-calculus
and sketch the parametrix construction.
Note that, up to now,
there does not seem to be a systematic way to construct
a calculus in general
 (i.e.\ finding the right 'Ansatz'). This is the hardest part, and it usually
involves a lot of trial and error.

In the Appendix we collect definitions and basic properties of some objects
which are characteristic for the geometric view of the \Mb-calculus.
\medskip

Prerequisites for Sections \ref{DGsecgeom} and \ref{DGsecana} are minimal
(basic analysis and manifold theory), except for \ref{DGsubsecdist} where
an acquaintance with distributions is assumed. In addition, in Section
\ref{DGsecpde} some basic functional analysis (e.g.\ compact operators) is
needed, and some vague ideas about elliptic operators are useful, 
though not strictly necessary. The many footnotes mostly give additional 
details and can be skipped at first reading.

The reader who wants only a quick impression of the \Mb-calculus should
at least skim the following definitions, remarks, and examples:
 \ref{DGdefmwc}, \ref{DGdefasymp}, \ref{DGdefnice}, \ref{DGremnice},
\ref{DGdefsing}, \ref{DGexsing}.3,
\ref{DGexpfpb}, \ref{DGexpft}, \ref{DGdefbdens}, \ref{DGdefexp1},
\ref{DGdefexp2}, \ref{DGdefbfibr}, \ref{DGthpft}, \ref{DGthpbt},
\ref{DGrempbt}.2, \ref{DGdefbdiff}, \ref{DGexbcalc}; and
Subsections \ref{DGsubsecbu},
 \ref{DGsubsecdist} (if unfamiliar with conormal
distributions), \ref{DGsubsecpdo} (if unfamiliar with classical \PDO s), 
and \ref{DGsubsecsmall}, \ref{DGsubsecfc}.

\medskip
Literature:
R.\ Melrose's 'green book' \cite{DGMel:APSIT} gives a detailed exposition
of the \Mb-calculus (in the narrow sense) on manifolds with boundary; 
the first papers on this were
\cite{DGMel:TBP}, \cite{DGMelMen:EOTCT}, and its extension to manifolds
with corners is discussed in \cite{DGLoy:PDCMWC}, \cite{DGMelPia:AKTMWC}.
Other 'calculi' (i.e.\ singularity structures of $P$, alias boundary
fibration structures) are analyzed in 
\cite{DGEpsMelMen:RLSPD}, \cite{DGHasMazMel:ASAE},
\cite{DGHasVas:RLTOACS}, \cite{DGMaz:HCCCM},
\cite{DGMaz:ETDEOI}, \cite{DGMazMel:ALHCLSSF},
\cite{DGMazMel:ASEI}, \cite{DGMazMel:PDOMWFB},
\cite{DGMcDon:LSWCS}, \cite{DGMel:APSIT},
\cite{DGMel:SSTLAES}, \cite{DGMel:GST}, \cite{DGMel:GOBS},
\cite{DGMelNis:HPDOIMWB}, \cite{DGMoo:HKAMWCS},
usually with applications to problems of geometric analysis.
\cite{DGMel:CCDMWC} gives a condensed presentation of the basic theorems
(Pull-Back and Push-Forward Theorem) on manifolds with corners.
In \cite{DGMel:PDOCSL} (an ICM-talk) boundary fibration structures are
introduced and a general strategy for constructing associated
pseudodifferential calculi is outlined. Of an expository nature are
also \cite{DGMel:GST}, \cite{DGMel:FCAPDO}.
The unfinished and long-awaited
book \Mnotes\ will be the ultimate source for all the details; currently
you can get it on the www, so be quick before it disappears again!
Comparisons with other approaches to singular analysis are made
in \cite{DGGriGru:SALPFT}, \cite{DGLauSei:PDAMWBCBCA} \inthisbook. See
\cite{DGLauSei:PDAMWBCBCA} for many references to
other approaches.
\medskip

Why an article about the \Mb-calculus, given all of these 
beautiful writings? I was told by some
that they would like to learn about the  \Mb-calculus, but find
it hard to get into the style in which it is usually presented: 
Often, things are
expressed in ways that many analysts are not used to.
My aim was to bridge this gap by explaining some of the basics and
highlighting some of the ideas which are usually hidden between the lines.
I hope to make this beautiful world accessible to a larger audience.

Some specific points in which this presentation differs from others are:
\begin{itemize}
\item
I propose a notion of 'asymptotic type' of a function as a blow-up
under which it becomes polyhomogeneous conormal (Definition \ref{DGdefsing}). 
While this is clearly
implicit in existing treatments, the explicit notion suggests naturally
the problem of determining a type of a push-forward (or pull-back)
of $f$ from
a type for $f$. Melrose's Push-Forward and Pull-Back Theorems answer 
this only partially
(see Remarks \ref{DGrempft}.2 and \ref{DGrempbt}).
\item
I give an alternative definition of the central notion of \Mb-fibration,
which I believe to be more intuitive (see Definition \ref{DGdefbfibr}).
\item
I discuss the relation between the notions of 'type' and 'regimes';
the former originates in pure mathematics (algebraic geometry)
while the latter is widely used in applied mathematics.
\item
For reasons of space I do not discuss the general notion of 
boundary fibration structures, nor any other instances
besides the \Mb-\PDO s. Also missing are
the \Mb-vector fields and associated \Mb-bundles
(but they are implicit in Principle \ref{DGprincb}) and the \Mb-Sobolev spaces.
\end{itemize}

\medskip

\begin{center} {\sc Acknowledgement} \end{center}
I am deeply grateful to R.\ B.\ Melrose for introducing me into this
world. Clearly, all the important ideas are due to him. I hope
he will not disagree with the particular slant given to some things here,
which reflect my own understanding and interests.

%% file: DGgeometry.tex
%%%%%%%%%%%%%%%%%%%%%%%%%%%%%%%%%%%%%%%%
%%%%%%%%%%%%%%%%%%%%%%%%%%%%%% Section Geometry
%%%%%%%%%%%%%%%%%%%%%%%%%%%%%%%%%%%%%%%%

\section{Geometry} \label{DGsecgeom}
We begin with the \Mb-calculus way of describing the asymptotic
behavior of smooth functions. Thus, we are given a non-compact
manifold $Z_0$ and want to find a 'good' way to describe how the value
$u(z)$ of a smooth function 
$u:Z_0\to\C$ (or $\R$) behaves when $z$ approaches the 'boundary' of
$Z_0$, i.e.\ leaves any compact subset of $Z_0$.
Here are a few examples, along with naive attempts to describe their
asymptotic behavior:
\begin{examples} \label{DGexfcns}
\ 
\begin{enumerate}
\item 
$Z_0=(0,\infty)$, $u(x)=1/x$: has a 'first order pole' at zero and
vanishes to first order at infinity.
\item 
$Z_0=(0,1)^2$, $u(x,y)=xy$: extends smoothly to a neighborhood of
$[0,1]^2$ (i.e.\ the asymptotics is given by Taylor expansion around
any boundary point of $[0,1]^2$).
\item 
$Z_0=(0,1)^2$,$u(x,y) = 1/xy$: similar to 2., except that negative
powers are allowed.
\item
$Z_0=\R^2\setminus\{(0,0)\}$, $u(x,y)=\sqrt{x^2+y^2}$:
'decays linearly to $(0,0)$ from all directions'.
\item 
$Z_0=(0,\infty)^2$, $u(x,y)=\sqrt{x^2+xy+y^3}$: smooth at the coordinate axes
except at $(0,0)$; near $(0,0)$: complicated ($x^2$ dominates for $x>y$, $xy$
for $y>x>y^2$, and $y^3$ for $y^2>x$).
\item 
$Z_0=\{(x_1,x_2,x_3)\in \R^3:\, x_1^2+x_2^2 = x_3^{2},\, x_3>0\}$,
$u(x_1,x_2,x_3) = x_3$: similar to 4.
\end{enumerate}
\end{examples}
(In Examples 4.-6.\ we did not consider the behavior at infinity.)

\vspace{1em}
\begin{list}{}{}
\item[{\em Observation 1:}]
It is useful to add certain 'boundary' points to $Z_0$, so that one can talk
for example about the behavior of a function 'at $(0,0)$' when actually
referring to its behavior in $Z_0\cap U$ for $U$ an arbitrarily small
neighborhood of $(0,0)$. 
(In Example 1 this means adding a 'point at infinity' also.)
\item[{\em Observation 2:}]
Descriptions of asymptotic behavior must refer to certain
coordinates (e.g.\ in Example 1 the standard coordinate on $\R$, both near
zero and near infinity; in Example 4 it is natural to use polar coordinates,
then $u$ is just a smooth function of $r\geq 0$, vanishing at zero; for
Example 5 it is less clear what 'good' coordinates would be).
\end{list}

\vspace{1em}
To explain the \Mb-calculus description of asymptotics, 
we proceed in two steps:
\begin{list}{}{}
\item[First step:]
We introduce what is considered model behavior: $Z_0$ is the interior
of a manifold with corners $Z$, and the functions  have
joint asymptotic expansions in all variables in the corners.
Melrose calls these functions 'polyhomogeneous conormal'. We prefer short
words here, so we call them 'nice'.
\item[Second step:]
We show how more general asymptotics may (often) be reduced
to this model case by
specifying an identification (diffeomorphism) of $Z_0$ with
the interior of some manifold with corners $W$.
The most common way to define such a diffeomorphism is by blow-up, which
we also discuss.
\end{list}

\vspace{1em}
We will see that this gives a very geometric way to describe the 
'asymptotic type' of a function. 
Good references for this section are \cite{DGMel:CCDMWC} and \Mbook\ 
(besides the all-encompassing unpublished \Mnotes).
\subsection{Manifolds with corners and nice functions}
\ \vspace{1ex}
\begin{definitions} \label{DGdefmwc}
\
\begin{enumerate}
\item 
A {\em manifold with corners} (mwc) is a topological space locally modelled
on pieces of the form $[0,\infty)^k\times \R^{n-k}$, for
various $k\in\{0,\ldots,n\}$ (in the same sense as
a manifold is modelled on pieces $\R^n$ and a manifold with boundary (mwb)
on pieces $\R^n$ and $[0,\infty)\times\R^{n-1}$).
\item
A mwc $Z$ is the union of its interior $\Zint$ and its boundary $\partial Z$.
The boundary is the union of the {\em boundary hypersurfaces} (bhs's) of $Z$
which are themselves mwc's.%
\footnote{\label{DGfnmwc}%
When mwc's are defined this way, a bhs may happen to be only immersed rather
than embedded, see Figure \ref{DGfigmwc}.
In the \Mb-calculus it is also always assumed that
the boundary hypersurfaces are embedded (and connected); 
this detail won't matter for a while.}
\item 
A {\em boundary defining function} (bdf) of a boundary hypersurface
$H$ of $Z$ is a function $\rho:Z\to[0,\infty)$ such that
$\rho^{-1}(0)=H$, $\rho$ is smooth up to the boundary,
 and $d\rho\not=0$ on $H$. (See below for the definition
of smoothness.)%
\footnote{
The assumption on $H$ to be embedded implies the existence
of a bdf for $H$.}
\end{enumerate}
\end{definitions}
\begin{figure}[htbp]
\input{figmwc.pstex_t}
\caption[Embedded boundary hypersurfaces?]{}\label{DGfigmwc}
\end{figure}
The simplest examples are $\R,\R_+:=[0,\infty),\R_+\times\R$ (all mwb)
and $\R_+^2$ (the simplest corner). 
Also, $[0,1]^n$ is a mwc, but $Z=\overline{Z_0}$ in Example
\ref{DGexfcns}.6 is not.
See Figure \ref{DGfigtriple} for a more complicated mwc.
On $\R_+$, a bdf is given by $\rho(x)=x$. But note that many others are
possible.
The cartesian product of two mwc's is again
a mwc (and the product of two mwb's is a mwc, but not a mwb, which is
one reason for introducing the notion of mwc; for another reason see
\ref{DGsubsubsecinv}).

Every mwc can be embedded in a manifold: For example $\R_+^k\times\R^{n-k}
\hookrightarrow \R^n$ in the obvious way; by definition,
an embedding is a map which looks locally
like that.

If $Z$ is a mwc then we will speak of a 'function $u$ on $Z$' even if $u$ 
is only defined in the interior of $Z$.

The role of bdf's is that they are the coordinates in terms of which the
asymptotic behavior of functions will be described. 

We now define nice functions. These should be thought
of as slightly more general than functions smooth up to the boundary, so
we discuss these shortly.

\subsubsection{Functions smooth up to the boundary}
These are, per definition, restrictions to $Z$ of
smooth functions on $M$, where $Z\hookrightarrow M$ is some embedding
into a manifold. However, it is desirable to characterize this intrinsically,
just using the values of the function on the interior $\Zint$. 
{\em Seeley's extension theorem} (see \Mnotes) says that 
\begin{quote}
$u$ is smooth up to the boundary iff all derivatives of all orders of $u$
are bounded on bounded subsets of $\Zint$.
\end{quote}
(Of course, a bounded subset is one whose closure in $Z$ is compact.)
For a characterization in terms of asymptotics see Remark 
\ref{DGremasymp}.3 below.

\subsubsection{Nice (polyhomogeneous conormal) functions}
These are functions that behave like sums of products of
terms like $x_i^\alpha\log^p x_i$, $i\in I$, near a corner defined by 
$\{x_i\geq 0, i\in I\}$. We call them nice since their behavior
under integration can be analyzed fairly easily (see Section 
\ref{DGsecana}), and a
discrete set of numbers is sufficient to describe their asymptotic behavior
completely. 

To set the stage, we consider a manifold with boundary first:

\begin{definition} \label{DGdefasymp}
Let $Z$ be a manifold with boundary $H=\partial Z$. 
\begin{enumerate}
\item An {\em index set} is a discrete subset $F\subset\C\times\N_0$
such that every 'left segment' $F\cap\{(z,p): \Re z< N\}, N\in\R$ is a finite set.
Also, it is assumed that $(z,p)\in F, p\geq q\implies (z,q)\in F$.
\item 
Given an index set $F$, a smooth function $u$ on $\Zint$ is called
polyhomogeneous conormal (in short, {\em nice})
with respect to $F$ if, on a tubular
neighborhood $[0,1)\times H$ of $H$, one has
\begin{equation} \label{DGeqasymp}
u(x,y) \sim \sum_{(z,p)\in F} a_{z,p}(y)x^z \log^p x\quad \text{ as }x\to0
\end{equation}
with $a_{z,p}$ smooth on $H$.%
\footnote{
The meaning of $\sim$ is this: Let $u_N(x,y)$ be the sum
of \eqref{DGeqasymp} restricted to $\Re z\leq N$. Then
\begin{equation} \label{DGeqremest}
|u(x,y) - u_N(x,y)| \leq C_Nx^N
\end{equation}
for all $N$, uniformly on compact subsets of $H$, plus analogous estimates
when  taking any number of $x\partial_x$ and $\partial_y$ derivatives.}%
\footnote{\label{DGfncoordind}%
It is easy to check that this definition is independent of the choice
of identification of a neighborhood of $H$ with $[0,1)\times H$,
if one assumes that $F$ satisfies 
\begin{equation}\label{DGeqcoordind}
(z,p)\in F\Rightarrow (z+1,p)\in F.
\end{equation}
In any case, any index set can be 'completed', i.e.\ enlarged to a smallest
index set having this property.}%
\footnote{
Instead, one could consider finite asymptotics, but it messes up the
notation. We prefer complete asymptotics (and $C^\infty$ functions)
so we can focus on more important things.
}
\end{enumerate}
\end{definition}

Thus, an index set tells us which terms $x^z\log^p x$ may occur in
the asymptotics of $u$ at the boundary. The finiteness condition on $F$ 
ensures that \eqref{DGeqasymp} makes sense.
$\log$'s are allowed since they often appear naturally, e.g.\ when integrating
smooth functions (see Example  \ref{DGexpft}.2).
Note that $x$ is a bdf for $H$ (at least near $H$, which is all that matters).

\begin{remarks}[Contents and limitations of Definition \ref{DGdefasymp}]
\label{DGremasymp}
\
\begin{enumerate}
\item In the simplest case of $Z=\R_+$, we allow functions like
$x^{-3},\log x$, but no 'fast oscillation' like $\sin 1/x$. But $e^{-1/x}$
is nice (for any index set $F$, e.g. $F=\emptyset$).
\item The exponents $\alpha$ are not allowed to depend on $y$; thus
the 'variable asymptotics' of Schulze 
(\cite {DGSchul:PDOMWS}, Section 2.3) is excluded.
\item Exercise: $u$ is nice with respect to the index set $0:=\{(n,0): n\in\N_0\}$
iff $u$ is smooth up to the boundary.
\end{enumerate}
\end{remarks}

We now turn to manifolds with corners, and this is where the story gets really
interesting.

For simplicity, we will only consider the mwc $\R_+^2$. The extension to
the general case is not difficult, see for example
\cite{DGEpsMelMen:RLSPD}, \cite{DGMazMel:ASEI}, \Mnotes, \cite{DGMel:CCDMWC}
 (for corners of higher codimension use induction, for additional
$\R$-variables assume smooth dependence).

An {\em index family $\Ecal$} for a mwc $Z$ is an assignment of an index set to
each bhs. For $Z=\Rplus^2$ we simply write $\Ecal=(E,F)$ if $E$ is
associated with the $x$-axis (i.e.\ $\{y=0\}$)
and $F$ is associated with $y$-axis (i.e.\ $\{x=0\}$). 

\begin{definition} \label{DGdefnice}
Let $(E,F)$ be an index family for $\R^2_+$.
 A function $u$ on $\R^2_+$ is polyhomogeneous conormal (in short, {\em nice})
with respect to
$(E,F)$ if it has an asymptotic expansion in $x$
as in \eqref{DGeqasymp},
%\begin{equation} \label{DGeqpoly}
%u(x,y) \sim \sum_{(z,k)\in E} a_{z,k}(y) x^\alpha \log^k x
%\end{equation}
where the coefficients $a_{z,p}$ are functions on $\R_+$ that are
nice with index set  $E$ (in the sense defined above for the mwb $\R_+$).%
\footnote{
Here, a good definition of $\sim$ is harder to come by since
the remainder in \eqref{DGeqremest} should be allowed to be singular in
$y$ (at $y=0$), but not too singular. One way around this is to require
$$ |u(x,y)-u_N(x,y)| \leq C_Ny^{-M}x^N$$
for some fixed $M$ and all $N$ (plus similar remainder estimates for 
the derivatives), {\em plus} an analogous expansion and estimate with $x$ 
and $y$ interchanged. Again, one has coordinate invariance. See also
\ref{DGsubsubsecinv}.}
\end{definition}
Again, it is easy to see that $u$ is smooth on $\R^2_+$ iff $u$ is
nice with respect to the index family $(0,0)$, with $0$ from
Remark \ref{DGremasymp}.3. Examples \ref{DGexfcns}.1-3 are nice.
(However, in Example 1 this describes only the behavior near zero;
see below for the behavior near infinity.)
\begin{remark} \label{DGremnice}
The MAIN POINT is that all coefficients in the expansion \eqref{DGeqasymp}
lie in the same singularity class. As a non-example, let us consider the  
function $u(x,y)=\sqrt{x^2+y^2}$ on $\R^2_+$. Clearly, $u$ extends smoothly
to the boundary except at $(0,0)$. Therefore, for each
fixed $y>0$, one has an asymptotic expansion
$$ \sqrt{x^2+y^2} \sim \sum_{i=0}^\infty a_i(y)x^{i} \quad\text{as }x\to 0,$$
and similarly with $x$ and $y$ interchanged.
But the coefficient functions $a_i(y)$ become more and more singular 
as $i\to\infty$. This is seen  easily by writing
\begin{align} \label{DGeqsqrt}
 \sqrt{x^2+y^2} & = y\sqrt{1+(x/y)^2} = y\sum_0^\infty c_i(\frac{x}y)^{2i}\\
 &= y + \frac12\frac{x^2}y - \frac18\frac{x^4}{y^3} + \ldots
\end{align}
with the Taylor series 
$\sqrt{1+t}=\sum_0^\infty c_i t^i = 1+t/2 - t^2/8+\ldots$ (for $|t|<1$).
It is easy to see from this that $u$ is not nice with respect to any index
family.

Therefore, niceness means having a 'joint' (or uniform) asymptotic expansion,
simultaneously as all variables tend to zero, in the corner.
\end{remark}

\subsection{More general asymptotic behavior. Asymptotic type.}
Though $u(x,y)=\sqrt{x^2+y^2}$ is not nice on $\Rplus^2$, 
we saw in \eqref{DGeqsqrt}
(and its analogue with $x$ and $y$ interchanged) that:
\begin{equation} \label{DGeqniceas}
\begin{split}
 u &\text{ is  nice as a function of $y$ and $x/y$ (for $x/y$ bounded), and} 
 \\
 u & \text{ is nice as a function of $x$ and $y/x$ (for $y/x$ bounded).}
 \\
\end{split} 
\end{equation}

We want to take this as a characterization of the 'asymptotic type' of
$u$. A beautiful way to do this is through the following construction:

\begin{definition} \label{DGdefsing}
Let $Z_0$ be any manifold. Let
$W$ be a compact 
mwc and $\beta:\Wint\to Z_0$ a diffeomorphism. Call a function $u$
on $Z_0$ {\em of (asymptotic) 
type $\beta$} if $\beta^*u$ is a nice function on $W$.%
\footnote{
Here $\beta^*u=u\circ\beta$, the pull-back. One can think of $\beta$ as a 
distortion lens, then $\beta^*u$ is simply $u$, looked at through this lens.
Of course, a function has many types.
}
\end{definition}
We then also say that  {\em $u$ is resolved by $\beta$}. Some people would
call $\beta$ a 'singular coordinate change'.
Of course, we may also specify an index family $\Ecal$ on $W$ and then
speak of 'type $\beta$ with index family $\Ecal$'.
The compactness of $W$ means that the asymptotics of $u$ 
is controlled 'in all directions'. We will freely consider non-compact $W$
as well, when we are only interested in the behavior of $\beta^*u$ on
a compact part of $W$.
\begin{examples} \label{DGexsing}
\
\begin{enumerate}
\item
To describe the behavior of  Example \ref{DGexfcns}.1 at infinity, one should
say in which sense $[0,\infty]$ is a mwc (actually, a mwb).
This may be done by
choosing any diffeomorphism $\beta:(0,\infty)\to(0,1)$ which is equal to
the identity near zero and to the map $x\mapsto 1-1/x$ near infinity%
\footnote{
This corresponds to the common usage in complex analysis,
where behavior of $u(z)$ 'at infinity' is described by behavior of $u(1/z)$ 
at zero.}.
Then $u$ from Example \ref{DGexfcns}.1
has type $\beta$ with index sets $\{(-1,0)\}$ at zero and $\{(1,0)\}$ at one.
\item
 Let $\beta:\R^2_+\to\R^2_+, (\xi,\eta)\mapsto (\xi\eta,\eta)$.
Then $u$ has type $\beta$ iff $u(\xi\eta,\eta) = v(\xi,\eta)$ with
$v$ nice. Writing $\xi\eta=x, \eta=y$ we see that this means exactly
that $u$ is nice as a function of $y$ and $x/y$ (for bounded $y$ and $x/y$).
\item
 (Polar coordinates) Let $W=\R_+\times [0,\pi/2]$, $Z=\R_+^2$ and 
\begin{equation}\label{DGeqpolar}
\beta (r,\theta) = (r\cos\theta,r\sin\theta).
\end{equation}
If $u(x,y)=\sqrt{x^2+y^2}$ on $\Rplus^2$ then
 $\beta^*u(r,\theta)=r$, so 
$u$ has type $\beta$. We will see below that \eqref{DGeqniceas}
is equivalent to $u$ having this type $\beta$ (see Remarks
\ref{DGremcoords}).
Of course, the same formula \eqref{DGeqpolar} works for
$Z_0=\R^2\setminus\{(0,0)\}$ and $W=\Rplus\times S^1$
(with $S^1=[0,2\pi]/0\sim2\pi$ the circle), which makes $u$ in 
Example \ref{DGexfcns}.4 of type $\beta$.

Note that by considering $\beta^*u$ we 'spread out' the values of $u$
near $0$ over a whole strip (a neighborhood of $\{0\}\times S^1$).
\item
 Let $\beta:\R_+\to\R_+, \xi\mapsto e^{-1/\xi}$. Then $\beta(\xi)=x$
iff $\xi=(\log 1/x)^{-1}$, so $u$ has type $\beta$ iff $u$ is nice
as a function of $(\log 1/x)^{-1}$. This is used in \cite{DGHasMazMel:ASAE},
for example.
\item
Exercise: Find a type for Examples \ref{DGexfcns}.5 and 6 (near $(0,0)$).
\end{enumerate}
\end{examples}

Often, $Z_0$ is given as dense subset of a mwc $Z$, but the functions
$u$ of interest are not nice (see Example \ref{DGexfcns}.4 when considered
on $\Rplus^2$, or Example 5). In this case, a space $W$
and map $\beta$ can often be obtained by a procedure called
'blow-up'. 
We discuss this next.

\subsection{Blow-up} \label{DGsubsecbu}
Blow-up is a way to obtain new mwc's from old. It is used to resolve
functions on a mwc and to desingularize
(algebraic) subsets of a mwc (see Subsection \ref{DGsubsecemb}).

The simplest non-trivial
case of a blow-up is given by polar coordinates (Example \ref{DGexsing}.3). 
We will discuss this case in some detail and then sketch the general construction.

\subsubsection{Blowing up $(0,0)$ in $\Rplus^2$}
Consider the 'polar coordinates map' \eqref{DGeqpolar}. 
Note that, for $p\in\R^2_+$, $\beta^{-1}(p)$ is a point unless $p=(0,0)$ when
it is the interval $\{0\}\times [0,\pi/2]$. Therefore, we say that $W$ is
obtained from $Z$ by 'blowing up $(0,0)$'. We write $W=[Z,(0,0)]$ and call
$\beta$ the blow-down map. The bhs's of $W$ are called $\lb=\{\theta=\pi/2\}$, 
$\rb=\{\theta=0\}$ (the left and right boundary) and
$\ff=\{r=0\}$ (the front face).
\begin{figure}[htbp]
\input{figpolar.pstex_t}
\caption{The blow up $[\R^2_+,(0,0)]$ of $\R^2_+$}\label{DGfigpolar}
\end{figure}
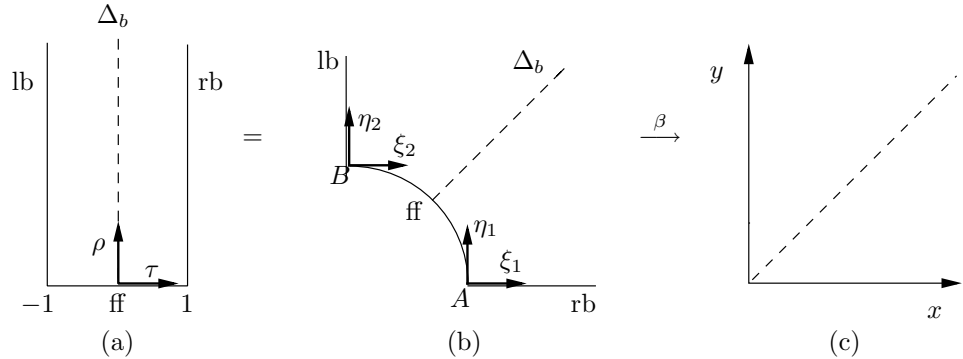

When drawing pictures, some people prefer to draw $W$ as in 
Figure \ref{DGfigpolar}(b) while others prefer \ref{DGfigpolar}(a).
In spirit these correspond roughly to using two different sets of coordinate
systems on $W$, which are often more convenient to use than $(r,\theta)$:%
\footnote{
The occurrence of the transcendental functions $\sin$ and $\cos$ in
\eqref{DGeqpolar}, with all their special properties (e.g. $\sin'=\cos$)
is rather accidental and usually distracts from what
really matters, e.g.\ the asymptotic behavior when approaching the boundary.
There is no way to completely erase such accidents, but the following
two options come close to it.}

\bigskip

{\em Coordinate systems on $[\R_+^2,(0,0)]$}
\begin{enumerate}
\item  (Projective coordinates)
For $y \ll x$, we have $r=\sqrt{x^2+y^2}\approx x$
and $\theta=\arctan y/x\approx y/x$. This suggests considering
\begin{subequations}\label{DGeqproj}
\begin{equation}
 \xi_1= x,\,\eta_1=\frac{y}x
\end{equation}
as coordinates on $W$. 
Indeed, from  $\xi_1=r\cos\theta,\eta_1=\tan\theta$ we see that they  
define local coordinates for $\theta\not=\pi/2$, i.e.\
 on $W\setminus \lb$, with $\xi_1$ a bdf for \ff\
and $\eta_1$ a bdf for \rb. Similarly,
\begin{equation}
\xi_2=\frac{x}y,\,\eta_2=y
\end{equation}
define coordinates on $W\setminus \rb$, with $\xi_2$ a bdf for \lb\
and $\eta_2$ a bdf for \ff.
In these two coordinates systems $\beta$ takes the simple form
\begin{align} 
\begin{split}
\beta_1(\xi_1,\eta_1)&=(\xi_1,\xi_1\eta_1)\\
\beta_2(\xi_2,\eta_2) &= (\xi_2\eta_2,\eta_2).\\
\end{split}
\end{align}
\end{subequations}
\item ('Rational polar coordinates')
Define (for $x,y>0$)
\begin{equation}\label{DGeqglobalcoord}
\rho=x+y, \quad \tau = \frac{x-y}{x+y}.
\end{equation}
Writing $x=r\cos\theta,y=r\sin\theta$ one easily sees that
$\rho=r a(\theta), \tau=b(\theta)$ with $a>{\rm const}>0$ and $b:[0,\pi/2]\to[-1,1]$
a diffeomorphism, thus $(\rho,\tau)\in\R^+\times[-1,1]$ 
may be regarded as new coordinates on $W$.

Solving \eqref{DGeqglobalcoord} for $x,y$ one obtains the form of
the blow-down map as
$$ \beta(\rho,\tau) = (\frac12\rho(1+\tau),\frac12\rho(1-\tau)).$$
Bdf's are given by $\rho$ for \ff, $1+\tau$ for \lb, and $1-\tau$ for \rb.
\end{enumerate}

\begin{remarks}\label{DGremcoords}
\
\begin{enumerate}
\item
Rather than beginning with polar coordinates one may define $W$ and $\beta$
directly by glueing two coordinate patches, i.e.\ 
\begin{equation}\label{DGeqident}
W=\R^2_+ \sqcup \R^2_+ /\sim 
\end{equation}
$$ \text{ where }
(\xi_1,\eta_1)\sim(\xi_2,\eta_2) :\Leftrightarrow \beta_1(\xi_1,\eta_1)=
\beta_2(\xi_2,\eta_2)
$$
with $\beta_{1/2}$ from (\ref{DGeqproj}c).
This identification is done precisely in order for $\beta$ to be injective on
$\Wint$. Note that injectivity is essential for the whole idea of defining asymptotic
types of functions on $Z_0$ using the map $\beta$.
\item
The first remark together with (\ref{DGeqproj}a) and (\ref{DGeqproj}b)
shows that \eqref{DGeqniceas} holds iff $\beta^*u$ is nice
on $[\R_+^2,(0,0)]$.
\item
\eqref{DGeqident} is the way that blow-up is usually defined in 
algebraic geometry, except that $\R_+$ is replaced by $\C$, so that
all spaces involved are smooth complex varieties without boundary.
\item 
An advantage of the $(\rho,\tau)$ coordinates is that they
are global on $W$. While the projective coordinates may feel cumbersome
at first since they are not global, they have several advantages:
The bdf's are simply the $\xi_i$ and $\eta_i$,
calculations tend to be very simple, and in many problems they occur naturally
(see \cite{DGGriGru:SALPFT}, for example).
\end{enumerate}
\end{remarks}

\subsubsection{More general blow-ups} \label{DGsubsubsecgenbu}
In general, given  mwc's $Z$ and $Y\subset Z$ (satisfying certain conditions), 
one constructs $W=[Z,Y]$, the 'blow-up of $Z$ along $Y$', together with
a smooth map $\beta:W\to Z$, the 'blow-down map', 
which is a diffeomorphism $\Wint\to \Zint\setminus Y$.
This is done as follows:
\begin{enumerate}
\item
Blow-up of an interior point of a 2-dimensional mwc replaces it by
a circle, see Example \ref{DGexsing}.3.
\item
In higher dimensions, blow-up of an interior point replaces it by a sphere. For 
example, $$[\R^n,0] = \R_+\times S^{n-1},\quad \beta(r,\omega) = r\omega$$
for $r\in\Rplus,\omega\in S^{n-1}$.
Similarly,  $[\R_+^n,0] = \R_+\times S^{n-1}_+$ where
$S^{n-1}_+=S^{n-1}\cap\R^n_+$.
$(r,\omega)$ provide polar coordinates on $[\R^n,0]$.
Projective coordinates are
$\xi_1=x_1,\xi_2=x_2/x_1,\ldots,\xi_n=x_n/x_1$
on $\{\omega_1\not=0\}$, and similarly on all other $\{\omega_i\not=0\}$.
\item
More generally, one can blow up closed submanifolds $Y\subset Z$:
Assume first that $Y$ lies in the interior of $Z$. Locally, the pair
$(Y,Z)$ is just $(\R^k\times\{0\}^{n-k},\R^n)$, and we simply set
$$[\R^n,\R^k\times\{0\}^{n-k}] = \R^k\times [\R^{n-k},0],$$
with $\beta$ as in 2. above. One can check that this is independent of
the coordinates chosen (up to diffeomorphism that intertwines  the $\beta$'s)
and can therefore be glued together to a global blow-up $\beta:[Z,Y]\to Z$.
A more intuitive model for the space $[Z,Y]$ is  $Z\setminus U_r(Y)$, where
$U_r(Y) = \{z\in Z: \dist(z,Y)<r\}$ with respect to some Riemannian
metric on $Z$, for $r$ sufficiently small (at least when $Y$ is compact).
(But then $\beta$ is more complicated to write down.)
\item
The construction from 3. can be extended directly to mwc's 
$Y$ hitting the boundary of $Z$, if this 'hitting' is transversal
in a suitable sense. (As a non-example, consider the parabola
$Y=\{x,x^2\}\subset Z=\R\times\R_+$ and try to define a blow-up!)
The exact condition is that
near any $p\in Y$ coordinates can be chosen such that $p=0$ and, locally,
$Z=\Rplus^k\times\R^{n-k}$ and $Y=Z\cap S$ for some coordinate subspace $S$. 
%% the pair $(Y,Z)$  locally  look like
%$((\R_+^l\times \{0\}^{k-l}) \times (\R^m\times\{0\}^{n-k-m}),\R_+^k\times\R^{n-k})$
%for some $k,l,m$. 
Such $Y$ are called {\em p-submanifolds}.
$[Z,Y]$ is a mwc.
\end{enumerate}

\begin{remarks}\label{DGrembu}
\
\begin{enumerate}
\item
As already indicated, all these blow-ups are defined invariantly (i.e.\ no
choices are made, beyond $Z$ and $Y$, to define $[Z,Y]$ up to diffeomorphism
that preserves $\beta$).
\item 
The construction above yields an {\em elementary blow-up}.
Sometimes, one needs an {\em iterated blow-up}. That is, one chooses
a p-submanifold $Y'$ in $[Z,Y]$ and considers $[[Z,Y],Y']$, or iterates even
further.
This is needed when describing
more complicated asymptotics of functions. For example, blowing
up first 
$(0,0)$ in $\Rplus^2$ and then the point $B$ in Figure \ref{DGfigpolar}(b)
resolves the function $u(x,y)=\sqrt{x^2+xy+y^3}$ from Example
\ref{DGexfcns}.5 (exercise!). 
Also, the 'triple \Mb-space' in the \Mb-pseudodifferential
calculus is an iterated blow-up (see Figure \ref{DGfigtriple}).
\item
When using an (iterated) blow-up $\beta:W\to Z$ to describe asymptotic
behavior of functions we have:
\begin{enumerate}
\item[(a)]
If $u$ is nice then $u$ is of type $\beta$.
\item[(b)]
The behavior of $u$ on a compact part of $Z$ is reflected by the behavior of
$\beta^*u$ on a compact part of $W$. (Cf.\ the remark on compactness
after Definition \ref{DGdefsing}.)
\end{enumerate}
(a) follows from the fact that $\beta$ is a \Mb-map (see below), and 
(b) is just the properness of $\beta$ (i.e.\ $\beta^{-1}({\rm compact})
= {\rm compact}$).
\end{enumerate}
\end{remarks}

\subsubsection{\Mb-maps}
An important property of blow-down maps is that they are \Mb-maps. 
We define these now.
Recall that if $w$ is any point in a mwc $W$ then
a neighborhood of $w$ can be identified with $\R_+^k\times\R^{n-k}$, with
$w$ corresponding to $0$ ($k$ depends on $w$).
\begin{definition} \label{DGdefbmap} %% b-maps
A map $f:W\to Z$ between mwc's is a {\em \Mb-map} at $w\in W$ if
for some (and therefore any) identification of neighborhoods of $w$ and $z=f(w)$ with
$\R_+^k\times\R^{n-k}$ and $\R_+^{k'}\times\R^{n'-k'}$, respectively,
sending $w,z$ to zero, the map has 'product type', i.e.
$$ f=(f_1,\ldots,f_{n'}),$$
\begin{equation}\label{DGeqbmap}
 f_i (x_1,\ldots,x_k,x_{k+1},\ldots,x_n) = 
     a_i(x) \prod_{j=1}^k x_j^{\alpha_{ij}},\quad \text{ for }i=1,\ldots,k', 
\end{equation}
with $a_i$ smooth and non-vanishing near zero, and non-negative integers
$\alpha_{ij}$.

$f$ is a \Mb-map if it is a \Mb-map at every point.%
\footnote{
Melrose calls such maps {\em interior} \Mb-maps. For a general \Mb-map
he allows that instead of \eqref{DGeqbmap} one has $f_i\equiv 0$ for
some $i$, i.e.\ that $f(W)\subset \partial Z$ (assuming $W$ connected).
In this article we never use these general \Mb-maps.
}
\end{definition}
In particular, \Mb-maps are smooth up to the boundary.
Examples \ref{DGexsing}.2 and 3 are \Mb-maps, while 4 is not. 
\begin{remarks} (Intuition and properties of \Mb-maps)
\label{DGrembmap}
\
\begin{enumerate}
\item
Condition \eqref{DGeqbmap} is an algebraic counterpart to
the weaker geometric condition that, near $z$, the zero set $f_i^{-1}(0)$ is
a union of bhs's through $z$ (more globally: The preimage of any bhs of
$Z$ is a union of bhs's of $W$). See \ref{DGsubsecpf} for a more detailed
discussion of the boundary geometry of \Mb-maps.
\item
The composition of \Mb-maps is a \Mb-map. Projective coordinates show
that the blow-down map for an elementary blow-up (and therefore for any
blow-up) is a \Mb-map.
\item
 If $\beta$ is a \Mb-map then $\beta^*({\rm nice})=$ nice
 since
$\log xy = \log x + \log y$ and $\log a $ is smooth for $a>0$.
For a more precise statement see the 'pull-back theorem', Theorem \ref{DGthpbt}.
\item
If $\beta$ is just smooth then $\beta^*({\rm smooth})=$ smooth, but
in general $\beta^*(\rm{nice})$ is not nice; for example, for
$\beta:\R_+^2\to\R_+,(x,y)\mapsto x+y$ and $u(t)=\log t$ we get
$\beta^*u(x,y) = \log (x+y)$ which is not nice on $\Rplus^2$. Similarly,
$\sqrt{x+y}$ is not nice.
\end{enumerate}
\end{remarks}

\subsection{Embedded blow-up} \label{DGsubsecemb}
So far, we have not addressed Example \ref{DGexfcns}.6.
$Z=\overline{Z_0}$ is not a manifold with corners, so the blow-up
construction above does not apply directly to the construction of
an appropriate 'blow-up space' $W$. However, $Z_0$ is embedded in $\R^3$,
which {\em is} a manifold. So one may blow up 0 (the singular point of $Z$)
 in $\R^3$ and then take $W$ to be the closure of the preimage of 
$Z_0$:
$$ W:=\overline{\beta^{-1}(Z_0)},\quad\text{ with } \beta:[\R^3,0]\to\R^3
\text{ the blow-down map.} $$
Using polar coordinates on $[\R^3,0]=\R_+\times S^2$, i.e.\
$\beta(r,\omega) = r\omega$, $\beta^{-1}(x) = (|x|,x/|x|)$, we get
$$ W = \Rplus\times C, \quad C=\{\omega_1^2+\omega_2^2=\omega_3^2,\,
\omega_3>0\}\subset
S^2. $$
$C$ is a smooth curve -- a circle -- on $S^2$, so $W$ is a mwb
and $\beta:\Wint\to Z_0$ a diffeomorphism, and
$$ \beta^*x_3 = r\omega_3.$$
Since $\omega_3$ is a smooth function on $S^2$, this is a nice function
by Definition \ref{DGdefasymp}.

This procedure is called 'embedded blow-up' (or embedded desingularization).
Hironaka showed in his famous 'resolution of singularities' work that
such an embedded blow-up exists for any (semi-)algebraic set 
(and (semi-)algebraic function on it) in $\R^n$,
and can be obtained by  an iterated blow-up. 
(These authors use the 'projective'
blow-up, but it should be easy to transfer the 
result to our situation.) See \cite{DGHir:RSAVFCZ}, \cite{DGBieMil:CDCZBUMSLI}.

\subsection{Invariance, regimes, etc.} \label{DGsubsecremgeom}
Here we collect some more remarks on the idea of 'asymptotic type'.

\subsubsection{On invariance} \label{DGsubsubsecinv}
Definition \ref{DGdefnice} (and its generalization to any mwc) is
coordinate invariant if all index sets in $\Ecal$ satisfy the
condition \eqref{DGeqcoordind}. This means: Let $\xtilde,\ytilde$
be any other bdf's for the $y$- and $x$-axis in $\Rplus^2$, respectively.
(In particular $\xtilde,\ytilde$ define coordinates near $(0,0)$.)
Then $u$ is nice with index family $\Ecal$ when expressed
in terms of $\xtilde,\ytilde$ iff it is in terms of $x,y$.
The reason is that both $x/\xtilde$ and $y/\ytilde$ are smooth and non-zero
up to the boundary.

Thus, although coordinates (i.e.\ bdf's) are needed to write down the 
particular asymptotics of $u$,
\begin{quote} %\label{DGeqniceclass} (LABEL)
{\em the class of nice functions on a mwc with a given index family $\Ecal$ 
is defined independent of coordinates,}
\end{quote}
and therefore defined purely by the geometry (the mwc) and the discrete
set $\Ecal$.

In contrast, there is no 'natural' class of coordinate functions
describing approach to $(0,0)$ in $\R^2$: Both $\sqrt{x^2+y^2}$
and $\sqrt{x^2+2y^2}$ would be equally good candidates as
'defining functions of $(0,0)$', but their quotient does not extend
smoothly to $(0,0)$.

This shows the special role played by mwc's and is one reason for their
central role in the \Mb-calculus.

\subsubsection{On invariance, II} \label{DGsubsubsecinv2}
Because of the invariance of the blow-up construction,
\ref{DGsubsubsecinv} 
  can be generalized to types other than nice. For example,
the following data:
\begin{itemize}
\item 
a compact mwc $Z$ and a p-submanifold $Y\subset Z$ (see 
\ref{DGsubsubsecgenbu}, point 4), and
\item
an index family $\Ecal$ on $[Z,Y]$, satisfying \eqref{DGeqcoordind}
\end{itemize}
define the class of functions on $Z\setminus Y$ which have
type $\beta:[Z,Y]\to Z$ with index family $\Ecal$.
Again, this is a piece of discrete data ($\Ecal$) and a piece of
geometric data (which is actually also discrete, since it is natural
to consider diffeomorphic pairs $(Z,Y)$ as equal).

\subsubsection{On 'regimes' and 'matching conditions'}
\label{DGsubsubsecreg}
Characterizations like \eqref{DGeqniceas} are often expressed in terms
of so-called {\em regimes}: In the regime $y/x<C$, $x<C$, $u$ has a
certain asymptotics and in the regime $x/y<C$, $y<C$ it has another.
Of course, these two pieces of data are not independent: Since both 
asymptotics describe the same function, certain relations 
(called {\em matching conditions}) hold between
their coefficients. 

The notion of 'type $\beta$' beautifully and economically
combines regimes and matching conditions
into a single geometric picture; in the case of $[\Rplus^2,0]$ this is the content
of Remark \ref{DGremcoords}.1. 

The correspondence between the 'regime' language and the mwc picture
can be described roughly as follows:
\begin{center}
\begin{tabular}{cp{.36\textwidth}cp{.36\textwidth}}
& regime & $\longleftrightarrow$ & minimal face \\
& matching condition between regimes A,B 
  & $\longleftrightarrow$ 
  & hypersurface containing  the faces corresponding to A,B
\end{tabular}
\end{center}

(A face of a mwc is a non-empty intersection of hypersurfaces, and 
faces are ordered with respect to inclusion.)
\begin{example} In Example \ref{DGexfcns}.5 there are three regimes
(corresponding to each of the three terms being dominant), and
these correspond to the three corners in the mwc used to resolve it
(see Remark \ref{DGrembu}.3).
\end{example}

\subsubsection{How many blow-ups to make?} \label{DGsubsubsecmanyblow}
In a given problem (usually involving differential equations) one often
expects certain type of asymptotic (or singular) behavior for
the solution (for example, from making a model calculation).
This may indicate on which blown-up space one should best consider
 the problem, in order to stay in the realm of nice functions.

However, one has to be careful not to blow up too much:
Although nice functions remain nice after blow-up (Remark \ref{DGrembu}.3a),
differential operators become 'worse'! Thus, one needs to find a 
balance between these forces. We will not address this important
problem any further. This is one of the difficulties in solving the
Main Problem in the Introduction. See the references given there for
solutions in some cases.

%% file: figmwc.pstex_t
\begin{picture}(0,0)%
\epsfig{file=figmwc.pstex}%
\end{picture}%
\setlength{\unitlength}{3947sp}%
\begin{picture}(5609,1634)(2686,-3586)
\put(3076,-3586){\makebox(0,0)[lb]{\smash{All bhs's embedded}}}
\put(6476,-3586){\makebox(0,0)[lb]{\smash{Bhs not embedded}}}
\end{picture}

%% file: figpolar.pstex_t
\begin{picture}(0,0)%
\epsfig{file=figpolar.pstex}%
\end{picture}%
\setlength{\unitlength}{3947sp}%
\begingroup\makeatletter\ifx\SetFigFont\undefined%
\gdef\SetFigFont#1#2#3#4#5{%
  \reset@font\fontsize{#1}{#2pt}%
  \fontfamily{#3}\fontseries{#4}\fontshape{#5}%
  \selectfont}%
\fi\endgroup%
\begin{picture}(5800,2150)(301,-2630)
\put(1751,-1403){\makebox(0,0)[lb]{\smash{\SetFigFont{10}{12.0}{\rmdefault}{\mddefault}{\updefault}$=$}}}
\put(5432,-2700){\makebox(0,0)[lb]{\smash{\SetFigFont{10}{12.0}{\rmdefault}{\mddefault}{\updefault}(c)}}}
\put(836,-645){\makebox(0,0)[lb]{\smash{\SetFigFont{10}{12.0}{\rmdefault}{\mddefault}{\updefault}$\Delta_b$}}}
\put(803,-2073){\makebox(0,0)[lb]{\smash{\SetFigFont{10}{12.0}{\rmdefault}{\mddefault}{\updefault}$\rho$}}}
\put(859,-2700){\makebox(0,0)[lb]{\smash{\SetFigFont{10}{12.0}{\rmdefault}{\mddefault}{\updefault}(a)}}}
\put(1361,-2463){\makebox(0,0)[lb]{\smash{\SetFigFont{10}{12.0}{\rmdefault}{\mddefault}{\updefault}$1$}}}
\put(301,-1069){\makebox(0,0)[lb]{\smash{\SetFigFont{10}{12.0}{\rmdefault}{\mddefault}{\updefault}$\lb$}}}
\put(3034,-2700){\makebox(0,0)[lb]{\smash{\SetFigFont{10}{12.0}{\rmdefault}{\mddefault}{\updefault}(b)}}}
\put(2298,-1660){\makebox(0,0)[lb]{\smash{\SetFigFont{10}{12.0}{\rmdefault}{\mddefault}{\updefault}$B$}}}
\put(3056,-2429){\makebox(0,0)[lb]{\smash{\SetFigFont{10}{12.0}{\rmdefault}{\mddefault}{\updefault}$A$}}}
\put(3814,-2463){\makebox(0,0)[lb]{\smash{\SetFigFont{10}{12.0}{\rmdefault}{\mddefault}{\updefault}$\rb$}}}
\put(2219,-935){\makebox(0,0)[lb]{\smash{\SetFigFont{10}{12.0}{\rmdefault}{\mddefault}{\updefault}$\lb$}}}
\put(2777,-1883){\makebox(0,0)[lb]{\smash{\SetFigFont{10}{12.0}{\rmdefault}{\mddefault}{\updefault}$\ff$}}}
\put(3446,-935){\makebox(0,0)[lb]{\smash{\SetFigFont{10}{12.0}{\rmdefault}{\mddefault}{\updefault}$\Delta_b$}}}
\put(1472,-1069){\makebox(0,0)[lb]{\smash{\SetFigFont{10}{12.0}{\rmdefault}{\mddefault}{\updefault}$\rb$}}}
\put(1138,-2240){\makebox(0,0)[lb]{\smash{\SetFigFont{10}{12.0}{\rmdefault}{\mddefault}{\updefault}$\tau$}}}
\put(914,-2463){\makebox(0,0)[lb]{\smash{\SetFigFont{10}{12.0}{\rmdefault}{\mddefault}{\updefault}$\ff$}}}
\put(2476,-1292){\makebox(0,0)[lb]{\smash{\SetFigFont{10}{12.0}{\rmdefault}{\mddefault}{\updefault}$\eta_2$}}}
\put(3201,-1961){\makebox(0,0)[lb]{\smash{\SetFigFont{10}{12.0}{\rmdefault}{\mddefault}{\updefault}$\eta_1$}}}
\put(2699,-1459){\makebox(0,0)[lb]{\smash{\SetFigFont{10}{12.0}{\rmdefault}{\mddefault}{\updefault}$\xi_2$}}}
\put(3368,-2184){\makebox(0,0)[lb]{\smash{\SetFigFont{10}{12.0}{\rmdefault}{\mddefault}{\updefault}$\xi_1$}}}
\put(357,-2463){\makebox(0,0)[lb]{\smash{\SetFigFont{10}{12.0}{\rmdefault}{\mddefault}{\updefault}$-1$}}}
\put(4372,-1403){\makebox(0,0)[b]{\smash{\SetFigFont{10}{12.0}{\rmdefault}{\mddefault}{\updefault}$\overset{\beta}\longrightarrow$}}}
\put(4729,-991){\makebox(0,0)[b]{\smash{\SetFigFont{10}{12.0}{\rmdefault}{\mddefault}{\updefault}$y$}}}
\put(6100,-2496){\makebox(0,0)[b]{\smash{\SetFigFont{10}{12.0}{\rmdefault}{\mddefault}{\updefault}$x$}}}
\end{picture}

%% file: DGanalysis.tex
\section{Analysis} \label{DGsecana}
In this section we discuss two of the basic processes of analysis:
pull-back and push-forward, and how they affect asymptotic behavior of smooth
functions (as discussed in the previous section) and conormal
distributions, which we also introduce.

What are pull-back and push-forward, and why are they important?
Pull-back is composition, push-forward is integration. They are important
since they may be used as building blocks for other operations.
This allows to carry out recurring ugly calculations (e.g.\ those
involving Fourier transform) once in the proof of theorems about pull-back
and push-forward, and then never look at them again.
Let us illustrate this in two simple but central examples:%
\footnote{
In the examples, $\R$ may be replaced by any manifold, equipped with a 
fixed density. 
For the moment we naively neglect the distinction between functions,
distributions and the respective densities; also, we neglect such 
tedious matters as integrability.
}
\begin{examples}\label{DGexpfpb}
\
\begin{description}
\item [Applying an operator to a function]
Let $\pi_1,\pi_2:\R^2\to\R$ be the projections onto the first and second
coordinate. If $v$ is a function on $\R$ then its pull-back
$\pi_2^*v=v\circ\pi_2$ is the function $(x,y)\mapsto v(y)$ on $\R^2$. 
If $u$ is a function on $\R^2$ then its push-forward $\pi_{1*}u$ is the
function 
\begin{equation}\label{DGeqpf1}
\pi_{1*}u(x) = \int u(x,y) dy.
\end{equation}
If $A$ is an operator, acting on functions $v$ on $\R$, with
integral kernel $A(x,y)$ then
\begin{equation}\label{DGeqaction}
 (Av)(x)  = \int A(x,y)v(y)\, dy = \pi_{1*}(A\cdot\pi_2^*v).
\end{equation}
Though this may look like an exercise in
formal nonsense, it shows
that mapping properties of $A$ may be read off from the structure 
(e.g.\ asymptotic type) of
the function (distribution) $A$, if one understands how such structure
is affected by pull-back and push-forward.%
\footnote{
Also, one needs to understand how structure is affected by
 multiplication. This is trivial for nice functions, geometrically non-trivial
for functions with different asymptotic types, and analytically non-trivial
for distributions. See Subsection 
\ref{DGsubsecdist} for the latter case.}

\item[Composition of operators]
If $A,B$ are operators, acting on functions on $\R$, with
 integral kernels $A(x,y),B(y,z)$, then 
$C=A\circ B$ has integral kernel $C(x,z)=\int A(x,y)B(y,z)\,dy$, i.e.
\begin{equation}\label{DGeqcomp}
C=\pi_{2*}(\pi_3^*A\cdot\pi_1^*B)
\end{equation}
where $\pi_1,\pi_2,\pi_3:\R^3\to\R^2$ are the projections leaving out
the first, second and third variable, respectively.
Again, understanding how pull-back, push-forward  and product
affect the structure
of distributions allows to predict, for example, whether a class of operators
with a given structure of its kernels is closed under composition.
\end{description}
\end{examples}

Another place where one needs to understand the behavior of distributions
under push-forward is in the 'specializations' mentioned in 
Footnote \ref{DGfnspec} in the Introduction, since many of them are obtained
from the full kernel of $Q$ by integration (i.e.\ push-forward).

The  maps used for pull-back and push-forward in the examples
are rather trivial projections, 
so it's legitimate to ask: Why be so formal, why
not talk simply of 'integration in $y$' instead of 'push-forward by $\pi_1$'?
The answer is given by:
\begin{principle}
The push-forward of a complicated function by a simple map should be 
analyzed by rewriting it as push-forward of a simple function by a complicated
map.
\end{principle}
The point is that the 'complication' of the map lies mainly in its global
geometry, so by a partition of unity the problem can be reduced to the sum of
relatively simple local problems. (In contrast, the 'complication' of the
function is local, typically.)
Melrose's Push-Forward Theorem gives the result of this analysis, for the
case of smooth functions.
In Subsection \ref{DGsubsecpf} we discuss all these matters, starting from
an example.
We also sketch the idea of a proof of the Push-Forward Theorem in the
special case that the target space is $\Rplus$.

When integrating one needs measures. Therefore, push-forward
is best defined as acting on measures (or densities) rather than
functions. The
push-forward of a smooth density may be not smooth, and (what's equivalent)
the pull-back of a distribution is not always defined.
For the reader unfamiliar with these matters, 
we give the precise definitions and a short discussion
of pull-back and push-forward, and how they act on smooth functions,
distributions and smooth and distributional densities, in the Appendix.%
\footnote{\label{DGfndens}%
The reader who prefers to neglect the distinction between
functions and densities is invited to do so, but will probably begin
to acknowledge their usefulness when making computations herself.
}

In Subsection \ref{DGsubsecpb} we state Melrose's
Pull-Back Theorem, which tells how pull-back by a \Mb-map affects nice
functions. This  is rather trivial in comparison to the Push-Forward Theorem.

Finally, in Subsection \ref{DGsubsecdist} we introduce conormal 
distributions and discuss how pull-back,
 push-forward and multiplication affect them. As an
illustration, we define pseudodifferential operators and 
study their composition.

%%%%%%%%%%%%%%%%%%%%%%%%%%%%%%%%%%%%%%%%%%%%%%%%%%
%%%%%%%%%%%%%%%%%%%%%%%%%%%%%%%%%%% Push-Forward Subsec.
%%%%%%%%%%%%%%%%%%%%%%%%%%%%%%%%%%%%%%%%%%%%%%%%%%

\subsection{Push-forward and asymptotic type}\label{DGsubsecpf}
We begin by analyzing a few examples of push-forward under the
projection $\R_+^2\to\R_+,(x,y)\mapsto x$.
In other words, we set
\begin{equation}\label{DGeqpft}
\utilde(x) = \int_0^\infty u(x,y)\,dy,\quad x>0.
\end{equation}
Assuming that $u$ is smooth in $(0,\infty)^2$ and $\supp u$ is bounded,
we ask how the behavior of $u$ near the boundary of $\R_+^2$ affects the
behavior of $\utilde$ near 0.
\begin{examples}\label{DGexpft}
\
\begin{enumerate}
\item
If $u$ is smooth up to the boundary then  so is $\utilde$ (by 
first-year analysis).
More generally (and just as easy),
$$ u \text{ nice with index family }(E,F) 
   \Rightarrow \utilde\text{ nice with index set }F$$
if the integral \eqref{DGeqpft} exists at all, i.e.\ if
\begin{equation}\label{DGeqint1}
\Re z>-1 \text{ for } (z,p)\in E.
\end{equation}
\item
If $u(x,y)=y^{-1}v(x/y,y)$ with $v$ smooth on $\R^2_+$ and compactly supported  then 
\begin{equation}\label{DGeqhyperbel}
 \utilde(x)=\int_0^\infty v(\frac{x}y,y)\,\frac{dy}y
\sim \sum_{i=0}^\infty (a_ix^i + b_ix^i\log x)\quad\text{ as }x\to 0,
\end{equation}
i.e., $\utilde$ is nice, but not smooth (the index set is
 $\N_0\times\{0,1\}$).%
\footnote{
Proof: Taylor expand $v(\xi,\eta)$ at $\xi=0$ (for each fixed $\eta$), then
Taylor expand each coefficient and the remainder at $\eta=0$ to obtain, for
any $N$,
\begin{equation}\label{DGeqexpansion}
v(\xi,\eta) = \sum_{\alpha=0}^{N-1} \xi^\alpha \eta^Na_\alpha(\eta)
+ \sum_{\beta=0}^{N-1} \eta^\beta \xi^N b_\beta(\xi)
+ \sum_{\alpha,\beta=0}^{N-1} c_{\alpha,\beta} \xi^\alpha \eta^\beta
+ \xi^N \eta^N r(\xi,\eta)
\end{equation}
with $a_\alpha,b_\beta,r$ smooth up to the boundary. Assume $\supp v
\subset [0,C]^2$. Then in the integral \eqref{DGeqhyperbel}
one may replace $\int_0^\infty$ by $\int_{x/C}^1$. To obtain the asymptotics,
simply integrate \eqref{DGeqexpansion} term by term,
using the substitution $z=x/y$ in the second sum.

Note that the $\log$-terms only come from the terms $\alpha=\beta$ in 
the third sum.
}
More generally, if $v$ has index sets $(E,F)$ then $\utilde$ is
nice with index set%
\footnote{
Same proof, after the (non-trivial) analysis argument that our
definition \eqref{DGdefnice} of niceness
implies an expansion like \eqref{DGeqexpansion}. Alternatively, one may define
niceness by this expansion.
}
\footnote{\label{DGfndyy}%
Why did we write the integral \eqref{DGeqhyperbel} with
$dy/y$ instead of simply $dy$? Since then the result 
\eqref{DGeqextunion} is beautifully symmetric!
Cf. '\Mb-densities' below.
}
\begin{equation} \label{DGeqextunion}
E\extunion F := E\cup F \cup \{(z,p'+p''+1):\, (z,p')\in E, \, (z,p'')\in F\}.
\end{equation}

\item
For $u(x,y)=\sqrt{x^2+y^2}$ explicit integration shows (restricting
to $y\leq 1$ for integrability -- this does not affect the essential point)
\begin{equation}\label{DGeqlog}
 \utilde(x)= \text{ (smooth near zero) } - \frac12 x^2\log x.
\end{equation} 
Thus, $\utilde$ is nice, and again a logarithm appears. 
\end{enumerate}
\end{examples}
The common feature of Examples 2 and 3 is that $u$ has asymptotic type 
$\beta$, where $\beta$ is the blow-up of 0 in $\Rplus^2$.
\medskip

{\em Claim: This already suffices to explain the similarity of the results
\eqref{DGeqhyperbel} and \eqref{DGeqlog}.}

\proof
We first show this by a simple calculation and then explain how it
may be seen directly by 'looking at pictures'.

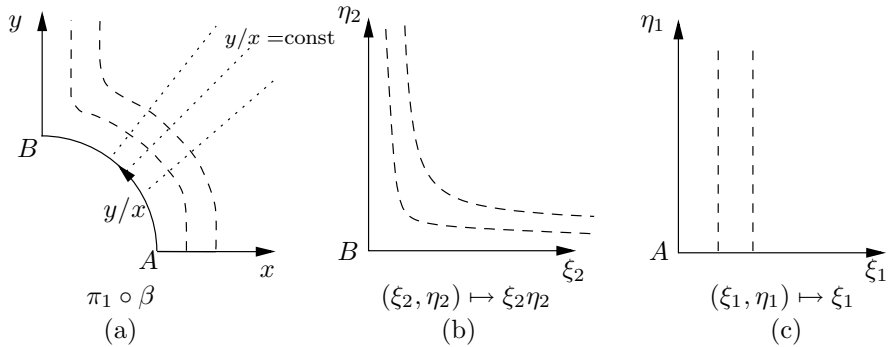
\begin{figure}[htbp]
\input{figpf.pstex_t}
\caption{Level lines for push-forward under $[\Rplus^2,(0,0)]\to\Rplus$}
\label{DGfigpf}
\end{figure}

{\em Calculation: } Consider any $u$ of type $\beta$,  so that $w=\beta^*u$ 
is nice on $W=[\Rplus^2,0]$, and assume $w$ has no
logarithmic terms in its expansions.
We split up the integral
$$
\int_0^\infty\!\! u(x,y)\,dy = A+B
$$
'smoothly near $y=x$'. That is, with any cut-off function $\phi\in 
C_0^\infty(\Rplus)$ which equals one near 0, and with $\psi=1-\phi$, we set
\begin{align*}
A=\int_0^\infty\!\! u(x,y)\,\phi(y/x)\,dy &
    = \int_0^\infty\!\! x u(x,x\eta_1)\,\phi(\eta_1)\,d\eta_1 
   = \int_0^\infty\!\! x w_1(x,\eta_1)\,\phi(\eta_1)\,d\eta_1\\
B=\int_0^\infty\!\! u(x,y)\,\psi(y/x)\,dy &
    = \int_0^\infty\!\! w_2(x/y,y)\,\psi(y/x)y\,\frac{dy}y.
\\
\end{align*}
Here, $w_1$ is just $w$ expressed in projective coordinates near the point
$A$ in Figure \ref{DGfigpf}(a)
(i.e.\ $w_1(\xi_1,\eta_1)=u(\xi_1,\xi_1\eta_1)$ or  $w_1 = \beta_1^*u$
with $\beta_1$ from (\ref{DGeqproj}c));
the integral $A$ is like Example 1 (with $u(\xi,\eta)=\xi w_1(\xi,\eta)
\phi(\eta)$).
Similarly, $w_2$ is just $w$ expressed in projective coordinates
near $B$ in Figure \ref{DGfigpf}(a) (i.e.\ $w_2(\xi_2,\eta_2)=u(\xi_2\eta_2,\eta_2)$
or $w_2=\beta_2^*u$); the integral $B$ is like Example 2 
(with $v(\xi,\eta) = w_2(\xi,\eta)\psi(\xi^{-1})\eta$).

Since by assumption $w_{1/2}$ have no $\log$'s in their expansions
and in Example 1 no logarithms are created, we conclude:
The log terms in Examples 2 and 3 are of the same nature.%
\footnote{
But we also see that for general $\beta$-singular $u$ infinitely many
$\log$ terms will appear. For $u(x,y)=\sqrt{x^2+y^2}$ only
one $\log$-term appears (see \eqref{DGeqlog}); this is due to the fact
that in $\sqrt{x^2+y^2}=y\sqrt{1+(x/y)^2}=x\xi_2^{-1}\sqrt{1+\xi_2^2}$
only one power of $x$ occurs. 
Such fine points are lost under (regular) coordinate changes and therefore
invisible in the geometric setup of the Push-Forward Theorem.
}

{\em Pictures: }
We now show how the same result can be 'seen' geometrically.
Since $\beta^*u=w$ and $\beta$ is a diffeomorphism in the interior, we have
$u=\beta_*w$, so
\begin{equation} \label{DGeqpfcomp}
\utilde = \pi_{1*}u = \pi_{1*}\beta_*w = f_*w
\end{equation}
with $f=\pi_1\circ\beta:W\to\Rplus$.
This says simply that $\utilde(x)$ equals  the integral of $w$ over
the fiber $f^{-1}(x)$ for each $x$ (let's postpone the question of
measures for a moment). This is clear since the values of $w$ on
$f^{-1}(x)$ are precisely the values of $u$ on $\pi_1^{-1}(x)$, which are
integrated to obtain $\utilde(x)$.

Some fibers of $f$ are shown in Figure \ref{DGfigpf}(a), some of $\pi_1$ 
in Figure \ref{DGfigpf}(c)
and some of $g(\xi,\eta)=\xi \eta$ in Figure \ref{DGfigpf}(b). 
Pictorially, we see:%
\footnote{
This can be made precise by expressing $f$ in projective
 local coordinates \eqref{DGeqproj}:
\begin{itemize}
\item On $W\setminus\lb$ ('near A') $f$ is expressed
as $f_1(\xi_1,\eta_1)=\xi_1$ (using $\beta_1$ in (\ref{DGeqproj}c)),
i.e.\ $f_1=\pi_1$,
\item on $W\setminus\rb$ ('near B') as $f_2(\xi_2,\eta_2)=\xi_2\eta_2$
(using $\beta_2$), i.e. $f_2=g$.
\end{itemize}
}
\begin{itemize}
\item Near A, Figure \ref{DGfigpf}(a) looks like Figure \ref{DGfigpf}(c),
\item near B, Figure \ref{DGfigpf}(a) looks like Figure \ref{DGfigpf}(b).
\end{itemize}

Therefore, push-forward of $w$ by $f$ is the sum of push-forward
(of $w$ near $A$)
by $\pi_1$ and push-forward (of $w$ near $B$) by $g$, and this was
precisely the calculation above. This also explains why the cut-off
had to be chosen as a smooth function of $y/x$, see Figure \ref{DGfigpf}(a).

In summary, we may say that the $\log$ terms in Examples 2 and 3 arise
from the fact that the fibers of $g$ and of $f=\pi_1\circ\beta$ approach the
corner as in Figure \ref{DGfigpf}(b) for $x\to 0$.
\qed
\medskip

In Melrose's Push-Forward Theorem these considerations are 
generalized to any \Mb-map $f:W\to Z$:
Under certain conditions on $f$,
it says that the push-forward of a nice density $\mu$
on $W$ is a nice density on $Z$, and computes the index sets of the
latter from the index sets of the former and the 'boundary geometry'
of $f$.

The conditions on $f$ are best understood if we consider the special case
$Z=\Rplus$ first. Before we can state them, we need some definitions.

From now on, we assume that all index sets satisfy 
$(z,p)\in E\implies (z+1,p)\in E$, the condition ensuring coordinate 
invariance of niceness, and that bhs's are embedded and connected
(see Definition~\ref{DGdefmwc} and Footnote~\ref{DGfnmwc}).

\bigskip

{\em Densities on manifolds with corners.}
If $W$ is a mwc, then a density on $W$ is, by definition, a density on the
interior $\Wint$ (concerning densities see the Appendix and 
Footnote~\ref{DGfndens}). The notion of niceness carries over to densities
immediately, e.g.\ on $\Rplus^2$:
\begin{definition}
A density $\mu=u\, dxdy$ on $\Rplus^2$ is {\em nice} with index sets $E,F$ 
if $u$ is nice with index sets $E,F$.%
\footnote{
As usual, this should be checked for coordinate independence.
But only under coordinate changes $(x,y)\mapsto (\xtilde,\ytilde)$ for
which $\xtilde,\ytilde$ are still bdf's of the coordinate axes!
Cf.\ \ref{DGsubsubsecinv}.
}
\end{definition}

The following is a slight variant in book-keeping, which makes lots of things
more transparent%
\footnote{\label{DGfnbdens}%
Examples:
\begin{enumerate}
\item
$\mu$ locally integrable $\Longleftrightarrow$ $\Re z>0$ whenever
$(z,p)\in E\cup F$
(rather than $-1$).
\item
The transformation under projective coordinates becomes especially simple:
Say $\xi_1=x,\eta_1=y/x$, then
\begin{equation}\label{DGeqbdenstrf}
 \frac{dx}x\frac{dy}{y} = \frac{d\xi_1}{\xi_1}\frac{d\eta_1}{\eta_1}.
\end{equation}
\item See Footnote \ref{DGfndyy} after Example \ref{DGexpft}.2.
\end{enumerate}
}
(though it may seem artificial to the uninitiated):
\begin{definition} \label{DGdefbdens}
A {\em \Mb-density}
 on $\Rplus^2$ is just a density, except that we write it as
$\mu=u(x,y) \frac{dx}x\frac{dy}y$ instead. When talking about smoothness or
the index family of $\mu$ then we mean smoothness or the index family
 of $u$ in
such a representation.
\end{definition}

Of course, a \Mb-density on $\Rplus\times\R$ is of the form
$u(x,y)\,\frac{dx}x dy$. That is, the $\frac{dx}x$ factor only occurs
in the variables $x$ defining some bhs. It is easy to see that the index 
family of a \Mb-density is well-defined on any mwc.

\bigskip

{\em Boundary geometry of a \Mb-map $f:W\to \Rplus$.}
\begin{definition} \label{DGdefexp1}
Let $f:W\to\Rplus$ be a \Mb-map. For any bhs $G$ of $W$ define $e_f(G)$ to be
the order of vanishing of $f$ at $G$.
\end{definition}

In other words,
in the local Definition \ref{DGdefbmap} with  $w\in G$, we set
$e_f(G)=\alpha_{1j_0}$ if $f(w)=0$ and $x_{j_0}$ is a bdf for $G$, and
$e_f(G)=0$ if $f(w)\not=0$. This is clearly locally constant and therefore
constant on $G$ by connectedness, so $e_f(G)$ is well-defined.
Note that
\begin{equation} \label{DGequnion}
f^{-1}(0) = \bigcup \{G: e_f(G) > 0\}.
\end{equation}

\begin{theorem}[Push-Forward Theorem, special case $Z=\Rplus$]
\label{DGthpft1}
Let $W$ be a manifold with corners and $f:W\to\Rplus$ a \Mb-map which is
a fibration over $(0,\infty)$.%
\footnote{
I.e.\ $f:f^{-1}((0,\infty))\to(0,\infty)$ is a fibration in the sense
of Footnote \ref{DGfnfibr} in the Appendix, 
except that $L$ is allowed to be a mwc.}
Let $\Ecal$ be an index family for $W$.
Assume that $f,\Ecal$ satisfy 
the integrability condition \eqref{DGeqint} below.
 
If $\mu$ is a compactly supported \Mb-density on $W$, nice
with index family $\Ecal$, then $f_*\mu$ is a \Mb-density on $\Rplus$,
nice with index family $f_\#\Ecal$ (defined in \eqref{DGeqindexset} below).
\end{theorem}
The integrability condition is:
\begin{equation}\label{DGeqint}
\inf \Ecal(G)>0  \text{ whenever } e_f(G) = 0
\end{equation}
where for any index set $E$
\begin{equation}\label{DGeqinfdef}
\inf E := \inf \{\Re z: (z,p)\in E\}
\end{equation}
(which is actually a minimum).%
\footnote{
Geometrically,  $e_f(G)=0$ means that $f>0$ on $G$, so 
 the fibers $f^{-1}(x)$, $x>0$, 
 will hit only these $G$, and actually transversally
as in Figure \ref{DGfigpf}(a) at the $x$-axis.
So \eqref{DGeqint} generalizes \eqref{DGeqint1}
and comes from the fact that 
$\int_0^1 x^z \,\frac{dx}x$ exists iff $\Re z>0$. }
To define $f_\#\Ecal$, associate to every  face $F$ 
(i.e.\ non-empty intersection of boundary hypersurfaces) of $W$ 
the index set
$$\tilde{\Ecal}(F) = \Extunion_G \left\{(\frac{z}{e_f(G)},p):\,(z,p)\in
       \Ecal(G)\right\}$$
where the extended union (defined in \eqref{DGeqextunion}) is 
over all bhs's $G$ containing $F$ and having $e_f(G)>0$.
Then, define
\begin{equation}\label{DGeqindexset}
f_\#\Ecal = \bigcup_F \tilde{\Ecal}(F).
\end{equation}

\begin{remarks}
\
\begin{enumerate}
\item
$f$ needs to be a fibration in the interior to ensure that 
$f_*\mu$ is smooth in the interior.
\item
The definition of $f_\#\Ecal$ given above is a little more precise
than the one given  in \cite{DGMel:CCDMWC} (which may yield an index set that
is 'too big'). But the Push-Forward Theorem with this 'smaller' $f_\#\Ecal$
follows directly from Melrose's by introducing a suitable partition of unity.%
\footnote{
Clearly, in \eqref{DGeqindexset} it is enough to take the union over all 
{\em minimal} faces (with respect to inclusion), for example the corners 
$A$, $B$ in Figure \ref{DGfigpf}(a). Thus, any 'regime' on $W$ (see
\ref{DGsubsubsecreg}) contributes some asymptotic terms.}
\item 
The {\em proof} was essentially done above: Localize
as in the discussion of Example \ref{DGexpft}.3, this reduces to the
cases of Examples \ref{DGexpft}.1/2 (modulo 
replacing $x,y$ by powers $x^\nu,y^\mu$ with $\nu,\mu>0$
determined by the $e_f(G)$, and modulo straight-forward 
generalization to higher dimensions.)
\item 
See the article \cite{DGGriGru:SALPFT}
 \inthisbook\ for a discussion of the relation of
the Push-Forward Theorem (with $Z=\Rplus$)
and the 'Singular Asymptotics Lemma' by Br\"uning and Seeley 
(\cite{DGBruSee:RSA}).
\end{enumerate}
\end{remarks}

\medskip

{\em Push-Forward Theorem with  general target space.}
Here, some additional assumptions on the map $f$ are needed.
Before we can state these, we need to look a little closer at the 
geometry of \Mb-maps:
\medskip

{\em Boundary geometry of \Mb-maps.}
By definition, $f:W\to Z$ is a \Mb-map iff
\begin{equation} \label{DGeqfH}
f_H:=\rho_H\circ f:W\to\Rplus
\end{equation}
is a \Mb-map for all bhs's $H$  of $Z$, 
and bdf's $\rho_H$ of $H$.
So we can define:
\begin{definition} \label{DGdefexp2}
The {\em exponent matrix} of a \Mb-map $f:W\to Z$ is the set of integers
$$ e_f(G,H) = e_{f_H}(G),\quad G\text{ bhs of }W,\, H\text{ bhs of }Z.$$
\end{definition}
Thus, $e_f(G,H)\not=0$ iff $f(G)\subset H$, and in this case if $p\in W$ 
has small distance $\eps$ from $G$ and distance $\geq$ const $>0$ from all
other bhs's of $W$, then $f(p)$ has distance of order $\eps^{e_f(G,H)}$ 
from $H$ (say in Euclidean metric for any local coordinate systems based 
at points of $G$ and $H$).

Referring to the Definition \ref{DGdefbmap} of \Mb-maps,
we have $e_f(G,H)=\alpha_{ij}$
in \eqref{DGeqbmap} if $G=\{x_j=0\}$ and $H=\{x_i'=0\}$ locally.

Recall that a {\em face} of a mwc $W$ 
is a non-empty intersection of boundary hypersurfaces,
or $W$ itself. Each face is a mwc.
A \Mb-map $f$ induces a map
 $$\fbar:\text{faces of }W\to\text{ faces of }Z$$
characterized by
\begin{equation} \label{DGeqftilde}
 x\in \interior{F} \implies f(x)\in \interior{(\fbar(F))}. 
\end{equation}
Alternatively, $\fbar(F)=$ the intersection of the bhs's $H$ of $Z$ satisfying
$f(F)\subset H$.

In summary, the 'combinatorics' of a \Mb-map $f$ can be described either
by giving the pairs $(G,H)$ with $f(G)\subset H$, or equivalently by
the map $\fbar$, or (a little more refined) by the matrix $e_f$.
\begin{definition}\label{DGdefbfibr}
A \Mb-map $f:W\to Z$ is a {\em \Mb-fibration} if for each face $F$ of $W$,
\begin{itemize}
\item[(a)]
$\codim \fbar(F)\leq \codim F$ 
(it is enough to require this of bhs's $F$), and
\item[(b)]
$f$ is a fibration $\interior{F}\to \interior{(\fbar(F))}$. %
\footnote{
Melrose's definition in \cite{DGMel:CCDMWC} looks different, but is 
equivalent: As is easily seen, condition (a) is equivalent to his
'\Mb-normality' (we assume all \Mb-maps to be interior), and, assuming (a),
condition (b) is equivalent to his '\Mb-submersion' condition, at least for
proper maps (for which submersion $\iff$ fibration), which is all that 
matters anyway.}
\end{itemize}
\end{definition}
For example, if $Z=\R_+$ then $f$ is a \Mb-fibration iff $f$ is a 
fibration over $(0,\infty)$ (here (a) is empty). The polar-coordinate
map is not a \Mb-fibration (since $\ff$ gets mapped to a codimension two
face), nor is any other non-trivial blow-down map.
An important example of a \Mb-fibration where condition (a) is non-empty
is given by the projection from the 'triple \Mb-space' in \eqref{DGeqbfibrex}.

For a \Mb-map $f:W\to Z$ and an index family $\Ecal$ on $W$, define
the index family $f_\#\Ecal$ on $Z$ by
\begin{equation}\label{DGeqindexfam}
f_\#\Ecal(H)= (f_H)_\#\Ecal,
\end{equation}
with the right hand side defined by \eqref{DGeqindexset} and
\eqref{DGeqfH}.
The integrability condition is:
\begin{equation} \label{DGeqint2}
\inf \Ecal(G)>0 \quad\text{whenever}\quad e_f(G,H)=0\ \;\;\forall\ H.
\end{equation}
(The latter condition means that $f(G)\not\subset\partial Y$.)
\begin{theorem}[Push-Forward Theorem]
\label{DGthpft}
Let $W,Z$ be manifolds with corners and $f:W\to Z$ a \Mb-fibration.
Let $\Ecal$ be an index family for $W$, and assume \eqref{DGeqint2}.

If $\mu$ is any compactly supported \Mb-density $\mu$ on $W$, nice with
index family $\Ecal$, then
the push-forward $f_*\mu$ is a \Mb-density on $Z$, nice with
index family $f_\#\Ecal$ (defined in \eqref{DGeqindexfam}).
\end{theorem}
See \cite{DGMel:CCDMWC} for a proof.
\begin{remarks}\label{DGrempft}
\
\begin{enumerate}
\item 
{\em Why the \Mb-fibration conditions are needed:}
For example, (a) is violated for the polar-coordinate map, and this map
does not preserve niceness (this was the reason for doing blow-ups
in the first place!).
Also, if
$f_*\mu$ is to be nice then the expansion coefficients in
the asymptotics at any face should be smooth in the interior of the face, so
one should require a fibration condition here, which explains (b).
(Thus, (a) ensures good behavior of $f_*\mu$ 
when approaching the boundary, while
(b) does so {\em in} the boundary and locally in the interior.)
\item 
{\em On determining the asymptotic type of $g_*\mu$ from
the asymptotic type of $\mu$}, for a map $g:Z\to Z'$ between mwc's:

For $Z'=\Rplus$ the answer is given essentially by Theorem \ref{DGthpft1}:
If $g:Z\to\Rplus$ is a fibration in the interior and $\mu$ has
type $\beta:W\to Z$ then applying the theorem 
to $f=g\circ\beta$ shows that $g_*\mu$ is nice on $\Rplus$.
Compare \eqref{DGeqpfcomp} where $g=\pi_1$. 
Clearly, this only works if $(\beta^{-1})_*$
maps nice densities on $Z$ to nice densities on $W$;  for
blow-ups $\beta$ this is clearly true, see \eqref{DGeqbdenstrf}
in Footnote \ref{DGfnbdens}.

For general $Z'$ the problem is: Given $g:Z\to Z'$ and a blow-up
$\beta:W\to Z$, find a blow-up $\beta':W'\to Z'$ such that densities
of type $\beta$ are pushed forward to densities of type $\beta'$.
By the Push-Forward Theorem, this would be satisfied if 
$\gtilde=(\beta')^{-1}\circ g\circ\beta:W\to W'$
was a \Mb-fibration.
\begin{equation} \label{DGeqcd}
\begin{CD}
W @>\gtilde>>W'\\
@V\beta VV @VV\beta' V \\
Z @>g>> Z'
\end{CD}
\end{equation}
 Note that even for $\gtilde$ to be a well-defined map
implies restrictions on $\beta'$. The problem 
when this is possible, and how to find $\beta'$, seems to be difficult.
\item
The support condition on $\mu$ in Theorem \ref{DGthpft} merely excludes
problems of non-integrability at infinity. Clearly, it could be weakened
to: $f$ is proper on $\supp\mu$.
(We need this extension when we discuss \PDO s.)
\end{enumerate}
\end{remarks}

\subsection{Pull-back and asymptotic type} \label{DGsubsecpb}
Though just as important as the Push-Forward Theorem, this is child's play
in comparison. The main point was already mentioned in Remark
\ref{DGrembmap}.1. Doing the book-keeping on the index sets easily yields:

\begin{theorem}[Pull-Back Theorem] \label{DGthpbt}
Let $f:W\to Z$ be a \Mb-map.
Then, for any function $v$ on $Z$ which is nice with 
index family $\Fcal$, the pull-back $f^*v$ is a nice function
on $W$ with index family $f^\#\Fcal$ defined by: If $G$ is a bhs of $W$ then
\begin{multline*}
f^\#\Fcal(G) = \Big\{(q+\sum_H e_f(G,H)z_H, \sum_H p_H):\,q\in\N_0
         \text{ and}\\
\text{ for each bhs $H$ of $Z$:}
\begin{cases}
(z_H,p_H)\in \Fcal(H) &\text{if } e_f(G,H)\not=0,\\
(z_H,p_H)=(0,0) &\text{if } e_f(G,H)=0.
\end{cases}
\Big\}.
\end{multline*}
\end{theorem}

\begin{remarks} \label{DGrempbt}
\
\begin{enumerate}
\item
{\em On determining the asymptotic type of $g^*v$ from the
asymptotic type of $v$.}
Here one is given $g$ and a blow-up $\beta'$  in diagram \eqref{DGeqcd},
and needs to find a blow-up $\beta$
such that $g^*v$ has type
$\beta$ whenever $v$ has type $\beta'$. By the Pull-Back Theorem,
this is satisfied if $\gtilde$ is
well-defined and a \Mb-map (and $g$ is surjective).
Note that $(\beta')^{-1}\circ g$ is usually only defined on the interior since
$\beta'$ is not a diffeomorphism on the boundary, so $W$ has to be chosen
'big' enough so that $\gtilde$ may extend continuously from the interior
to all of $W$.
\item
{\em The triple \Mb-space.}
As example consider the case relevant for composition 
in the \Mb-\PDO\ calculus:
$g=\pi_1:\Rplus^3\to\Rplus^2$ is the projection $\pi_1(x_1,x_2,x_3)=(x_2,x_3)$.
The solution is easy: If $\beta':W'\to\Rplus^2$ is any blow-up then let
$W=\Rplus\times W'$, $\beta=\id_{\Rplus}\times\beta':\Rplus\times W'\to
\Rplus\times\Rplus^2$. In the case relevant for us, where $W'=[\Rplus^2,0]$
is just the blow-up of zero, $W$ is the blow-up of the $x_1$-axis.

However, in the composition problem $W$ and $\beta$ need to work for
several maps $g$  {\em simultaneously}, and this makes
the problem more interesting. 
Let $\pi_i:\Rplus^3\to\Rplus^2$ be the projection that forgets the
$i$'th coordinate, for $i=1,2,3$.
\medskip
\begin{quote}
{\bf Problem:}
Find a blow-up $\beta:W\to\Rplus^3$ such that whenever $v$ has type
$\beta':[\Rplus^2,0]\to\Rplus^2$ then $\pi_i^*v$ has type $\beta$ for
$i=1,2,3$.
\end{quote}
\medskip

In other words, $\pitilde_i=(\beta')^{-1}\circ\pi_i\circ\beta:W\to[\Rplus^2,0]$
must be a \Mb-map for $i=1,2,3$.
It is clear that at least all three coordinate axes must be blown up.
The most naive thing to try is to blow up one axis (say the $x_1$-axis)
and then (the preimages of)
the other two. However, it is easily seen
that $\pi_2$ and $\pi_3$ are still not well-defined on the resulting space.

But there is a beautiful solution which even preserves the symmetry:
First, blow up zero in $\Rplus^3$. Then, blow up the preimages of
the three coordinate axes (in any order, since they are separated now!).
the result is called triple \Mb-space $X^3_b$
and shown in Figure \ref{DGfigtriple}.
\begin{figure}[htbp]
\input{figtriple.pstex_t}
\caption{The triple \Mb-space (and projection $\pitilde_3$)}
\label{DGfigtriple}
\end{figure}
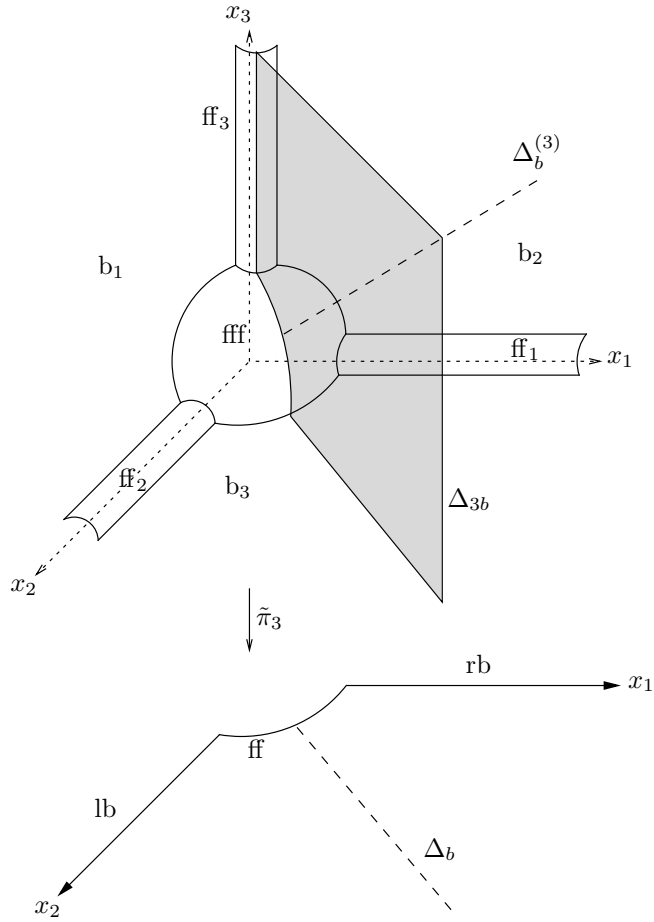

Let us convince ourselves pictorially that the maps
\begin{equation} \label{DGeqbfibrex}
\pitilde_i:X_b^3\to X_b^2:= [\Rplus^2,(0,0)]
\end{equation}
are well-defined and \Mb-fibrations.
By symmetry, it is enough to consider $\pitilde_3$.
It is well-defined since the $x_3$-axis was blown up.
Denote the bhs's of $X_b^3$ by $\bface_1,\bface_2,\bface_3$ (the 'old'
bhs's from $\Rplus^3$), $\ff_1,\ff_2,\ff_3$ (the front faces of the axis
blow-ups), and $\fff$ (the front face of the point blow-up) as in
Figure \ref{DGfigtriple}.
The bhs's are mapped as follows:
\begin{align} \label{DGeqbhs}
\begin{split}
 \ff_2, \bface_1& \mapsto \lb\\
 \ff_1, \bface_2&\mapsto \rb\\
 \fff,\ff_3&\mapsto\ff\\
\bface_3&\mapsto X^2_b,
\end{split}
\end{align}
and all these maps are onto. Also, $\interior{(X_b^3)}\to\interior{(X_b^2)}$.
Therefore, the preimage of each bhs is a union of bhs's, which almost
shows that $\pitilde_3$ is a \Mb-map (see  Remark \ref{DGrembmap}.1; of course
one may check the full condition \eqref{DGeqbmap} by direct calculation).
Also, \eqref{DGeqbhs} defines the map $\pitildebar_3$
(see \ref{DGeqftilde}) on bhs's and
 this determines $\pitildebar_3$ for all faces.
Since all faces on the right of \eqref{DGeqbhs} have codimension at most one,
condition (a) in the Definition \ref{DGdefbfibr} of a \Mb-fibration
 is satisfied. 

Condition (b) is easily checked for each face: For example, 
$\interior{\bface_1}\to\interior{\lb}$ is basically the same map
as $[\Rplus^2,(0,0)]\to\Rplus$ from \ref{DGfigpf}, and so is
$\interior{\fff}\to\interior{\ff}$ near the boundary. All maps from
codimension two faces are either diffeomorphisms or constant, so they
are fibrations trivially.
\end{enumerate}
\end{remarks}

%%%%%%%%%%%%%%%%%%%%%%%%%%%%%%%%%%%%%%%%%%%%%%%%%%%%%%%%%%%%%%%
%%%%%%%%%%%%%%%%%%%%%%%%%%%%%%%%%%%% Distributions subsec
%%%%%%%%%%%%%%%%%%%%%%%%%%%%%%%%%%%%%%%%%%%%%%%%%%%%%%%%%%%%%%%

\subsection{Distributions} \label{DGsubsecdist}
So far, all singular behavior occurred at the boundary.
Now we turn to the description of singularities in the interior of
a mwc $Z$. This means talking about distributions%
\footnote{
There are also distributions whose singular support is 
contained in the boundary. We will not discuss them here (although
they are not really more difficult). 
See \Mnotes, for example.}.
In many situations only a very special class of distributions occurs, 
the 'step 1 polyhomogeneous conormal' (here called 'conormal') ones%
\footnote{
This is not the most general kind of what's usually called
conormal distributions, but they are easy
to define and sufficient for many purposes.}.
They are smooth outside a
submanifold, and at the submanifold have a special explicitly
describable kind of singular behavior, which is in some sense similar
to the behavior of a  nice function at the boundary.

In the case of manifolds most of
this material is quite standard (see e.g.\ \cite{DGHor:ALPDOIII},
\cite{DGSim:PDO});
we will briefly recall
the definition, give some examples and state the push-forward and pull-back 
theorems. As an illustration, we
use this to show that the set of (properly supported) classical
pseudodifferential operators on $\R$ is closed under composition.
The extension of the definition and basic properties of conormal
distributions to manifolds with corners is quite straight-forward if
the singular submanifold hits the boundary in a 'product-type' way. 

For lack of space we 
do not treat the transformation of the principal symbol under pull-back
and push-forward. However,
this is important for the composition formula for pseudodifferential
operators (see the references above). 
\subsubsection{Conormal distributions on manifolds}
\label{DGsubsubseccon}
\begin{definition} \label{DGdefcon}
Let $Z$ be a manifold and $Y\subset Z$ a submanifold. A distribution
$u\in \Dcal(Z)$ is {\em conormal} with respect to $Y$ if, for some $m\in\R$,
\begin{itemize}
\item
$u$ is smooth on $Z\setminus Y$, and
\item
in any local coordinate system $x:U\subset Z\to\R^n$ sending $Y\cap U$ to
$\R^k\times\{0\}^{n-k}\subset\R^n$ there is a representation
\begin{equation} \label{DGeqcon}
u(y,z) = \int_{\R^{n-k}} e^{iz\zeta} a(y,\zeta)\,d\zeta
\end{equation}
where $y=(x_1,\ldots,x_k), z=(x_{k+1},\ldots,x_n)$ and $a$ is a smooth
function on $(Y\cap U)\times \R^{n-k}$ with asymptotics 
\begin{equation}\label{DGeqconasymp}
a(y,\zeta) \sim \sum_{j=0}^\infty a_{m-j}(y,\zeta),
\end{equation}
as $|\zeta|\to\infty$,
where $a_l$ is homogeneous of degree $l$ in $\zeta$, for each $l$.%
\footnote{
The meaning of the asymptotics is that, for any $N$, if $a^{(N)}$ is the 
sum up to the term $a_{-N}$ then 
$|a(y,\zeta)-a^{(N)}(y,\zeta)| \leq C|\zeta|^{-N-1}$, plus analogous
estimates for all derivatives in $y$ and $\zeta$. 

The {\em order} of $u$ is defined to be $m+(n-2k)/4$, if $a_m\not\equiv 0$.
}
\end{itemize}
\end{definition}
Note that \eqref{DGeqcon} is simply the inverse Fourier transform in $z$
(i.e.\ 'transversal' to $Y$), with smooth dependence 
on the parameter $y\in Y$.
\begin{examples} \label{DGexcon}
\
\begin{enumerate}
\item
For $Z=\R$ and $Y=\{0\}$ the distributions $\delta$ and $p.v.\frac1x$
are conormal, and also all of their derivatives and anti-derivatives, 
which include for example $x_+^\alpha$ for $\alpha\in \N_0$
(the function vanishing for $x\leq0$ and equal to $x^\alpha$ for $x>0$).
Any $(-1)$-homogeneous distribution is a linear combination of $\delta$
and $p.v.\frac1x$,
so a conormal distribution (with $m\in\Z$) may be thought
of as 'series' of such terms, of increasing regularity.%
\footnote{
More generally, one can define conormality without reference
to the Fourier transform: \eqref{DGeqcon}
and \eqref{DGeqconasymp} are equivalent to the existence of
distributions $u_s(y,\cdot)\in\Dcal'(\R^{n-k}) \cap
C^\infty(\R^{n-k}\setminus 0)$, 
homogeneous
of degree $s$ and depending smoothly on the parameter $y$, such that
$$ u - \sum_{j=0}^N u_{j+n-k-m}\in C^\infty(\R^k,C^{N'}(\R^{n-k}))$$
(locally) for all $N$, where $N'=N-C$, with $C$ only depending on $n$. 
This may look weaker than the definition
above, but is actually equivalent (exercise!).
}
\item
For $Z=\R^n$, $n>1$, and $Y=\{0\}$ there is much more freedom since
now the space of $l$-homogeneous distributions is infinite-dimensional for
each $l$. The simplest example is $\delta$ again.
\item
For $Z=\R^n\times\R^n$ (with coordinates $w,w'\in\R^n$) and
$Y=\{w=w'\}$ (the diagonal) the conormal distributions are
the integral kernels of classical pseudodifferential operators on $\R^n$
since \eqref{DGeqcon} precisely amounts to their 'usual' definition, in the
coordinates $y=w$, $z=w-w'$, see for example \cite{DGShu:PDOST}, Section 3.7.
 (And similarly for \PDO s on any manifold.)
The order of the conormal distribution is the order of the operator
in the usual sense.
For example, the differential operator $P=\sum_\alpha a_\alpha(w)(\partial
/\partial w)^\alpha$ has kernel $P(w,w') = \sum_\alpha a_\alpha(w)
\delta^{(\alpha)}(w-w')$.
\end{enumerate}
\end{examples}
\begin{remarks}\label{DGremcondist}
\
\begin{itemize}
\item 
It's not obvious, but the definition is actually independent of the chosen
coordinate system $x$, see \cite{DGHor:ALPDOIII}. 
Of course, $a$ will depend on the choice of coordinates, but its leading term
$a_m$ is invariant if considered as section of the conormal bundle of $Y$. 
It is called the {\em principal symbol} of $u$. One easily sees that it depends
only on the restriction of $u$ to arbitrarily small neighborhoods of $Y$.
\item
The definition carries over immediately to distribution densities or,
more generally, to distributions with values in any bundle over $Z$.
\end{itemize}
\end{remarks}

We now consider push-forward and pull-back of distributions under a
smooth map $f:W\to Z$. The proofs of the following theorems are quite
easy, given the coordinate invariance of Definition \ref{DGdefcon}. 
They can be found in \Mnotes.
The push-forward of any distribution density $\mu$ on $W$ is a distribution
density on $Z$ (supposing, as usual, that $f$ is proper on $\supp \mu$), see
the Appendix.
The question arises whether conormality of $\mu$ with respect to a
submanifold $X\subset W$  implies
conormality of $f_*\mu$.
The answer is no in general%
\footnote{
Example: If $f$ is bijective and smooth then  
$\singsupp f_*\mu \supset f(\singsupp \mu)$. 
The latter need not be (contained in) a submanifold even if
$\singsupp\mu$ is, therefore
$f_*\mu$ may be not conormal even if $\mu$ is.}; it is a very tricky problem
to determine precise conditions when it is true.
It depends essentially on the behavior of the fibers of
$f$ and their tangency to  $X$. We only consider the simplest case of a
fibration whose fibers meet $X$ transversally in isolated points only,
which is enough for many purposes.

\begin{theorem}[Push-forward of conormal distributions]
\label{DGthpfdist}
Let $f:W\to Z$  be a fibration between manifolds. Let $X$ be a submanifold
of $W$ such that for each $x\in X$, the tangent spaces to $X$ and to the
fiber $f^{-1}(f(x))$ through $x$ intersect only in zero.%
\footnote{
Equivalently, $d(f_{|X})$ is injective.}
Let $\mu$ be a distribution density on $W$, conormal with respect to $X$,
and assume $f$ is proper on $\supp\mu$.
\begin{enumerate}
\item[(a)]
If $f_{|X}$ is a diffeomorphism onto $Z$ then $f_*\mu$ is smooth.
\item[(b)]
Otherwise, $f(X)$ is a proper submanifold of $Z$
and $f_*\mu$ is conormal with respect to $f(X)$.
\end{enumerate}
\end{theorem}
Thus, the 'vertical' (in fiber direction) singularities get integrated out,
while the others remain, as in the simple example of
 $f:\R^2\to\R$, $(x,y)\mapsto x$:\\
For $X=\{(x,0):x\in\R\}$ one has, for example,
$f_*(\delta(y)\,dxdy) = 1\cdot dx$ (case (a)), and for $X=\{(0,0)\}$
one has $f_*(\delta(x)\delta(y)\,dxdy) = \delta(x)\,dx$ (case (b)).
(Calculations done using \eqref{DGeqpf}.)
See Remark \ref{DGrempft}.3 concerning the support condition.

For pull-back the situation is different: While 
the pull-back for general distributions is only defined under fibrations
(and conormality is always preserved under a fibration),
a {\em weaker}  condition
on $f$ already allows to pull back conormal distributions:
\begin{theorem}[Pull-back of conormal distributions]
 \label{DGthpbdist}
Let $Y\subset Z$ be a submanifold, and assume that $f:W\to Z$ 
is transversal{$\,$}%
\footnote{
I.e.\ if $z=f(x)\in Y$ then $T_zZ$ is spanned by $T_zY$ and 
$df(T_xW)$.
}
 to $Y$.Then $f^{-1}(Y)$ is a submanifold of $W$, and if
$u$ is a distribution on $Z$, conormal with respect to $Y$, then $f^*u$
is a distribution on $W$ which is conormal with respect to $f^{-1}(Y)$.
\end{theorem}
We saw in Examples \ref{DGexpfpb} that we also need to {\em multiply}
distributions. A complete discussion of when this is possible would lead us 
too far astray, so we'll just sketch the procedure which allows to define a
product in this context.

For two functions $u_1,u_2$ on a manifold $Z$, we can translate 
the trivial identity
$u_1(x)u_2(x) = u_1(x)u_2(y)_{|x=y}$ into geometric terms as
$$ u_1u_2 = i^* (u_1\times u_2) $$
where $i:Z\to Z\times Z,x\mapsto (x,x)$ is the diagonal inclusion
and $(u_1\times u_2)(x,y)=u_1(x)u_2(y)$ defines the direct product
of $u_1$ and $u_2$ as function on $Z\times Z$.

When trying to generalize this to distributions $u_1,u_2$ on $Z$, 
we first note that  the direct product is well-defined as
a distribution on $Z\times Z$ (since $u_1$ and $u_2$
'depend on different sets of variables' in $Z\times Z$). 
The problem arises with the pull-back: $i$ is certainly not a fibration
(it's not even surjective), so one would hope to apply Theorem
\ref{DGthpbdist}. But this fails since usually $u_1\times u_2$ is not
conormal, even if $u_1$ and $u_2$ are conormal;
for example for $Z=\R$ and $u_1=u_2=p.v.\frac1x$ one gets
$(p.v. \frac1x)(p.v.\frac1y)$, which is not
conormal (since its singular support, the union
of both coordinate axes, is not a manifold).
\medskip

The following theorem allows a way out (at least in some situations):
\begin{theorem}[Direct product of conormal distributions]
\label{DGthproddist}
Let $u_i$ be distributions on manifolds $Z_i$, 
conormal with respect to submanifolds
$Y_i$, for $i=1,2$.
Then the direct product $u_1\times u_2$ can be written
$u_1\times u_2 = v+w$
where $v$ is  conormal with respect to
$Y_1\times Y_2$ and $w$ has wave front set contained in any given conic
neighborhood of
$(Z_1\times N^*Y_2) \cup (N^*Y_1\times Z_2)$.
\end{theorem}
We will not define wave front sets here, see \cite{DGHor:ALPDOI}. 
$N^*X_i$ denotes the conormal bundle.
The point is that in the applications we have in mind
(Examples \ref{DGexpfpb}) the product is integrated in the end,
and then the position of $WF(w)$ guarantees that the 
term resulting from $w$  is smooth (by a generalization of Theorems
\ref{DGthpfdist}(a) and \ref{DGthpbdist}), so the singularities are determined
only by the conormal part $v$. 

\subsubsection{Composition of pseudodifferential operators}
\label{DGsubsubsecdistpdo}
Let us check how these results show that the composition of
two pseudodifferential operators on a manifold $X$
is a pseudodifferential operator%
\footnote{\label{DGfnprop}%
The support condition in Theorem \ref{DGthpfdist}
translates into the condition that at least one of the factors
is properly supported (i.e.\ the  two projections $X^2\to X$ are proper on
the support of the integral kernel). 
We will neglect this in the following discussion. See also Remark
\ref{DGremclasspdo}.3.}. For simplicity we take $X=\R$ although the general
case works precisely the same way.

Thus, we are given distributions $A,B$ on $\R^2$, conormal with
respect to the diagonal (see Example \ref{DGexcon}.3). 
To avoid confusion later on, we will
write $A\in\Dcal'(X_1\times X_2)$, $B\in\Dcal'(X_2\times X_3)$ 
although $X_1=X_2=X_3=\R$. The diagonals will be denoted
$\Delta_A\subset X_1\times X_2$ and $\Delta_B\subset X_2\times X_3$.

The composition of $A$ and $B$ has integral kernel given by
\eqref{DGeqcomp}. Here, the product should be expanded, as explained
above, as diagonal pull-back of the direct product. However, matters
can be simplified slightly. The $\pi_1,\pi_3$ pull-backs can be omitted
 since clearly one also has
\begin{equation} \label{DGeqCeq}
 C= \pi_{2*}(i^* (A\times B))
\end{equation}
where $i$ is the embedding
$$ i: X_1\times X_2\times X_3 \to X_1\times X_2\times X_2\times X_3,
\quad (x_1,x_2,x_3)\mapsto (x_1,x_2,x_2,x_3).$$

Using Theorem \ref{DGthproddist} we write 
\begin{equation} \label{DGeqprodsum}
A\times B= v+w,
\end{equation}
with $v$ conormal with respect to $\Delta_A\times\Delta_B$
and $WF(w)$ close to $(X_1\times X_2\times N^*\Delta_B) \cup
(N^*\Delta_A\times X_2\times X_3)$.

We first analyze $\pi_{2*}(i^*v)$:
$i$ is transversal to $\Delta_A\times\Delta_B$ since
the image of $di$ is $\{(\alpha,\beta,\beta,\gamma)\}$ and the
tangent space of $\Delta_A\times\Delta_B$ is $\{(\delta,\delta,\eps,\eps)\}$
(all free variables in braces range over $\R$),
and these two subspaces clearly span $\R^4$.
Therefore, Theorem \ref{DGthpbdist} shows that $i^*v$ is conormal
with respect to $\Delta':=i^{-1}(\Delta_A\times\Delta_B) = \{x_1=x_2=x_3\}$,
the space diagonal.
Finally, the tangent spaces to $\Delta'$ and the fiber of $\pi_2$ 
are $\{(\alpha,\alpha,\alpha)\}$ and $\{(0,\beta,0)\}$, so they have
zero intersection, and $\pi_2(\Delta')=\Delta_C$, the diagonal in $X_1\times
X_3$. Therefore, Theorem \ref{DGthpfdist}(b) applies, so
\begin{equation} \label{DGeqvcon}
 \pi_{2*}(i^*v) \text{ is conormal with respect to the diagonal.}
\end{equation}

Finally, we analyze $w$, using standard results on wave front sets. First, 
by Theorem 8.2.4 in \cite{DGHor:ALPDOI}, the pull-back $i^*w$ is defined
as a distribution if $WF(w)\cap N^*({\rm Im}\, i)=\emptyset$, and then
$WF(i^*w)\subset i^*(WF(w))$.
Now $N^*({\rm Im}\,  i)=\{(0,\alpha,-\alpha,0)\}$ at every point, 
and this has non-zero angle
with the fiber of $N^*\Delta_A\times X_2\times X_3$, which is
$\{(\alpha,-\alpha,0,0)\}$, and similarly with the fiber of
$X_1\times X_2\times N^*\Delta_B$. 
Therefore, by choosing $WF(w)$ close enough to these latter sets, 
 we may assume that
$WF(w)\cap N^*({\rm Im}\, i)=\emptyset$. Also, $WF(i^*w)$ is contained
in a small conic neighborhood of $i^*$ of these sets, i.e.\ of
(fiberwise)
$\{(\alpha,-\alpha,0)\}\cup\{(0,\alpha,-\alpha)\}$.
By another standard theorem (see \cite{DGSim:PDO}, ex.\ 6.7.8)
  the push-forward
$\pi_{2*}u$ of a distribution $u$ on $X_1\times X_2\times X_3$ is
smooth unless $WF(u)$ hits the conormal space to the fiber of $\pi_2$.
Since the latter is $\{(\alpha,0,\beta)\}$, this is clearly
not the case for $u=i^*w$, so finally we obtain:
$$ \pi_{2*}(i^*w)\text { is smooth.} $$
This together with \eqref{DGeqvcon}, \eqref{DGeqCeq} and \eqref{DGeqprodsum}
shows that $C$ is conormal with respect to the diagonal, i.e.\ the integral
kernel of a pseudodifferential operator.

\subsubsection{Conormal distributions on manifolds with corners}
\label{DGsubsubsecconmwc}
The definition of conormal distributions depends in an essential way on the
fact that 'normal slices' to $Y$ in $Z$ look the same at every point of $Y$.
Therefore, it would be problematic to try to define conormality for
$Z=\Rplus^2$, $Y=\{(x,x):x\in\Rplus\}$: At zero, there is not even a
reasonable candidate for a 'normal slice'!

However, if $Z$ looks like $\R^{n-k}\times Y$ near $Y$ then Definition 
\ref{DGdefcon} makes sense literally even if $Y$ is a mwc, when we
require the $y$-dependence to be smooth up to the boundary everywhere.
For example, this is the case for $Z=\Rplus\times\R$, $Y=\Rplus\times\{0\}$ or
$Z=[\Rplus^2,0]$, $Y=\Delta_b:=\{\theta=\pi/4\}$ (using polar coordinates on
$Z$, see Figure \ref{DGfigpolar}(b)).

Thus, we have defined {\em distributions on a mwc $Z$ which are conormal
with respect to an interior p-submanifold, smoothly at the boundary}
(cf.\ \ref{DGsubsubsecgenbu}.4; 'interior' means that 
$Y\not\subset\partial Z$).
One may actually allow nice (rather than smooth) behavior at the
boundary, with respect to a given index family. This gives 
{\em nice conormal distributions.} The precise definitions are
quite straight-forward, we leave them as exercise to the reader.
(See \Mnotes.)

The {\em push-forward} for nice conormal distributions may be analyzed
by a combination of Theorems \ref{DGthpft} and \ref{DGthpfdist}. Since 
the point of Theorem \ref{DGthpft} was to allow maps $f$ more general than
fibrations, the assumption on $f$ should be: $f$ is a \Mb-fibration, and 
a fibration (of mwc's) in some neighborhood of $X$, satisfying the additional
transversality condition in Theorem \ref{DGthpfdist}. Then the result
of the push-forward is nice conormal again. The proof is straightforward by 
use of a partition of unity, and is left as an exercise.

Similarly {\em pull-back} of nice conormal distributions is easy by
combining Theorems \ref{DGthpbt} and \ref{DGthpbdist}. $f$ needs to
be a \Mb-map transversal to $Y$, then the pull-back by $f$ of a nice
conormal function is nice conormal with respect to $f^{-1}(Y)$.

%% file: figpf.pstex_t
\begin{picture}(0,0)%
\epsfig{file=figpf.pstex}%
\end{picture}%
\setlength{\unitlength}{3947sp}%
\begingroup\makeatletter\ifx\SetFigFont\undefined%
\gdef\SetFigFont#1#2#3#4#5{%
  \reset@font\fontsize{#1}{#2pt}%
  \fontfamily{#3}\fontseries{#4}\fontshape{#5}%
  \selectfont}%
\fi\endgroup%
\begin{picture}(5525,2177)(151,-5407)
\put(205,-4265){\makebox(0,0)[lb]{\smash{\SetFigFont{10}{12.0}{\rmdefault}{\mddefault}{\updefault}$B$}}}
\put(151,-3449){\makebox(0,0)[lb]{\smash{\SetFigFont{10}{12.0}{\rmdefault}{\mddefault}{\updefault}$y$}}}
\put(1728,-5026){\makebox(0,0)[lb]{\smash{\SetFigFont{10}{12.0}{\rmdefault}{\mddefault}{\updefault}$x$}}}
\put(2217,-3395){\makebox(0,0)[lb]{\smash{\SetFigFont{10}{12.0}{\rmdefault}{\mddefault}{\updefault}$\eta_2$}}}
\put(3631,-5026){\makebox(0,0)[lb]{\smash{\SetFigFont{10}{12.0}{\rmdefault}{\mddefault}{\updefault}$\xi_2$}}}
\put(2217,-4917){\makebox(0,0)[lb]{\smash{\SetFigFont{10}{12.0}{\rmdefault}{\mddefault}{\updefault}$B$}}}
\put(2489,-5189){\makebox(0,0)[lb]{\smash{\SetFigFont{10}{12.0}{\rmdefault}{\mddefault}{\updefault}$(\xi_2,\eta_2)\mapsto\xi_2\eta_2$}}}
\put(4555,-5189){\makebox(0,0)[lb]{\smash{\SetFigFont{10}{12.0}{\rmdefault}{\mddefault}{\updefault}$(\xi_1,\eta_1)\mapsto\xi_1$}}}
\put(4120,-3449){\makebox(0,0)[lb]{\smash{\SetFigFont{10}{12.0}{\rmdefault}{\mddefault}{\updefault}$\eta_1$}}}
\put(5533,-5026){\makebox(0,0)[lb]{\smash{\SetFigFont{10}{12.0}{\rmdefault}{\mddefault}{\updefault}$\xi_1$}}}
\put(749,-5407){\makebox(0,0)[lb]{\smash{\SetFigFont{10}{12.0}{\rmdefault}{\mddefault}{\updefault}(a)}}}
\put(2869,-5407){\makebox(0,0)[lb]{\smash{\SetFigFont{10}{12.0}{\rmdefault}{\mddefault}{\updefault}(b)}}}
\put(4935,-5407){\makebox(0,0)[lb]{\smash{\SetFigFont{10}{12.0}{\rmdefault}{\mddefault}{\updefault}(c)}}}
\put(640,-5189){\makebox(0,0)[lb]{\smash{\SetFigFont{10}{12.0}{\rmdefault}{\mddefault}{\updefault}$\pi_1\circ\beta$}}}
\put(4174,-4917){\makebox(0,0)[lb]{\smash{\SetFigFont{10}{12.0}{\rmdefault}{\mddefault}{\updefault}$A$}}}
\put(967,-4972){\makebox(0,0)[lb]{\smash{\SetFigFont{10}{12.0}{\rmdefault}{\mddefault}{\updefault}$A$}}}
\put(749,-4645){\makebox(0,0)[lb]{\smash{\SetFigFont{10}{12.0}{\rmdefault}{\mddefault}{\updefault}$y/x$}}}
\put(1510,-3558){\makebox(0,0)[lb]{\smash{\SetFigFont{8}{9.6}{\rmdefault}{\mddefault}{\updefault}$y/x=$const}}}
\end{picture}

%% file: figtriple.pstex_t
\begin{picture}(0,0)%
\epsfig{file=figtriple.pstex}%
\end{picture}%
\setlength{\unitlength}{3947sp}%
\begingroup\makeatletter\ifx\SetFigFont\undefined%
\gdef\SetFigFont#1#2#3#4#5{%
  \reset@font\fontsize{#1}{#2pt}%
  \fontfamily{#3}\fontseries{#4}\fontshape{#5}%
  \selectfont}%
\fi\endgroup%
\begin{picture}(3877,5801)(1651,-4947)
\put(2977,-1353){\makebox(0,0)[lb]{\smash{\SetFigFont{10}{12.0}{\rmdefault}{\mddefault}{\updefault}$\fff$}}}
\put(4394,-2382){\makebox(0,0)[lb]{\smash{\SetFigFont{10}{12.0}{\rmdefault}{\mddefault}{\updefault}$\Delta_{3b}$}}}
\put(4792,-1440){\makebox(0,0)[lb]{\smash{\SetFigFont{10}{12.0}{\rmdefault}{\mddefault}{\updefault}$\ff_1$}}}
\put(4834,-836){\makebox(0,0)[lb]{\smash{\SetFigFont{10}{12.0}{\rmdefault}{\mddefault}{\updefault}$\bface_2$}}}
\put(2201,-922){\makebox(0,0)[lb]{\smash{\SetFigFont{10}{12.0}{\rmdefault}{\mddefault}{\updefault}$\bface_1$}}}
\put(2330,-2261){\makebox(0,0)[lb]{\smash{\SetFigFont{10}{12.0}{\rmdefault}{\mddefault}{\updefault}$\ff_2$}}}
\put(5396,-1482){\makebox(0,0)[lb]{\smash{\SetFigFont{10}{12.0}{\rmdefault}{\mddefault}{\updefault}$x_1$}}}
\put(3193,-3125){\makebox(0,0)[lb]{\smash{\SetFigFont{10}{12.0}{\rmdefault}{\mddefault}{\updefault}$\pitilde_{3}$}}}
\put(2996,-2304){\makebox(0,0)[lb]{\smash{\SetFigFont{10}{12.0}{\rmdefault}{\mddefault}{\updefault}$\bface_3$}}}
\put(4511,-3430){\makebox(0,0)[lb]{\smash{\SetFigFont{10}{12.0}{\rmdefault}{\mddefault}{\updefault}$\rb$}}}
\put(5528,-3520){\makebox(0,0)[lb]{\smash{\SetFigFont{10}{12.0}{\rmdefault}{\mddefault}{\updefault}$x_1$}}}
\put(3140,-3961){\makebox(0,0)[lb]{\smash{\SetFigFont{10}{12.0}{\rmdefault}{\mddefault}{\updefault}$\ff$}}}
\put(4245,-4582){\makebox(0,0)[lb]{\smash{\SetFigFont{10}{12.0}{\rmdefault}{\mddefault}{\updefault}$\Delta_b$}}}
\put(3001,689){\makebox(0,0)[lb]{\smash{\SetFigFont{10}{12.0}{\rmdefault}{\mddefault}{\updefault}$x_3$}}}
\put(1801,-4936){\makebox(0,0)[lb]{\smash{\SetFigFont{10}{12.0}{\rmdefault}{\mddefault}{\updefault}$x_2$}}}
\put(1651,-2911){\makebox(0,0)[lb]{\smash{\SetFigFont{10}{12.0}{\rmdefault}{\mddefault}{\updefault}$x_2$}}}
\put(2176,-4336){\makebox(0,0)[lb]{\smash{\SetFigFont{10}{12.0}{\rmdefault}{\mddefault}{\updefault}$\lb$}}}
\put(2851, 14){\makebox(0,0)[lb]{\smash{\SetFigFont{10}{12.0}{\rmdefault}{\mddefault}{\updefault}$\ff_3$}}}
\put(4801,-211){\makebox(0,0)[lb]{\smash{\SetFigFont{10}{12.0}{\rmdefault}{\mddefault}{\updefault}$\Delta_b^{(3)}$}}}
\end{picture}

%% file: DGpde.tex
\section{Partial Differential Equations} \label{DGsecpde}
We now turn to the \Mb-calculus in the narrow sense: the construction
of  parametrices for the 'simplest' class of non-uniformly elliptic
differential operators, the \Mb-differential
(or 'totally characteristic') operators on a
manifold with corners. Since this is an introductory article,
we only consider manifolds with boundary (and mostly even just
$\Rplus$). In this case, the operators are also called Fuchs
type or cone operators.
They have been studied by many authors, see e.g.\ 
\cite{DGBruSee:RESORSO}, \cite{DGLes:OFTCSAM},
\cite{DGSchul:PDOMWS}; some of them used pseudodifferential operator 
(\PDO) techniques. For a more complete list of references see
\cite{DGLauSei:PDAMWBCBCA}.

In \cite{DGMelMen:EOTCT} and \Mbook\ the central results (on action and
composition) were proved by direct calculation. There is a more systematic
way to do this, using the Pull-Back and Push-Forward Theorem (as 
indicated in Figure \ref{DGfigbcalc} and
Examples \ref{DGexpfpb}), and this shows more clearly
the geometric reasons for the precise form of
these results. Melrose alludes to this
often in \Mbook, and (probably) proceeds like this
 in the part of \Mnotes\ that is not publicly
available yet. Therefore, we will follow this systematic method here,
but will be sketchy otherwise.

\begin{definition}\label{DGdefbdiff}
A differential operator on a manifold with boundary $X$ is a
{\em \Mb-differential operator} of order $m$, in symbols $P\in\Diff_b^m(X)$, 
if it has smooth coefficients and,
in any coordinate system around a boundary point in which
$x$ denotes a boundary defining function and $y=(y_1,\ldots,y_{n-1})$ 
coordinates in the boundary, it has the form
\begin{equation}\label{DGeqbop}
P= \sum_{j+|\alpha|\leq m} a_{j\alpha}(x,y)\, (x\partial_x)^j 
\partial_y^\alpha 
\end{equation}
($\alpha$ runs over multi-indices in $\N^{n-1}$, 
and $\partial_x=\partial/\partial x$ etc.) with 
functions $a_{j\alpha}$ that are smooth up to the boundary.%
\footnote{
For Fuchs type operators or cone operators one usually
multiplies this with $x^{-\mu}$ for some positive number $\mu$. This
is inessential for parametrix constructions since the factor
can be transferred to the parametrix. However, it makes an
essential difference for the spectral theory (the analysis
of the resolvent), which we don't consider here. See for example 
\cite{DGBruSee:RESORSO}, \cite{DGGil:FAEHTNSAECO}, \cite{DGLoy:SREPO}.
}
\end{definition}

On manifolds without boundary the power of the \PDO\ calculus
derives from two facts:
\begin{itemize}
\item The (principal) symbol of a \PDO\ describes the operator 
fairly precisely (up to lower order), and 
\item
symbols are easy to invert.
\end{itemize}
We discuss this shortly in Subsection \ref{DGsubsecpdo}.
Recall from Example \ref{DGexcon}.3 
that a \PDO\ on a manifold $X$ is just a distribution on 
$X\times X$ which is conormal with respect to 
\begin{equation} \label{DGeqdiag}
\Delta = \{(p,p):\,p\in X\} \subset X\times X \quad (\text{the diagonal}),
\end{equation}
and that the (principal) symbol is determined purely by the singular behavior
at $\Delta$ (see Remark \ref{DGremcondist}).
 
Therefore, for an extension of this theory to \Mb-(pseudo)-differential
operators (on $\Rplus$, say)
one expects that it should be useful to define a class
of distributions on $\Rplus\times\Rplus$ by giving precise descriptions
of their behavior in all possible limits. These are%
\footnote{
As already in Section \ref{DGsecgeom} we always restrict attention 
to compact parts of the spaces involved. Therefore, we do not consider
limits 'at $\infty$'.
}:
\begin{enumerate}
\item[(a)]
Approaching (and at) the interior of the diagonal $\Delta$.
\item[(b)]
Approaching $\partial(\Rplus^2) \setminus (0,0)$.
\item[(c)]
Approaching the corner $(0,0)$.
\end{enumerate}
While there are obvious candidates for (a) and (b) (conormal singularity
and complete asymptotics ('niceness' in Definition \ref{DGdefasymp}),
respectively),
it is unclear what a good description for (c) might be: Simultaneously
one needs to describe the behavior of the 
conormal singularity on $\Delta$
when approaching
zero, and of the coefficients in the asymptotics of (b). 

This is most elegantly solved by blowing up the corner $(0,0)$,
which has the effect of separating the sets $\Delta$ and
$\partial\Rplus^2\setminus(0,0)$. Let us illustrate this by a simple
example:
\begin{example} \label{DGexbcalc}
On $\Rplus$, consider the simplest non-trivial \Mb-equation
$$ \left(x\frac{d}{dx} + c\right)u(x) = v(x),$$
for some fixed $c\in\C$. Multiplying by $x^{c-1}$ we get
$\frac{d}{dx}(x^c u(x)) = x^{c-1}v(x)$, which can be integrated to yield
$$ u(x) = x^{-c} a_0 + x^{-c}\int_0^x (x')^{c-1} v(x')\, dx'$$
with $a_0 = (x^cu(x))_{|x=0}$ (assuming $(x')^{c-1}v(x')$ is integrable near
zero). For simplicity, consider only the solution with
$a_0=0$. We can write it as
$$ u(x) = \int_0^\infty K(x,x') v(x')\, \frac{dx'}{x'}$$
with
\begin{equation}\label{DGeqmodel}
 K(x,x') = \left(\frac{x'}x\right)^c H(x-x').
\end{equation}
Here, $H(t)=1$ for $t>0$ and $H(t)=0$ for $t\leq 0$.
Thus, $K$ is (the kernel of) an inverse of the operator
$x\frac{d}{dx}+c$ (on suitable function spaces), so it should be an
example of a \Mb-pseudodifferential operator.

As expected, $K$  is singular at the two coordinate axes (actually only at
$\{x'=0\}$)and at the diagonal $\Delta=\{x=x'\}\subset\Rplus^2$
(unless $c$ happens to be a positive integer).
Looking at $K$ on the blow-up space $[\Rplus^2,(0,0)]$ 
(i.e.\ at $\beta^*K$ for $\beta:[\Rplus^2,0]\to\Rplus^2$ the blow-down map)
means rewriting $K$ in terms of coordinates on this space; using projective
coordinates $x$ and $$s=x'/x,$$ for example, we get:
$$ \beta^*K(x,s) = s^c H(1-s).$$
This is much nicer than \eqref{DGeqmodel} since:
\begin{itemize}
\item
The submanifolds at which $\beta^*K$ is singular are disjoint.%
\footnote{
This is better than needed. 'Normal crossings' (i.e.\ locally
looking like coordinate subspaces in a suitable coordinate system) would be
enough. This is satisfied by the boundary faces of a mwc (and by the
four distinguished submanifolds $\lb,\rb,\ff,\Delta_b$ of
$[\Rplus^2,(0,0)]$, see below), but not by the singular support of $K$.
}
(They are $\lb=\{s=\infty\}$, $\rb=\{s=0\}$ and $\Delta_b=\{s=1\}$, see
Figure \ref{DGfigpolar}(b).)%
\footnote{
Strictly speaking, one should also check $\beta^*K$ in the 
$x',x/x'$ coordinates (in the sequel we will neglect this when it gives 
no information). Alternatively, you may use instead
the coordinates $\rho=x+x'$, $\tau=(x-x')/(x+x')$
(see \ref{DGsubsecbu}), then
$$ \beta^*K (\rho,\tau) = \left(\frac{1-\tau}{1+\tau}\right)^c H(\tau)$$
and $\lb=\{\tau=-1\}$, $\rb=\{\tau=1\}$, $\Delta_b=\{\tau=0\}.$
}
\item
$\beta^*K$ is nice outside $\Delta_b$.
\item
$\beta^*K$ has a conormal singularity at $\Delta_b$, smoothly up to
the boundary (i.e.\ up to $\{x=0\}$).
\end{itemize}
\end{example}

Thus, we are lead to define (kernels of) \Mb-\PDO s on a manifold $X$
with boundary as a certain class of distributions on a blown-up 
space $X^2_b$ (where $X^2_b=[\Rplus^2,(0,0)]$ for $X=\Rplus$). These
are then considered as kernels on $X^2$ by use of
the identification of the interiors of $X^2_b$ and $X^2$ via the
blow-down map.

What should we expect the {\em symbol} of a \Mb-\PDO\  to be?
$X^2_b$ has four distinguished submanifolds: $\lb,\rb,\ff$ and
\begin{equation} \label{DGeqbdiag}
\Delta_b = \overline{\beta^{-1}({\rm interior}(\Delta))} = \{s=1\},
\end{equation}
and \Mb-\PDO s are characterized by their behavior at them.
Just as the symbol in the boundaryless case describes the leading
behavior at $\Delta$ (which is the only distinguished submanifold
then), one might expect the symbol in the \Mb-calculus to describe
the leading behavior at $\Delta_b,\ff,\lb,\rb$.
In fact, it turns out that  only the
former two are needed since their vanishing implies compactness
(between spaces that are determined by the latter two!).

The \Mb-calculus is introduced in two steps. They are motivated by 
the construction of a parametrix of an elliptic \Mb-differential
operator $P$.
First,
the 'small calculus' $\Psi_b^*$ is constructed; it allows inversion
of the symbol on $\Delta_b$ and thus, by the usual iterative procedure,
inversion of $P$ modulo errors that are smooth on $\Delta_b$.
Since this game is played away from $\lb,\rb$, operators in $\Psi_b^*$
are assumed to vanish to infinite order there.
However, as we saw in the example above, the inverse of even the simplest 
\Mb-differential operator is not of this type. This is reflected in the
fact that the 'remainders' in the parametrix construction, i.e. elements
of $\Psi_b^{-\infty}$, are not compact operators, even if
$X$ is compact. Therefore, in a second
step the 'full calculus' is introduced.
This allows inversion of the symbol at \ff\ (the 'conormal symbol'
or 'indicial operator'); the price to pay is non-trivial asymptotic
behavior at $\lb,\rb$.

For the sake of presentation we mostly work with the simplest manifold with
boundary, $\Rplus$; most ideas may be understood already in this case.%
\footnote{
Of course, one must resist the temptation to use simpler arguments
only suited for ordinary differential equations.
}
At the end of the chapter we add some remarks on the changes necessary
when dealing with a general manifold with boundary.

\subsection{Classical pseudodifferential operators}
\label{DGsubsecpdo}
We shortly summarize the essential ingredients of the classical \PDO\
calculus, and how they are used to find parametrices, i.e.\
approximate inverses, for elliptic
(pseudo-)differential operators on a compact manifold $X$.
Extensive treatments can be found in \cite{DGHor:ALPDOIII} and
\cite{DGShu:PDOST}, for example.
A similar (and more general) axiomatic treatment was given in
\cite{DGSchul:TIOSSS}.
\medskip

Given, for each $m\in\R$:
\begin{boldlist}
\item[(Op)]
Classes $\Psi^m$
of distributions on $X^2$, with $\Psi^m\subset\Psi^{m+1}$.
($\Psi^m$ is taken to be the set of
distributions on $X^2$ conormal with respect to
the diagonal $\Delta$, of order $m$.)
\item[(Symb)]
Classes $S^m$ of symbols, with $S^m\subset S^{m+1}$, and
symbol maps 
$\sigma:\Psi^m\to S^{[m]} := S^m/S^{m-1}$.
($S^m$ is taken to be the set of smooth
 functions on $T^*X$ with complete asymptotic expansions 
 in homogeneous components
 of order $\leq m$. Homogeneity refers to the covariable $\xi$, and the
 asymptotics is for $\xi\to\infty$.

$\sigma(u)$ is 
defined by the principal symbol of $u\in \Psi^m$, which
is a function on $N^*\Delta$ (see Remark \ref{DGremcondist}), 
using the canonical identification
$N^*\Delta \cong T^*X$. Note that $S^{[m]}\cong$ smooth functions on 
$T^*X\setminus 0$, homogeneous in $\xi$ of order $m$.)
\end{boldlist}
\medskip

The essential properties of these objects are:
\begin{boldlist}
\item[(Alg)]
$\bigcup_m \Psi^m$ and $\bigcup_m S^m$ are graded algebras,
and $\sigma$ is an algebra homomorphism.
(The products are taken as composition, defined by \eqref{DGeqcomp},
 and pointwise multiplication, 
respectively, and 'graded' means
 $P\in\Psi^m,Q\in\Psi^l\implies PQ\in \Psi^{m+l}$, and
similarly for symbols;
the main point is that $\sigma$ respects products.)
\item[(Exact)]
The sequence
\begin{equation} \label{DGeqses}
0\to \Psi^{m-1} \hookrightarrow \Psi^m \stackrel{\sigma}{\to} S^{[m]} \to 0
\end{equation}
is exact for every $m$.
This means:
\begin{enumerate}
\item[{\bf(E1)}]
For each $a\in S^m$ there is $P\in\Psi^m$ with $\sigma(P) = a \mod S^{m-1}$.
\item[{\bf(E2)}]
If $P\in\Psi^m$ and $\sigma(P)=0$ then $P\in\Psi^{m-1}$.
\end{enumerate}
\end{boldlist}
Part of (Alg) was checked in \ref{DGsubsubsecdistpdo}, and (Exact)
is straight-forward from the definitions.

Finally, we define:
An element $P\in\Psi^m$ is {\em elliptic} if $\sigma(P)$ is invertible
(then its inverse lies in $S^{-m}$ necessarily).
A {\em parametrix} of order $k$ for $P\in\Psi^m$ is a $Q\in\Psi^{-m}$
such that both $PQ-\Id$ and $QP-\Id$ lie in $\Psi^{-k}$.

The main fact is:
\begin{theorem}[Parametrix construction for elliptic \PDO]\label{DGthparam}
(Alg) and \mbox{(Exact)} above imply:
If $P\in\Psi^m$ is elliptic then it has a parametrix of any order.
\end{theorem}
\begin{quote}
{\small
Let us quickly recall the proof:
By ellipticity of $P$,  $\sigma(P)^{-1}$ is invertible with inverse
in $S^{-m}$. By
(E1), there is $Q\in\Psi^{-m}$ with $\sigma(Q) = \sigma(P)^{-1}$.
Then, by (Alg), $\sigma(PQ-\Id)=\sigma(P)\sigma(Q)-\sigma(\Id)=0$,
so by (E2) (with $m=0$) we have $R:=\Id-PQ\in\Psi^{-1}$.
Thus, $Q$ is a 'right' parametrix of order 1.
 Set 
$Q_k=Q(\Id+R+\ldots+R^{k-1})$, then $PQ_k = (\Id-R)(\Id+R+\ldots+R^{k-1})
=\Id-R^k$ and $R^k\in\Psi^{-k}$, so $Q_k$ is a right parametrix of order $k$.
By the same procedure we get a left parametrix $Q_k'$  order $k$.
Evaluating $Q_k'PQ_k$ in two ways one obtains that $Q_k-Q_k'\in\Psi^{-m-k}$,
and from this  that $Q_k$ is also a left parametrix
of order $k$.
\qed
}
\end{quote}

This may be refined slightly:
One also has
\begin{boldlist}
\item[(AC)] Asymptotic completeness: If $P_i\in\Psi^{m-i}$ for $i\in\N_0$
then there is $P\in\Psi^m$ with $P-\sum_{i=0}^N P_i\in\Psi^{m-N-1}$
for all $N$.
\end{boldlist}
This clearly implies that elliptic elements have parametrices of
order $\infty$ (usually just called parametrices). This improvement
is mainly cosmetic and not needed in most applications.

Note that these arguments were purely formal and did not use any properties
beyond (Alg), (Exact) and (AC). Therefore, the same result holds with different
choices of $\Psi^*,S^*,\sigma$.

However, in order to apply Theorem \ref{DGthparam} to problems of
differential equations, one needs:
\begin{boldlist}
\item[(Diff)]
$\Diff^m\subset\Psi^m$, where $\Diff^m$ denotes the differential
operators of order $m$ (with smooth coefficients).
\item[(Ell)]
The 'usual' elliptic operators one is usually interested in
(Dirac, Laplace) are elliptic in the sense above.%
\footnote{ 
For example, this is not satisfied if one takes $\Psi^m$ 
as above but lets $S^m=\Psi^m$, $\sigma=\id$,
which satisfies all other requirements!
}
\end{boldlist}

Finally, in order to make all of this
useful for analysis
(e.g.\ for proving regularity of solutions of elliptic PDE) one
needs:
\begin{boldlist}
\item[(Map)] \hspace*{-5pt} Mapping properties of $P\in\Psi^m$ 
(e.g. continuity on Sobolev spaces).
\item[(Neg)]
The remainders, i.e.\ the elements in $\Psi^{-\infty}$,
 are actually 'negligible'
(e.g.\ compact, trace class, smoothing, etc.).
\end{boldlist}
Of course,  the usual \PDO\ calculus has all of these properties.

\begin{remarks} \label{DGremclasspdo}
\
\begin{enumerate}
\item 
The motivation for the definition of $\Psi^*$
 lies in solving differential equations by Fourier transform,
which gives precise solutions for constant coefficient equations.
In this case, it suffices to invert the symbol. The rest is just
the algebra that's needed to make this method work for
non-constant coefficient equations.%
\footnote{
In other words, inversion of the symbol $\sigma(P)(x_0,\xi)$ corresponds
to inversion of (the principal part of) the constant coefficient operator
$P(x_0,D)$ obtained by freezing coefficients at $x_0$ (which acts on
$T_{x_0}X$), and the parametrix construction shows how to patch these local
inverses together.
}
Note that the Fourier transform does not appear explicitly. It is
stowed away in the definition of conormal distributions, and is
only used in the proofs of pull-back and push-forward theorem for these (cf.\
the hierarchy in Figure \ref{DGfigbcalc}).
\item 
When trying to construct $\Psi^*$, one has to find a compromise between
opposing forces: It has to be large enough to contain elliptic differential
operators and their parametrices, but small enough for $\Psi^{-\infty}$
to be actually negligible.
\item
If $X$ is not compact then composition of $P,Q\in\Psi^*$ may  be 
undefined (since the integral $\int P(x,y)\,Q(y,z)\,dy$ may diverge
'at $\infty$'). The simplest remedy is to replace $\Psi^*$ by 
$$\hspace*{4em} \Psi^*_{\rm{prop}} = \{P\in\Psi^*:\, K\Subset X 
\implies (\supp P)\cap (X\times K) \ \text{ is compact}\}, $$
the set of properly supported \PDO s.
Then everything goes through as before, except that remainders are still
smoothing but not compact, so that (Neg) should be localized:
\vspace{1ex}
\begin{boldlist}
\item[(Neg')]
If $P\in\Psi^{-\infty}_{\rm{prop}}$ and $\phi\in C_0^\infty(X)$ then
$\phi P$ (i.e.\ $P$ followed by multiplication by $\phi$) is negligible
(i.e.\ compact, trace class, etc.).
\end{boldlist}
\end{enumerate}
\end{remarks}

%%%%%%%%%%%%%%%%%%%%%%%%%%%%%%%%%%%%%%%%%%%%%%%%%%%%%
%%%%%%%%%%%%%%%%%%%%%%%%%%%%% Small calc. subsec
%%%%%%%%%%%%%%%%%%%%%%%%%%%%%%%%%%%%%%%%%%%%%%%%%%%%%

\subsection{The small \Mb-calculus} \label{DGsubsecsmall}
The small \Mb-calculus takes care of the conormal singularity on
the diagonal and its behavior near the boundary of the diagonal. 
However, it does
not admit non-trivial asymptotics near the boundary of $X^2$
away from the corner $(0,0)$.

We will first motivate the definition of the small \Mb-calculus by calculating
the kernels of \Mb-differential operators as distributions on $X^2_b$,
then discuss the annoying but non-negligible question of (half-)densities
(but see Footnote \ref{DGfndens}),
and finally define the small calculus and check its properties, in particular
composition.

Recall that we work on $X=\Rplus$, and that $X^2_b=[\Rplus^2,(0,0)]$.
\subsubsection{Kernels of \Mb-differential operators}
\label{DGsubsubsecmdiff}
For a start, let's be naive and consider as integral
kernel of an operator $P$ a distribution $K_P(x,x')$ on $X^2$ such that 
\begin{equation}\label{DGeqnaive}
(Pu)(x)=\int K_P(x,x')u(x')\,dx' \quad\text{('naive kernel').}
\end{equation}
The simplest \Mb-differential operator is the identity $\Id$.
Its kernel (on $X^2$) is $\delta(x-x')$. To obtain the kernel as distribution
on $X^2_b$, we simply rewrite this in terms of coordinates on $X^2_b$,
for example of the projective coordinates $x,s=x'/x$. Since $\delta$ is
homogeneous of degree $-1$, we get
$$ K_\Id = \delta(x-x') = \delta(x(1-\frac{x'}{x})) = \frac1x\delta(1-s).$$
Next, let us consider $x\partial_x$. Since $\delta'$ is homogeneous of
degree $-2$, we get
$$ K_{x\partial_x} = x\delta'(x-x') =  \frac1x\delta'(1-s).$$
Iterating this, one sees easily:
\begin{theorem}[\Mb-differential operators as naive kernels]
\label{DGthnaive}
The 'naive' kernels $($as in \eqref{DGeqnaive}$)$
 of \Mb-differential operators on $\Rplus$ are precisely the distributions
of the form
\begin{equation}\label{DGeqnaive2}
 \frac1x \sum_{\rm finite} a_j (x) \delta^{(j)}(1-s)
\end{equation}
with $a_j$ smooth up to $x=0$.
\end{theorem}

\subsubsection{Densities, half-densities and their \Mb-rethren}
Since we want to switch between different coordinate systems (as
above), we have to know how objects transform under such changes;
to  put it differently, we should define things invariantly. The idea of
'integral kernel' involves integration and therefore measures. There are
three obvious possibilities to take care of these:
\begin{itemize}
\item
Integral kernel is a function, acts on densities.
\item
Integral kernel is a density, acts on functions.
\item
Integral kernel is a 'density in $x'$', acts on functions.
\end{itemize}
In the first case, the result $\int K(x,x')\mu(x')$ would be a function.
Therefore, this would describe an operator mapping densities to functions.
But then two such operators cannot be composed%
\footnote{
unless one identifies functions with densities by (non-canonical)
choice of a fixed density, which is what we wanted to avoid in the first
place!}. The same problem occurs in the second case: Here functions are mapped
to densities. The third possibility avoids this problem, here functions are
mapped to functions, but now the symmetry between $x$ and $x'$ is broken.
There is a way to overcome this flaw as well:
\begin{itemize}
\item
Integral kernel is a half-density, acts on half-densities.
\end{itemize}
For a very short introduction to half-densities, see the Appendix.
Now there is only one kind of objects: Operator kernels are half-densities,
and they map half-densities to half-densities, so they can be composed
without making additional choices. Also, an integral kernel is a 'symmetric'
object, i.e.\ $x,x'$ can be interchanged freely.%
\footnote{
The pedantic reader will notice that now the formula 
\eqref{DGeqaction} translating action into the pull-back/push-forward world
does not make sense directly: Although $\pi_2^*v$ is well-defined for a 
half-density $v$ since $\pi_2$ is a fibration, the product $A\cdot \pi_2^*v$
is not quite a density on $X^2$ (a half-density in the first factor is missing),
so its push-forward is not defined. This can be remedied easily:
Fix any non-vanishing half-density $\mu$ on $X$. Then the push-forward
$\pi_{1*}(\pi_1^*\mu\cdot A \cdot \pi_2^*v)$ is a well-defined density on
the first factor $X$, and dividing it through $\mu$ one obtains the result,
which is immediately checked to be independent of the choice of $\mu$
(see \Mnotes). A similar remark applies to composition \eqref{DGeqcomp}.
}

Since (half-)densities are differential objects, we should, by Principle  
\ref{DGprincb} in the Introduction,
use \Mb-half-densities instead.
Let us illustrate this by an example:
Consider the identity again. Now using the same letter for
kernels and operators, we have
\begin{alignat*}{2}
 \Id &= \delta(x-x') |dx\,dx'|^{1/2} &&\text{as half-density}\\
   & = x\delta(x-x') \left|\dxx\,\dxxprime\right|^{1/2} &&
     \text{as \Mb-half-densities}\\
  & = \delta(1-s)  \left|\dxx\dss\right|^{1/2} &&
     \text{as \Mb-half-densities on $X^2_b$}
\end{alignat*}
where in the second line we used $(x')^{1/2}\delta(x-x') = x^{1/2}\delta(x-x')$
and in the last line the calculation before Theorem \ref{DGthnaive},
together with \eqref{DGeqbdenstrf}.
In general, one sees easily that the
$x^{-1}$ factor in Theorem \ref{DGthnaive} is always canceled when we use
half-densities, i.e.
\begin{equation}\label{DGeqdiffkernel}
\sum_j a_j(x) (x\partial_x)^j = \sum_j a_j(x) \delta^{(j)}(1-s)
          \dxshalf.
\end{equation}
Recall that $\{s=1\}$ is just the \Mb-diagonal $\Delta_b$.
An expression like $\sum a_j(x)\delta^{(j)}(1-s)$ is called
a Dirac distribution on $\Delta_b$.
We may now restate Theorem \ref{DGthnaive} in a very simple form
(Lemma 4.21 in \Mbook):
\begin{theorem}[\Mb-differential operators as \Mb-half-density kernels]
\label{DGthdiff}
When considered as \Mb-half-densities on $X^2_b$,
the kernels of \Mb-differential operators on $X=\Rplus$ are precisely the
Dirac distributions on $\Delta_b$ which are smooth up to the boundary.
\end{theorem}

\subsubsection{Definition and properties of the small \Mb-calculus}
Theorem \ref{DGthdiff} suggests:
\begin{definition} \label{DGdefsc}
The {\em small \Mb-calculus} $\Psi_b^m(X)$, $m\in\R$, is defined as
the set of (\Mb-half-density-valued) distributions $u$ on $X^2_b$ satisfying
\begin{enumerate}
\item[(a)]
$u$ is conormal of order $m$ with respect to $\Delta_b$, smoothly up to 
the boundary \ff,
\item[(b)]
$u$ vanishes to infinite order at \lb\ and \rb.%
\footnote{
Requiring $u$ to vanish in a neighborhood of \lb\ and
\rb\ would work as well, but would be somewhat less natural in the
context of the full calculus.
}
\end{enumerate}
\end{definition}

We now want to see that the '\PDO\ machine'  from Subsection
\ref{DGsubsecpdo} works for $\Psi^*_b$ as well. 
We address the essential properties listed there, in
varying degrees of completeness.
\begin{boldlist}
\item[(Symb)]
For a \Mb-differential operator as in \eqref{DGeqbop} the principal
symbol is defined as 
$$p(x,y,\lambda,\eta) =\sum_{j+|\alpha|=m} 
       a_{j\alpha}(x,y) \lambda^j\eta^\alpha.$$
We don't give the definition for general \Mb-\PDO s here.
We only remark that, in order to make this defined invariantly, it should
be considered as function on a bundle called ${}^bT^*X$ 
(the \Mb-cotangent bundle) by Melrose. See \Mbook, Sections 2.2 and 4.10.

Ellipticity means that $p(x,y,\lambda,\eta)\not=0$ whenever
$(\lambda,\eta)\not=0$, $x\geq 0$.%
\footnote{\label{DGfnell}%
Schulze's definition (see \cite{DGSchul:PDBVPCSA}) of ellipticity requires,
in addition, that the 'conormal operator' (the operator family on the
boundary given locally by
$\sum_{j,\alpha} a_{j\alpha}(0,y)\lambda^j\partial_y^\alpha$,
$\lambda\in\C$) be invertible on a 'weight line' 
$\Re\lambda=\frac{n}2-\gamma$, where $\gamma$ is a parameter. 
This condition will be imposed also in the parametrix construction in the
full calculus, see Remarks~\ref{DGremfc} below. It is satisfied for all but
countably many values of $\gamma$ and ensures that $P$ is Fredholm between
suitable $\gamma-$weighted Sobolev spaces (if $X$ is compact). 
Lesch shows (\cite{DGLes:OFTCSAM}) that Fredholmness holds even without this
condition, for any closed extension of $P$ (whose domain may then not
be a Sobolev space).
}
\item[(Alg)]
Let us check that $\Psi_b^*$ is closed under composition.%
\footnote{
Assuming, as in \ref{DGsubsubsecdistpdo}, that there is no problem
with integrability at infinity, i.e.\ that at least one factor is properly
supported. Compare Footnote~\ref{DGfnprop} and 
Remark~\ref{DGremclasspdo}.3.
}
(That the symbol map preserves products is then 
done precisely as in the classical case.)
As we discussed in Remark \ref{DGrempbt}.2, 
for a systematic analysis of composition one needs to blow-up
the space $X^3$, and a good way to do this is to
first blow up zero and then the preimages of the
three coordinate axes. The resulting space is called 'triple \Mb-space'
$X^3_b$, see Figure \ref{DGfigtriple}.
Recall that the three projections $\pi_i:X^3\to X^2$,
$i=1,2,3$ lift to \Mb-fibrations $\pi_{ib}:X_b^3\to X_b^2$ (called
$\pitilde_i$ before). The composition formula \eqref{DGeqcomp}, rewritten
in terms of \Mb-spaces, reads:
\begin{equation}\label{DGeqbcomp}
P\circ Q = \pi_{2b*}(\pi_{3b}^*P\cdot \pi_{1b}^*Q).
\end{equation}
We give the argument purely in geometric terms.
Some details are left to the reader.
Let $P,Q$ be distributions on $X^2_b$, conormal with respect to $\Delta_b$
and supported near $\Delta_b$.
We only check condition (a) in Definition \ref{DGdefsc}.
Condition (b) will be discussed shortly in the context of the full calculus
(see the proof of Theorem \ref{DGthfcprop}).
{\small
\begin{enumerate}
\item
The maps $\pi_{ib}$ are transversal to $\Delta_b$. By the pull-back theorem for
conormal distributions (see \ref{DGsubsubsecconmwc})
$\Delta_{ib}:= \pi_{ib}^{-1}(\Delta_b)$ are p-submanifolds
(see Figure \ref{DGfigtriple} for $\Delta_{3b}$), and
$\pi_{3b}^*P,\pi_{1b}^*Q$ are conormal with respect to
$\Delta_{3b},\Delta_{1b}$, respectively.
\item
The behavior of the space $X^3_b$ and the maps $\pi_{ib}$  transversal to
$\Delta_b^{(3)}:=\Delta_{1b}\cap\Delta_{3b}$ (the space diagonal
in $X^3_b$) remains the same at the boundary (i.e.\ at \fff) as in the
interior (cf.\ the discussion in \ref{DGsubsubsecconmwc}), and
$\pi_{2b}(\Delta_b^{(3)}) = \Delta_b$. 
Therefore, the discussion of product and push-forward in 
\ref{DGsubsecdist} carries over literally, as far as the part $v$
in \eqref{DGeqprodsum}
is concerned.
\item
The map $\pi_{2b}$ is a fibration near $\Delta_{1b}$ and $\Delta_{3b}$
and satisfies the transversality condition in the push-forward
Theorem \ref{DGthpfdist}.
Case (a) of that theorem and the arguments after \eqref{DGeqvcon}
 apply, therefore part $w$ in
\eqref{DGeqprodsum} is smooth.
Since $P$ and $Q$ are supported near $\Delta_b$, $P\circ Q$ is supported
near $\Delta_b$.
\item
It remains to check the smoothness of the conormal singularity of
$P\circ Q$ up to the boundary, i.e.\ at \ff. By the Pull-Back Theorem \ref{DGthpbt}
(or rather its trivial extension to the present context), $\pi_{3b}^*P$ and
$\pi_{1b}^*Q$ have index set $0$ at \fff\  (recall that $0$ stands for
smooth behavior, see Remark \ref{DGremasymp}.3)%
\footnote{
Here one needs to check that the (half-)density factors do not introduce
extra powers of bdf's. Compare Footnote \eqref{DGeqbdenstrf}.
}, then so does their
product, and then the Push-Forward Theorem \ref{DGthpft} shows that $P\circ Q$ has index set $0$
at \ff.
\end{enumerate}
}
\item[(Exact)]
Checking the short exact sequence \eqref{DGeqses} is almost
trivial (as in the classical case), once symbols are defined.
\item[(AC)]
Asymptotic completeness is easy, just the same as in the classical case.
\item[(Diff)]
We have $\Diff_b^m\subset \Psi_b^m$ by construction.
\item[(Map)]
See the section on the full calculus.
\item[(Ell)]
Typical elliptic \Mb-operators are the Laplacian and Dirac operators
on Riemannian manifolds with infinite cylindrical ends, or (up to
a conformal factor) with conical points.
\item[(Neg)]
This is the main difference  between the small calculus and the
classical \PDO\ calculus, and the point that makes an extension
(the full calculus) necessary: 
\begin{equation}\label{DGeqnonc}
\text{\em Elements in $\Psi^{-\infty}_b$ are not necessarily compact.}
\end{equation}
This is due to the $1/x$ factor that always occurs in \Mb-densities%
\footnote{
in other words, the $1/x$-factor in Theorem \ref{DGthnaive},
so this is not an artifact of the \Mb-formalism!
}.
More precisely, since $X=\Rplus$ is not compact we should look at
(Neg') instead: Choose $\phi\in C_0^\infty(\Rplus)$ with $\phi(0)=1$.
Then:
\begin{equation} \label{DGeqnoncomp}
 \phi P\in\Psi_b^{-\infty}\text{ is compact }\iff P_{|\ff} \equiv 0.
\end{equation}
This is not hard to check directly, see \Mbook, Section 4.14.
We'll just check the even easier:
$$ \phi P\in\Psi_b^{-\infty}\text{ is Hilbert-Schmidt }\iff P_{|\ff} \equiv 0. $$
\begin{quote}
{\small \proof
Note that the condition on a half-density $\alpha$ to be $L^2$
(square-integrable) is independent of the choice of a measure, since
$|\alpha|^2$ is a density. Similarly, for an operator $A$ acting on
half-densities on a space $X$, being Hilbert-Schmidt is equivalent to
$\int_{X^2}|A|^2<\infty$ for the Schwartz kernel.
Assume for (notational) simplicity that $P$ is supported in a compact set
disjoint from \lb\ and \rb.
When we use projective coordinates $x,s=x'/x$ on $X^2_b$
then $P\in\Psi_b^{-\infty}$ means $P(x,s) = p(x,s)|\dxx\dss|^{1/2}$
with $p$ smooth in $x\in\Rplus,s\in\Rplus$; $\phi P$ has kernel 
$\phi(x)p(x,s)$, which by assumption
is supported in  $x\leq C$, $C^{-1}\leq s\leq C$,
for some positive $C$. Therefore,
\begin{equation} \label{DGeqhs}
\hspace*{2em}
\int_{X^2}|\phi P|^2 = \int_0^C \left( \int_{C^{-1}}^C |\phi(x)p(x,s)|^2\,
\dss\right)\,\dxx,
\end{equation}
and this is finite iff the smooth function in parentheses vanishes at
$x=0$, i.e.\ iff $p(0,s)=0$ for all $s$.
\qed}
\end{quote}
\end{boldlist}

In summary, we have seen that the small calculus has all properties
required for the \PDO\ machine described in Subsection 
\ref{DGsubsecpdo}, except (Neg). Therefore, applying the machine one can
find a parametrix $Q\in\Psi_b^{-m}$ for an elliptic element 
$P\in\Psi_b^{m}$, that is an inverse modulo $\Psi_b^{-\infty}$.
However, since this error term is not compact
(even after localization), this is not enough
to draw many conclusions about the analytic properties of $P$.

Therefore, more work is needed: The full calculus.

%%%%%%%%%%%%%%%%%%%%%%%%%%%%%%%%%%%%%%%%%%%%%%%%%%%%%%%%%%%
%%%%%%%%%%%%%%%%%%%%%%%%%%%%%%%%%%%%%%% Full calc. subsec
%%%%%%%%%%%%%%%%%%%%%%%%%%%%%%%%%%%%%%%%%%%%%%%%%%%%%%%%%%%%

\subsection{The full \Mb-calculus}\label{DGsubsecfc}
Introducing the full \Mb-calculus in detail would exceed the scope of this
article. We refer the reader to \Mbook, chapter 5, for an extensive
treatment. Here we will only state its definition and main properties,
check how it acts on nice function (which gives, again, a nice
illustration of Pull-Back and Push-Forward Theorem), and then outline how this
definition
 arises when one tries to
 improve the parametrix construction for
elliptic \Mb-differential operators.

\begin{definition} \label{DGdeffc}
The {\em full \Mb-calculus} on $X=\Rplus$ is the collection of
spaces $\Psi^{m,\Ecal}_b$ defined as follows:
Let $m\in\R$ and let $\Ecal=(E_\lb,E_\rb)$ be an index family for
$X^2$. Then a distribution $u$ on $X^2_b$ is in $\Psi^{m,\Ecal}_b$
iff $u=u_1+u_2+u_3$ with
\begin{enumerate}
\item[(a)]
$u_1\in\Psi_b^m$, the small calculus,
\item[(b)]
$u_2$ is nice on $X^2_b$, with index family $(E_\lb,0,E_\rb)$
at $(\lb,\ff,\rb)$ ($0$ is the 'smooth' index set, see Remark 
\ref{DGremasymp}.3),
\item[(c)]
$u_3=\beta^*v$, where $\beta:X^2_b\to X^2$ is the blow-down map
and $v$ is nice on $X^2$ with index family $\Ecal$.
\end{enumerate}
\end{definition}
Thus, the conormal singularity at $\Delta_b$ is the same as in
the small calculus, but in addition one allows non-trivial asymptotic
expansions at \lb\ and \rb, plus an additional  'residual' term
which is even nice on $X^2$.%
\footnote{
$u_3$ can not be absorbed into $u_2$ since it has index
family $(E_\lb,E_\lb\cup E_\rb,E_\rb)$ by the Pull-Back Theorem; however, it is
much better than just any function with these index sets.}

The main properties of the full calculus are:
\begin{theorem}\label{DGthfcprop}
The full calculus acts on nice functions and is conditionally
closed under composition,
more precisely:
\begin{description}
\item[Action]
Let $P\in\Psi_b^{m,\Ecal}$ and $w$ be a nice function on $\Rplus$
with index set $F$. If
$$ \inf E_\rb + \inf F > 0$$
then  $Pw$ is nice with index set $E_\lb\extunion F$.
\item[Composition]
Let $P\in\Psi_b^{m,\Ecal}$, $Q\in\Psi_b^{m',\Fcal}$ and
$\Fcal=(F_\lb,F_\rb)$. If
$$ \inf E_\rb + \inf F_\lb > 0$$
then $P\circ Q\in \Psi_b^{m+m',(E_\lb\extunion F_\lb,E_\rb\extunion
F_\rb)}$.
\end{description}
(To avoid problems at infinity, assume $P$ to be properly supported.)
\end{theorem}

\proof
These are Proposition 5.52 and Theorem 5.53 in \Mbook. The proofs
there avoid the systematic use of Pull-Back and Push-Forward Theorem. 
Let us check the statement on action in the systematic way.
Write the kernel of $P$ as $u_1+u_2+u_3$ as in Definition
\ref{DGdeffc}. Let  us assume $P=u_2$ for simplicity, this is
the most interesting part.

Let $\pi_{1/2}:X^2\to X$ be the projections onto the first and second
factor and $\pi_{ib}=\pi_i\circ\beta:X^2_b\to X$ their analogues on 
$X_b^2$. Then, formula \eqref{DGeqaction} becomes
$$ Pw = \pi_{1b*}(P\cdot\pi_{2b}^*w).$$
Using the Pull-Back and Push-Forward Theorem, we can now read off the result
from Figure \ref{DGfigaction}:
\begin{figure}[htbp]
\input{figaction.pstex_t}
\caption{Applying an operator to a function in the full calculus}
\label{DGfigaction}
\end{figure}
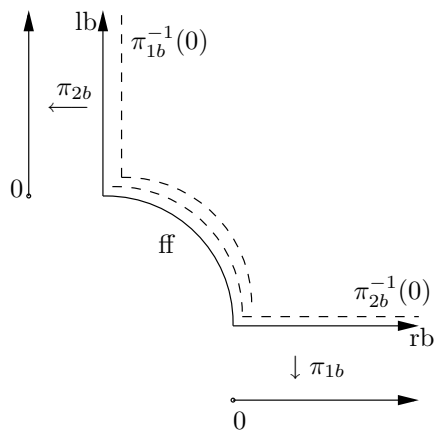

Since $\pi_{2b}^{-1}(0)=\ff\cup\rb$, the Pull-Back Theorem shows that
$\pi_{2b}^*w$ has index family $(0,F,F)$ (at $(\lb,\ff,\rb)$);
therefore, $P\cdot \pi_{2b}^*w$ has index family
$(E_\lb,F,E_\rb+F)$.%
\footnote{
Addition of index sets $E,F$ is defined in the obvious way:
$$ E+F=\{(z+w,k+l):\, (z,k)\in E, (w,l)\in F\}.$$
If $u,v$ are nice with index families $\Ecal,\Fcal$ then clearly $uv$
is nice with index family $\Ecal+\Fcal$.
} 
Finally, since $\pi_{1b}$ is a \Mb-fibration and
$\pi_{1b}^{-1}(0) = \lb\cup\ff$, the Push-Forward Theorem shows that
$\pi_{1b*}(P\cdot\pi_{2b}^*w)$ is nice with index set
$E_\lb\extunion F$, provided the integrability condition
$\inf (E_\rb+F)>0$ holds (since \rb\ is the only face which 
is not mapped to $0$ under $\pi_{1b}$).

The proof for composition proceeds similarly. Here one needs to know
in addition that $\pitilde_{2}:X^3_b\to X^2_b$,
the lift of the 'middle' projection $X^3\to X^2$ already used in 
\eqref{DGeqbcomp}, is a \Mb-fibration. This was checked in Remark
\ref{DGrempbt}.2. \qed

\begin{remarks}[Parametrix, why full calculus, etc.]
\label{DGremfc}
\
\begin{enumerate}
\item
Let us see how terms of type $u_2$ arise from improving the
parametrix 
construction in the small calculus.
Let $P\in\Diff_b$ be elliptic%
\footnote{
In the case $X=\Rplus$ this simply means $a_m(x)\not=0$ for all $x$,
in \eqref{DGeqbop}.
},
and $Q_1$ a parametrix in the small calculus, i.e.\
$PQ_1=\Id+R$ with $R\in\Psi_b^{-\infty}$ (and similarly
for $Q_1P$).
\begin{itemize}
\item
As we saw in \eqref{DGeqnoncomp}, 
the obstruction to compactness of the
remainder $R$ is the restriction of its Schwartz kernel to the front
face \ff.
Therefore, in order to improve the parametrix we must 'cancel'
this obstruction.
\item
For any  $A\in\Psi_b^*$ define the {\em indicial operator}
$I(A)\in\Psi_b^*$ by
\begin{equation}\label{DGeqindic}
\hspace*{2em} A=A(x,s)\dxshalf \implies I(A) = A(0,s)\dxshalf,
\end{equation}
i.e.\ by 'freezing coefficients' at $x=0$.%
\footnote{
For a satisfying discussion of $I(A)$ it now actually matters that its
kernel is not properly supported. Melrose deals with
this by compactifying, i.e.\ adding suitable points at $x=\infty$,
see Section 4.15 in \Mbook; we neglect this here.
}
Since $\ff = \{x=0\}$, we have
$$\hspace*{2em} A_{|\ff} = 0 \iff I(A) = 0.$$
If $P$ is a \Mb-differential operator then \eqref{DGeqdiffkernel} shows that
\begin{equation}\label{DGeqdiffind}
\hspace*{2em} P=\sum_j a_j(x)(x\partial_x)^j\implies I(P)
  =\sum_j a_j(0) (x\partial_x)^j,
\end{equation}
i.e.\ $I(P)$ has constant coefficients as a \Mb-operator.
\item
Since $I(P)$ has constant coefficients it can be inverted
easily:  Substitute $x=e^t$, so that $x\partial_x=\partial_t$, and
then solve a constant coefficient ordinary differential equation.

Alternatively, this may be done using the Mellin transform
(which is just the Fourier transform under this change of variables),
since it transforms $I(P)$ into a multiplication operator.
The latter method works for any elliptic $A\in\Psi_b^*$ instead of $P$.

The calculation shows (see \Mbook, equation (5.28) and Lemma 5.16) 
that the inverse of $I(P)$ (for elliptic $P\in\Diff_b^*$)
has kernel of the type $u_1+u_2$. Two important points are:
\begin{itemize}
\item
For the index sets of $u_2$ one has% 
%(with $-E_\rb:=\{(-z,p):\,(z,p)\in E_\rb\}$):%
\footnote{
Melrose's definition of $\Spec_b$, (5.10) in \Mbook, differs from this
by a 90 degree rotation. Our convention fits better to our definition
of index sets (which is the same as Melrose's in \Mbook; but
in \Mnotes\ index sets are also rotated by 90 degrees; this is more 
consistent with conventions in (mathematical literature on) scattering
theory).
}
\begin{gather} \notag
E_\lb\cup (-E_\rb)=\Spec_b I(P):=\hspace*{.12\textwidth}\\ \label{DGeqSpec}
\hspace*{.2\textwidth} :=\{(z,l):\sum_j a_j(0)z^j 
\text{ has a zero of order}\geq l+1\text{ at }z\}
\end{gather}
(with $-E_\rb:=\{(-z,p):\,(z,p)\in E_\rb\}$).
\item
There is not a unique inverse, but a whole family, $(I(P)^{-1})_\gamma$, 
parametrized by
a real parameter $\gamma$. In terms of kernels, $\gamma$ determines how
$\Spec_b I(P)$ splits up into $E_\lb$ and $E_\rb$. Functional analytically,
$\gamma$ is the weight of a pair of suitable Sobolev spaces on which 
$I(P)$ is actually invertible with inverse $(I(P)^{-1})_\gamma$.
$\gamma$ is restricted to lie in $\R\setminus\{\Re z:\,(z,0)\in\Spec_b
I(P)\}$, and $(I(P)^{-1})_\gamma$ is locally constant on this set
(as a distribution). See \Mbook, Proposition 5.15, and 
Footnote \ref{DGfnell}.%
\footnote{
In Example \ref{DGexbcalc}, we have $I(P)=P$ and
$\Spec_b I(P) = \{(-c,0)\}$. The inverse we constructed there
had $E_\lb=\emptyset, E_\rb=\{(c,0)\}$. It 
corresponds to $\gamma>-\Re c$. As an exercise, check
that it maps $x^\gamma L^2_{\rm comp}(dx/x)\to x^\gamma
L^2_{\rm loc}(dx/x)$. The latter space may be improved to a suitably
defined $H^1$-Sobolev space.
}
\end{itemize}
\item
Finally, a parametrix for $P$ may be constructed as follows:
Set $Q=Q_1+Q_2$, where $Q_2$ is any (compactly supported) operator with
$I(Q_2)=-I(P)^{-1}I(R)$. Then one checks easily that
$PQ=\Id+R'$ with $I(R')=0$.%
\footnote{
Here one needs that $I$ preserves products. This is clear from
\eqref{DGeqdiffind} for \Mb-differential operators, but requires a 
little work in general.

In fact, this is the reason why we identify $P_{|\ff}$ with the
operator $I(P)$: This is done precisely in order to make $P\mapsto
P_{|\ff}$ an algebra homomorphism (the 'second symbol map', see
below)!
}
Since $R$ is smoothing, $Q_2$ is of the type $u_2$. $R'$ is in
$\Psi_b^{-\infty,\Ecal_\gamma}$ and vanishes on \ff,
so it is compact by a similar argument as before (after localization).
More precisely, one has a parametrix for each admissible parameter $\gamma$,
and $R'$ is compact on the Sobolev spaces with weight $\gamma$.
\end{itemize}
\item 
This parametrix construction may be formalized 
analogous to the classical \PDO\ calculus in \ref{DGsubsecpdo}:
Regard the set of constant coefficient operators as second
symbol space (just like for the usual symbols, their inversion is
easier than that of $P$ itself) and $I$ as second symbol map. The central
fact is again a short exact  sequence:
$$ 0 \to\rho\Psi_b^{m,\Ecal}\hookrightarrow\Psi_b^{m,\Ecal}
\stackrel{I}{\to} \Psi_{b,\text{const coeff}}^{m,\Ecal} \to 0,$$
where now $\rho$ is a boundary defining function for \ff,
and $m\in\R$ and the index family $\Ecal$ are (almost) arbitrary
(see (5.160) in \Mbook). 
The parametrix is obtained by  combined  use of both symbol maps
and both short exact sequences.
(The first symbol is defined from the singularity of $u_1$ at $\Delta_b$ as in 
the small calculus.)
\item
The parametrix constructed above does not contain a term of
type $u_3$. But $u_3$ is contained in the calculus since it arises
when composing two terms of type $u_2$. In particular, such
compositions are necessary when improving the parametrix further
(for example, making $R$ vanish to higher than first order at \ff, see the
proof of Theorem \ref{DGthparam}).  
The index set $\Ecal$ will have to be enlarged in this process.
Such a more precise parametrix is constructed in \Mbook, Sections
5.18-5.25.
\end{enumerate}
\end{remarks}

\subsection{General manifolds with boundary} \label{DGsubsecgenmwb}
We describe shortly the changes needed in small and full calculus
when considering a general manifold with boundary $X$ instead
of $\Rplus$. For simplicity, we assume $\partial X$ to be connected.
$x,y$ will denote local coordinates as in 
\eqref{DGeqbop}:
\begin{itemize}
\item
{\em Definition of $X^2_b$ and $X^3_b$: }
In $X^2$, one has boundary defining functions $x$ (for the boundary
of the first factor) and $x'$ (for the second).
The 'corner' in $X^2$ is the submanifold of codimension two
$(\partial X)^2 = \{x=x'=0\}$. Then
$$ X^2_b := [X^2,(\partial X)^2].$$
If $y,y'$ are local coordinates in the boundary of the first and second factor,
then local coordinates on $X^2_b$ are $x$, $s=x'/x$, $y$, $y'$ (and
$x'$, $x/x'$, $y$, $y'$). In other words,
everything is as before, only with $y,y'$ as parameters.
The \Mb-diagonal is defined by the first equation in \eqref{DGeqbdiag},
in coordinates
$$\Delta_b = \{s=1,\ y=y'\}.$$
$X^2_b$ again has three boundary hypersurfaces, denoted
$\lb$, $\rb$, $\ff$ as before and locally given by
$\{s=\infty\}$, $\{s=0\}$, $\{x=0\}$, respectively.

The triple \Mb-space $X^3_b$ is defined by first blowing up
$(\partial X)^3$ in $X^3$ and then the (now disjoint) preimages of
$(\partial X)^2\times X$, $\partial X\times X\times \partial X$,
$X\times (\partial X)^2$. Again, this means doing the same as for
$\Rplus$, carrying the $y$-variables along as parameters.
\item
{\em Small calculus:}
Theorem \ref{DGthdiff} and Definition \ref{DGdefsc} extend
literally, and also the discussion of properties (except that in 
\eqref{DGeqhs} an additional $dy\,dy'$ integration is needed).
\item
{\em Full calculus:}
Definition \ref{DGdeffc} and Theorem \ref{DGthfcprop} (and its
proof) and the first and last point of Remark \ref{DGremfc}.1
(the parametrix construction) extend literally, as well as
Remarks \ref{DGremfc}.2 and 3.
$I(A)$ is now defined on $\Rplus\times\partial X$,
and defined as in \eqref{DGeqindic}, except that $A$ depends on
$y,y'$ also (in local coordinates on the boundary)
and the half-density factor is $\left|\frac{dx}x\frac{ds}sdydy'
\right|^{1/2}$. Similarly, in \eqref{DGeqdiffind} the $a_j$ depend on
$y$ also and in $I(P)$ are replaced by $a_j(0,y)$.
But now $I(P)$ has constant coefficients only in $x$, not in $y$!%
\footnote{
This is the 'partial freezing of coefficients' mentioned in the
Introduction.
}
Therefore, it should be considered as ordinary differential operator
(in $x$) whose coefficients are partial differential operators 
on $\partial X$: 
$$I(P)=\sum A_j\cdot (x\partial_x)^j,\quad A_j\in\Diff (\partial X).$$
 Then the inversion using Mellin transform works
as before.\\ 
$\Spec_b$ is defined as in \eqref{DGeqSpec}, with
$a_j(0)$ replaced by $A_j$. Since  $I(P,z):=\sum_j A_j z^j$ 
is now an operator%
\footnote{
$I(P,z)$ is called 'conormal symbol' by some authors, e.g.\ Schulze
\cite{DGSchul:PDBVPCSA}, see Footnote \ref{DGfnell}.
}
on $\partial X$ for every $z$, the notion of 'zero'  must  be
interpreted suitably: as a point $z$ where $I(P,z)$ is not invertible.
('Order' may also be defined easily, see Section 5.2 in \Mbook.)
The role of $\Spec_b I(P)$ in determining the asymptotic type
of a parametrix at $\lb$ and $\rb$ is as before (except that,
for coordinate invariance, $E_\lb,E_\rb$ have to be 'completed',
cf.\ Footnote \ref{DGfncoordind}).
The only {\em essentially new features} are:
\begin{itemize}
\item
$\Spec_b I(P)$ may be an infinite set (but -- in case $\partial X$
is compact -- it is still discrete
and finite for $\Re z$ bounded, which is proved by
'analytic Fredholm theory'). 
\item
$\Spec_b I(P)$ is global on the boundary, i.e.\ determined
by the (global) solvability of some partial differential equation
on $\partial X$.
\item
The algebra of symbols is not commutative.
(But this does not matter in the parametrix construction since commutativity
was never used.)
\end{itemize}
\end{itemize}

%% file: figaction.pstex_t
\begin{picture}(0,0)%
\epsfig{file=figaction.pstex}%
\end{picture}%
\setlength{\unitlength}{3947sp}%
\begingroup\makeatletter\ifx\SetFigFont\undefined%
\gdef\SetFigFont#1#2#3#4#5{%
  \reset@font\fontsize{#1}{#2pt}%
  \fontfamily{#3}\fontseries{#4}\fontshape{#5}%
  \selectfont}%
\fi\endgroup%
\begin{picture}(2579,2643)(3451,-4319)
\put(5201,-3969){\makebox(0,0)[lb]{\smash{\SetFigFont{10}{12.0}{\rmdefault}{\mddefault}{\updefault}$\downarrow \pi_{1b}$}}}
\put(4209,-1928){\makebox(0,0)[lb]{\smash{\SetFigFont{10}{12.0}{\rmdefault}{\mddefault}{\updefault}$\pi_{1b}^{-1}(0)$}}}
\put(5959,-3794){\makebox(0,0)[lb]{\smash{\SetFigFont{10}{12.0}{\rmdefault}{\mddefault}{\updefault}$\rb$}}}
\put(4384,-3211){\makebox(0,0)[lb]{\smash{\SetFigFont{10}{12.0}{\rmdefault}{\mddefault}{\updefault}$\ff$}}}
\put(4851,-4319){\makebox(0,0)[lb]{\smash{\SetFigFont{10}{12.0}{\rmdefault}{\mddefault}{\updefault}$0$}}}
\put(3451,-2861){\makebox(0,0)[lb]{\smash{\SetFigFont{10}{12.0}{\rmdefault}{\mddefault}{\updefault}$0$}}}
\put(5609,-3502){\makebox(0,0)[lb]{\smash{\SetFigFont{10}{12.0}{\rmdefault}{\mddefault}{\updefault}$\pi_{2b}^{-1}(0)$}}}
\put(3684,-2336){\makebox(0,0)[lb]{\smash{\SetFigFont{10}{12.0}{\rmdefault}{\mddefault}{\updefault}$\longleftarrow$}}}
\put(3743,-2219){\makebox(0,0)[lb]{\smash{\SetFigFont{10}{12.0}{\rmdefault}{\mddefault}{\updefault}$\pi_{2b}$}}}
\put(3859,-1811){\makebox(0,0)[lb]{\smash{\SetFigFont{10}{12.0}{\rmdefault}{\mddefault}{\updefault}$\lb$}}}
\end{picture}

%% file: DGappendix.tex
\section*{Appendix: Pull-back, push-forward, densities etc.}
\label{DGsecpfpb}
\setcounter{theorem}{0}
\setcounter{section}{1}
\renewcommand{\thesection}{\Alph{section}}

Let $f:M\to N$ be a smooth map between manifolds.

The {\em pull-back} by $f$ of a function $v$ on $N$ is the function
$$f^*v = v\circ f$$
on $M$. Clearly, $f^*v$ is smooth if $v$ is. 

Pull-backs appear everywhere. Depending on context and personal taste, 
they may be interpreted
as 'plugging in', 'reinterpretation', or 'distortion'. 
For example, $v(xy)$
(plugging in $xy$ into $v$) is the pull-back of $v$ under the map
$f(x,y)=xy$; $\pi_2^*v$ from Example \ref{DGexpfpb}.1 is just $v$
reinterpreted as function on $\R^2$; and if $f$ is a diffeomorphism,
then $f^*v$ is just $v$ looked at through the 'distortion lens' $f$.
(For example, if $f:(0,\infty)\to(0,\infty),x\to x^2$ then 
the graph of $f^*v$ is obtained from the graph of $v$ by a stretching
for $x<1$ and a compression for $x>1$.)

Related, though quite different at first glance, is {\em push-forward}
by $f$. The idea is that, for a function $u$ on $M$, $(f_*u)(y)$ for 
$y\in N$ should be 'the integral of $u$ over the fiber $f^{-1}(y)$'.
Now this clearly depends on the choice of a measure%
\footnote{
We are a little sloppy about the use of the word 'measure': Contrary
to standard usage we only require that a measure be defined on 
bounded Borel sets (i.e.\ those contained in a compact set). 
Thus, $u(x)dx$ is a measure on $\R$ for any {\em locally} integrable
function $u$.} on this fiber (e.g.\ $dy$ in \eqref{DGeqpf1}).
Rather than to consider $u$ and this measure separately, or to
consider a measure on each fiber, it is more convenient to start
with a measure $\mu$ (Borel, complex) on all of $M$.
For measures, push-forward is a standard operation:
$f_*\mu$ is the measure on $N$ defined by 
$(f_*\mu)(V)=\mu(f^{-1}(V))$ (= measure of the union of all fibers over $V$), 
$V\subset N$.%
\footnote{
With our use of the word 'measure' one needs to require 'integrability' here.
This is guaranteed for example when $f$ is proper on the support of $\mu$.
We always assume this tacitly.}
In terms of integrals, this is equivalent to
\begin{equation}\label{DGeqpf}
\int_N (f_*\mu)\,\phi = \int_M \mu\, f^*\phi  \tag{App.1}
\end{equation}
for all $\phi\in C_0^\infty$.
If $f=\pi_1$ as in Example \ref{DGexpfpb}.1 then this easily gives
$$ \pi_{1*}(u(x,y)\,dxdy) = \left(\int_\R u(x,y)\,dy\right)\,dx, $$
which shows  that push-forward in the sense of measure theory is
integration over the fiber, as we intended.
 The additional
factor $dx$ may look cumbersome, but this is the only way to have invariance
with respect to coordinate changes on both $M$ and $N$.

\eqref{DGeqpf} shows that push-forward $f_*$ is dual to pull-back $f^*$,
under the duality of functions and measures. Also, \eqref{DGeqpf} is
not just formal nonsense but actually a recipe for calculation:

\begin{example}
Let $f:(0,\infty)^2\to (0,\infty)$ be such that $f(x,y)=xy$, and let
$\mu(x,y)=u(x,y)\,dxdy$. Then $(f^*\phi)(x,y)=\phi(xy)$, so
\begin{align*}
\int_{(0,\infty)^2} u f^*\phi &= \int_0^\infty \int_0^\infty u(x,y)\phi(xy)\,dy\,dx\\
&= \int_0^\infty \int_0^\infty u(x,\frac{t}x)\phi(t)\,\frac{dt}x\,dx,\\
\end{align*}
(using Fubini and changing variables $t=xy$ in the inner integral)
and comparison with \eqref{DGeqpf} gives%
\footnote{
Anyone who is still sceptical of densities should once try to 
calculate (or just make sense of) the notion of integrating a function over
the hyperbola $xy=t$!}
$$ f_*u(t) = \left( \frac1x \int_0^\infty u(x,t/x)\,dt \right)\,dx.$$
\end{example}

In the smooth context we usually deal not with arbitrary measures,
but rather with the more special densities ('smooth measures') and 
with the more general distributional densities, which we introduce next.

A {\em  (smooth)  density} $\mu$ on a manifold $M$ is a measure such that
for any local coordinate system
$x:U\subset M\to\R^n$ one can find a smooth function $u$ on $U$ such that
$\mu(U') = \int_{U'} u(x)dx$ for all measurable $U'\subset U$. 
In this case we write $\mu=u(x)\,dx$ for short.%
\footnote{
This definition shows how $u$ transforms under a change of
coordinates. Of course one can use this to define a line bundle over $M$,
usually denoted $\Omega_M$,
such that densities are just sections of this line bundle, see \Mbook,
Section 4.5. }

A {\em distribution density} on $M$ 
is an object which on every coordinate patch
$x:U\subset M\to\R^n$ looks like $u(x)dx$ for a distribution $u(x)$
on $U$, where $u$ transforms as for densities.%
\footnote{
If you want to define this more formally, you can exploit the 
idea of duality; then a distributional density on $M$ is simply an element 
in the dual space of $C_0^\infty(M)$.
}

The push-forward of a distribution density $\mu$ under the map $f$
may be defined by Equation \eqref{DGeqpf} again (where $\int$ is interpreted 
as the usual pairing of distributions and functions), which shows that 
$f_*\mu$ is a distribution density again.

So far, everything was quite straight-forward, the main problem was keeping
the dualities straight. Here comes a more substantial point:
\begin{center}
{\em The push-forward of a smooth density need not be a smooth density!}
\end{center}
For example, if $f:\R\to\R, x\mapsto x^3$ then 
$f_*dx = \frac13 y^{-2/3}\,dy$.%
\footnote{
Proof: $\int dx (f^*\phi)(x)=\int \phi(x^3)dx = \int
\phi(y) y^{-2/3}dy/3$, and use \eqref{DGeqpf}.
}

However, if $f$ is a {\em fibration} then $f_*\mu$ is a smooth density
whenever $\mu$ is (assuming integrability).%
\footnote{\label{DGfnfibr}%
Recall the definition of a fibration:
For any $x$ the preimage $L=f^{-1}(x)$ is a manifold and a neighborhood
of $L\subset M$ can be identified (via a diffeomorphism) with $U\times L$,
for some neighborhood $U$ of $x$, such that
$f(x',l)=x'$ for all $l\in L$, $x'\in U$. I.e., locally $f$ looks like (a higher-dimensional
version of) $\pi_1$ in
Example \ref{DGexpfpb}.1, and then the assertion is clear.
}

A {\em distribution} is a continuous functional on the set of
(compactly supported) smooth densities. Therefore, the pull-back
$f^*u$ is defined for a distribution $u$ if $f$ is a fibration
(as the adjoint operation to $f_*$).%
\footnote{
If $f$ is not a fibration, then $f^*u$ may not make sense at all.
For example, expressions like $\delta(x^3)$ or $\delta(xy)$ make no sense.
}
Actually, here it suffices that $f$ be a {\em submersion}, i.e.\
have surjective differential at every point. (This is weaker than, and
the local analogue of, $f$ being a fibration.)
$f^*u$ may then be defined by approximation of $u$ by smooth functions
(see \cite{DGHor:ALPDOI}, chapter VI.1).
\newpage

Table \ref{DGtabapp} gives an overview of push-forward and pull-back compatibilities.

\begin{table}[htb]
\setlength{\extrarowheight}{2pt}
\begin{tabular}{l||l|l|l}
& $\alpha$ smooth  & $\alpha$ distributional & 
$\alpha$ conormal distributional \\ \hline\hline
Push-forward $f_*\alpha$ & $f$ fibration & any $f$ & $f$ fibration,\\
&&& transversality condition \\ \hline
Pull-back $f^*\alpha$ & any $f$ & $f$ submersion & $f(M)$ transversal\\
&&& to $\singsupp \alpha$ \\ \hline
\end{tabular}
\vspace*{1em}

\caption{Conditions on a smooth map $f:M\to N$ to define push-forward
or pull-back of a density (resp.\ function) $\alpha$, resulting in an object 
of same type as $\alpha$. For push-forward, it is always assumed that $f$ is
proper on $\supp \alpha$.}\label{DGtabapp} 
\end{table}

Finally, when talking about operators it is useful to have (smooth
or distributional) {\em half-densities}. These are defined as objects which
look like
$$ u(x)|dx_1\cdots dx_n|^{1/2} $$
in any coordinate system $x$, with $u$ smooth respectively a distribution.
Accordingly, $u$ transforms with the square root of the Jacobian under
coordinate changes. The reader who prefers more formal definitions may
consult \Mbook, Section 4.5 for the definition of a (trivial) line bundle
of which half-densities are sections.

%% file: DGart.bbl
\begin{thebibliography}{10}

\bibitem{DGBieMil:CDCZBUMSLI}
E.~Bierstone and P.~Milman.
\newblock Canonical desingularization in characteristic zero by blowing up the
  maximum strata of a local invariant.
\newblock {\em Inv. Math.}, 128:207--302, 1997.

\bibitem{DGBruSee:RSA}
J.~Br{\"u}ning and R.~Seeley.
\newblock Regular singular asymptotics.
\newblock {\em Adv. in Math.}, 58:133--148, 1985.

\bibitem{DGBruSee:RESORSO}
J.~Br{\"u}ning and R.~Seeley.
\newblock The resolvent expansion for second order regular singular operators.
\newblock {\em J. Func. Anal.}, 73:369--429, 1987.

\bibitem{DGEpsMelMen:RLSPD}
C.~Epstein, R.~Melrose, and G.~Mendoza.
\newblock {Resolvent of the Laplacian on strictly pseudoconvex domains}.
\newblock {\em {Acta Math.}}, 167:1--106, 1991.

\bibitem{DGGil:FAEHTNSAECO}
J.~B. Gil.
\newblock {\em Full asymptotic expansion of the heat trace for
non-self-adjoint elliptic cone operators}.
\newblock Preprint 2000.

\bibitem{DGGriGru:SALPFT}
D.~Grieser and M.~J. Gruber.
\newblock Singular asymptotics lemma and push-forward theorem.
\newblock In J.~Gil, D.~Grieser, and M.~Lesch, editors, {\em Approaches to
  Singular Analysis}, Advances in Partial Differential Equations, Basel, 2000.
  Birkh\"auser.

\bibitem{DGHasMazMel:ASAE}
A.~Hassell, R.~Mazzeo, and R.~B. Melrose.
\newblock {Analytic surgery and the accumulation of eigenvalues}.
\newblock {\em Commun. Anal. Geom.}, 3(1):115--222, 1995.

\bibitem{DGHasVas:RLTOACS}
A.~Hassell and A.~Vasy.
\newblock The resolvent for {L}aplace-type operators on asymptotically conic
  spaces.
\newblock Preprint, lanl-server math/0002114, 2000.

\bibitem{DGHir:RSAVFCZ}
H.~Hironaka.
\newblock Resolution of singularities of an algebraic variety over a field of
  characteristic zero.
\newblock {\em Ann. Math.}, 79:109--326, 1964.

\bibitem{DGHor:ALPDOI}
L.~{H\"ormander}.
\newblock {\em The {Analysis} of {Linear} {Partial} {Differential} {Operators}
  {I}}.
\newblock Springer-Verlag, 1983.

\bibitem{DGHor:ALPDOIII}
L.~{H\"ormander}.
\newblock {\em The {Analysis} of {Linear} {Partial} {Differential} {Operators}
  {III}}.
\newblock Springer-Verlag, 1985.

\bibitem{DGLauSei:PDAMWBCBCA}
R.~Lauter and J.~Seiler.
\newblock Pseudodifferential analysis on manifolds with boundary --- a
  comparison of b-calculus and cone algebra.
\newblock In J.~Gil, D.~Grieser, and M.~Lesch, editors, {\em Approaches to
  Singular Analysis}, Advances in Partial Differential Equations, Basel, 2000.
  Birkh\"auser.

\bibitem{DGLes:OFTCSAM}
M.~Lesch.
\newblock {\em Operators of Fuchs type, conical singularities, and asymptotic
  methods}.
\newblock Number 136 in Teubner Texte zur Mathematik. Teubner--Verlag, Leipzig,
  1997.

\bibitem{DGLoy:PDCMWC}
P.~A. Loya.
\newblock {\em On the b-pseudodifferential calculus on manifolds with corners}.
\newblock PhD thesis, MIT, 1998.

\bibitem{DGLoy:SREPO}
P.~A. Loya.
\newblock {\em The structure of the resolvent of elliptic pseudodifferential
  operators}.
\newblock Preprint 2000.

\bibitem{DGMaz:HCCCM}
R.~Mazzeo.
\newblock The {H}odge cohomology of a conformally compact metric.
\newblock {\em J. Diff. Geom.}, 28:309--339, 1988.

\bibitem{DGMaz:ETDEOI}
R.~Mazzeo.
\newblock {Elliptic theory of differential edge operators. I}.
\newblock {\em Commun. Partial Differ. Equations}, 16(10):1615--1664, 1991.

\bibitem{DGMazMel:ALHCLSSF}
R.~Mazzeo and R.~B. Melrose.
\newblock {The adiabatic limit, Hodge cohomology and Leray's spectral sequence
  for a fibration}.
\newblock {\em J. Differ. Geom.}, 31(1):185--213, 1990.

\bibitem{DGMazMel:ASEI}
R.~Mazzeo and R.~B. Melrose.
\newblock {Analytic surgery and the eta invariant}.
\newblock {\em Geom. Funct. Anal.}, 5(1):14--75, 1995.

\bibitem{DGMazMel:PDOMWFB}
R.~Mazzeo and R.~B. Melrose.
\newblock Pseudodifferential operators on manifolds with fibred boundaries.
\newblock {\em Asian J. Math.}, 2:833--866, 1998.

\bibitem{DGMcDon:LSWCS}
P.~McDonald.
\newblock {\em The Laplacian on spaces with conical singularities}.
\newblock PhD thesis, MIT, 1990.

\bibitem{DGMel:DAMWC}
R.~B. Melrose.
\newblock Differential analysis on manifolds with corners.
\newblock Book in preparation. http://www-math.mit.edu/$\sim$rbm/book.html.

\bibitem{DGMel:TBP}
R.~B. Melrose.
\newblock {Transformation of boundary problems}.
\newblock {\em Acta Math.}, 147:149--236, 1981.

\bibitem{DGMel:PDOCSL}
R.~B. Melrose.
\newblock {Pseudodifferential operators, corners and singular limits}.
\newblock In {\em Proc. Int. Congr. Math., Kyoto/Japan 1990, Vol. I}, pages
  217--234, 1991.

\bibitem{DGMel:CCDMWC}
R.~B. Melrose.
\newblock {Calculus of conormal distributions on manifolds with corners}.
\newblock {\em Int. Math. Res. Not.}, 3:51--61, 1992.

\bibitem{DGMel:APSIT}
R.~B. Melrose.
\newblock {\em {The Atiyah-Patodi-Singer index theorem}}.
\newblock A. K. Peters, Ltd., Boston, Mass., 1993.

\bibitem{DGMel:SSTLAES}
R.~B. Melrose.
\newblock {Spectral and scattering theory for the Laplacian on asymptotically
  Euclidean spaces}.
\newblock In {Ikawa, Mitsuru}, editor, {\em Spectral and scattering theory.
  Proc. Taniguchi international workshop, Sanda, Hyogo, Japan}, volume 161 of
  {\em Lect. Notes Pure Appl. Math.}, pages 85--130, Basel, 1994. Marcel
  Dekker.

\bibitem{DGMel:GST}
R.~B. Melrose.
\newblock {\em {Geometric scattering theory}}.
\newblock Cambridge Univ. Press., 1995.

\bibitem{DGMel:FCAPDO}
R.~B. Melrose.
\newblock {Fibrations, compactifications and algebras of pseudodifferential
  operators}.
\newblock In {Hoermander, Lars} et~al., editors, {\em Partial differential
  equations and mathematical physics. The Danish-Swedish analysis seminar,
  Copenhagen, Denmark, Lund, Sweden, March 17-19, May 19-21, 1995.
  Proceedings.}, volume~21, pages 246--261, Boston, MA, 1996. Birkh\"auser.

\bibitem{DGMel:GOBS}
R.~B. Melrose.
\newblock {Geometric optics and the bottom of the spectrum}.
\newblock In {Colombini, Ferruccio} et~al., editors, {\em {Geometrical optics
  and related topics}}, volume~32, pages 267--293, {Boston, MA}, 1997.
  {Birkh\"auser}.

\bibitem{DGMelMen:EOTCT}
R.~B. Melrose and G.~A. Mendoza.
\newblock Elliptic operators of totally characteristic type.
\newblock MSRI Preprint, 1983.

\bibitem{DGMelNis:HPDOIMWB}
R.~B. Melrose and V.~Nistor.
\newblock Homology of pseudo-differential operators {I}. {M}anifolds with
  boundary.
\newblock Preprint, 1996.

\bibitem{DGMelPia:AKTMWC}
R.~B. Melrose and P.~Piazza.
\newblock Analytic $k$-theory on manifolds with corners.
\newblock {\em Adv. Math.}, 92:1--26, 1992.

\bibitem{DGMoo:HKAMWCS}
E.~A. Mooers.
\newblock {Heat kernel asymptotics on manifolds with conic singularities.}
\newblock {\em J. Anal. Math.}, 78:1--36, 1999.

\bibitem{DGSchul:TIOSSS}
B.-W. Schulze.
\newblock {\em Topologies and invertibility in operator spaces with symbolic
  structures}, volume 111 of {\em Teubner {T}exte zur {M}athematik}, pages
  259--288.
\newblock Teubner Verlag, Leipzig, 1989.

\bibitem{DGSchul:PDOMWS}
B.-W. Schulze.
\newblock {\em Pseudo-Differential Operators on Manifolds with Singularities},
  volume~24 of {\em Stud. Math. Appl.}
\newblock North-Holland, Amsterdam, 1991.

\bibitem{DGSchul:PDBVPCSA}
B.-W. Schulze.
\newblock {\em Pseudo-Differential Boundary Value Problems, Conical
  Singularities and Asymptotics}.
\newblock Akademie Verlag, Berlin, 1994.

\bibitem{DGShu:PDOST}
M.~A. Shubin.
\newblock {\em Pseudodifferential operators and spectral theory}.
\newblock Springer-Verlag, 1987.

\bibitem{DGSim:PDO}
S.~R. Simanca.
\newblock {\em Pseudo-differential operators}.
\newblock Number 236 in Pitman Res. Notes in Math. Longman, 1990.

\end{thebibliography}
